\documentclass[a4paper,11pt,twoside]{article}
\author{
	\normalsize Fran\c{c}ois Delarue \\[8pt]
	\small Laboratoire J.A.Dieudonn\'e \\
	\small UMR CNRS-UNS No 7351 \\
	\small Universit\'e de Nice Sophia-Antipolis \\
	\small Parc Valrose \\
	\small France-06108 NICE Cedex 2\\
	\small delarue@unice.fr
	\and
	\normalsize William Salkeld \\[8pt]
	\small Brown University, 182 Geogre street \\
	\small Providence, RI 02906 \\
	\small  william\_salkeld@brown.edu
}
\usepackage[nodayofweek]{datetime}

\usepackage[utf8]{inputenc}
\usepackage[T1]{fontenc}
\usepackage[english]{babel}

\usepackage{charter}

\usepackage{xcolor}
\usepackage{verbatim}

\usepackage{hyperref}

\usepackage[top=3.cm, bottom=4.0cm, left=2.2cm, right=2.2cm]{geometry}
\usepackage{amsmath}
\usepackage{amsthm}	
\usepackage{amsfonts}	
\usepackage{amssymb,bbm}
\usepackage{shuffle}
\usepackage{mathrsfs}
\usepackage{pifont}

\usepackage{enumerate}
\usepackage{enumitem}
\usepackage[nobysame, alphabetic,initials]{amsrefs}

\usepackage{stmaryrd}
\usepackage{appendix}
\usepackage{etoolbox}		

\usepackage{graphicx}

\usepackage{
	ulem		
} 

\numberwithin{equation}{section}

\theoremstyle{plain}
\newtheorem{theorem}{Theorem}[section]
\newtheorem{lemma}[theorem]{Lemma}
\newtheorem{proposition}[theorem]{Proposition}

\newtheorem{definition}[theorem]{Definition}

\newtheorem{remark}[theorem]{Remark}
\newtheorem{example}[theorem]{Example}

\allowdisplaybreaks

\newcommand{\bE}{\mathbb{E}}
\newcommand{\bF}{\mathbb{F}}

\newcommand{\bN}{\mathbb{N}}
\newcommand{\bP}{\mathbb{P}}

\newcommand{\bR}{\mathbb{R}}

\newcommand{\bW}{\mathbb{W}}

\newcommand{\cA}{\mathcal{A}}
\newcommand{\cB}{\mathcal{B}}
\newcommand{\cC}{\mathcal{C}}

\newcommand{\cE}{\mathcal{E}}
\newcommand{\cF}{\mathcal{F}}
\newcommand{\cG}{\mathcal{G}}
\newcommand{\cH}{\mathcal{H}}
\newcommand{\cI}{\mathcal{I}}
\newcommand{\cJ}{\mathcal{J}}

\newcommand{\cL}{\mathcal{L}}
\newcommand{\cM}{\mathcal{M}}

\newcommand{\cP}{\mathcal{P}}

\newcommand{\cR}{\mathcal{R}}
\newcommand{\cS}{\mathcal{S}}

\newcommand{\fD}{\mathfrak{D}}

\newcommand{\fF}{\mathfrak{F}}

\newcommand{\fL}{\mathfrak{L}}

\newcommand{\fR}{\mathfrak{R}}

\newcommand{\fT}{\mathfrak{T}}

\newcommand{\fh}{\mathfrak{h}}

\newcommand{\fr}{\mathfrak{r}}

\newcommand{\scA}{\mathscr{A}}
\newcommand{\scC}{\mathscr{C}}

\newcommand{\scE}{\mathscr{E}}
\newcommand{\scF}{\mathscr{F}}
\newcommand{\scG}{\mathscr{G}}
\newcommand{\scH}{\mathscr{H}}
\newcommand{\scI}{\mathscr{I}}
\newcommand{\scJ}{\mathscr{J}}

\newcommand{\scL}{\mathscr{L}}
\newcommand{\scM}{\mathscr{M}}
\newcommand{\scN}{\mathscr{N}}
\newcommand{\scP}{\mathscr{P}}
\newcommand{\scQ}{\mathscr{Q}}
\newcommand{\scS}{\mathscr{S}}
\newcommand{\scT}{\mathscr{T}}

\newcommand{\rD}{\mathbf{D}}

\newcommand{\rp}{\mathbf{p}}
\newcommand{\rQ}{\mathbf{q}}

\newcommand{\rv}{\mathbf{V}}
\newcommand{\rw}{\mathbf{W}}

\newcommand{\rId}{\mathbf{1}}

\newcommand{\vertiii}{{\vert\kern-0.25ex \vert\kern-0.25ex \vert}}

\DeclareMathOperator{\supp}{supp}
\DeclareMathOperator{\Shuf}{Shuf}

\DeclareMathOperator{\lin}{Lin}
\DeclareMathOperator{\lion}{Lions}

\newcommand{\A}[1]{A_{#1}[0]}

\usepackage{pgf,tikz}
\usetikzlibrary{arrows}

\tikzstyle{vertex} = [fill, shape=circle,inner sep=2pt,]
\tikzstyle{edge} = [fill, line width = 0.5pt]
\tikzstyle{zhyedge1} = [opacity=.5,fill opacity=.5, line cap=round, line join=round, line width=27pt,color=black]
\tikzstyle{zhyedge2} = [opacity=.5,fill opacity=.5, line cap=round, line join=round, line width=25pt,color=white]
\tikzstyle{hyedge1} = [opacity=.5,fill opacity=.5, line cap=round, line join=round, line width=27pt]
\tikzstyle{hyedge2} = [opacity=.5,fill opacity=.5, line cap=round, line join=round, line width=25pt, color=white]

\tikzstyle{vertexS} = [fill, shape=circle,inner sep=1pt,]
\tikzstyle{edgeS} = [fill, line width = 0.25pt]
\tikzstyle{zhyedge1S} = [opacity=.5,fill opacity=.5, line cap=round, line join=round, line width=12pt,color=black]
\tikzstyle{zhyedge2S} = [opacity=.5,fill opacity=.5, line cap=round, line join=round, line width=10pt,color=white]
\tikzstyle{hyedge1S} = [opacity=.5,fill opacity=.5, line cap=round, line join=round, line width=12pt]
\tikzstyle{hyedge2S} = [opacity=.5,fill opacity=.5, line cap=round, line join=round, line width=10pt, color=white]

\pgfdeclarelayer{background}
\pgfsetlayers{background,main}

\title{An example driven introduction \\ to probabilistic rough paths}

\begin{document}
	
	\maketitle
	
	\begin{abstract} 
		In these notes, we provide an introduction to a new regularity structure used for solving rough mean-field equations. 
		
		The index set of this regularity structure is described a collection of novel objects which we refer to as Lions trees. These objects arise in Taylor expansions involving the Lions derivative and capture many of the desirable properties of mean-field dynamics. 
		
		This work represents a comprehensive generalisation of the ideas first introduced in \cite{2019arXiv180205882.2B} that promise powerful insights into how interactions with a collective determine the dynamics of an individual within this collective. 
	\end{abstract} 
	
	
	
	\vspace{0.3cm}
	
	
	Secondary: 60L30, 16T05, 46G05
	
	
	\noindent
	{\bf Acknowledgements}: William Salkeld wishes to thank the \textbf{London Mathematical Society} for the award of an \emph{Early Career Fellowship} (ECF-1920-29) which facilitated this research.
	
	Further, William Salkeld was supported by MATH+ project AA4-2 and by the US Office of Naval Research under the Vannevar Bush Faculty Fellowship N0014-21-1-2887.  

	\setcounter{tocdepth}{2}
	\tableofcontents
	
	\section{Introduction}
	
	This is one of a series of papers studying the infinitesimal dynamics of mean-field equations. The breadth of material developed necessarily requires a number of manuscripts but the focus of this paper has been to overview these ideas, providing readers with an accessible introduction. 
	
	As such, many fine details and generalisations will be skipped over and we hope that the reader will delay any scepticism for another paper in which these will be addressed thoroughly. 
	
	\subsection{Motivation}
	This work is part of a wider program of research at the intersection between two active fields in mathematics. On the one hand, rough path theory, as first developed in \cite{lyons1998differential}, has proved to be a robust tool in the study of stochastic differential systems, particularly those that are driven by highly oscillatory signals and are thus out of the scope of standard integration theories such as Young integration for regular enough paths or stochastic integration for semimartingales. On the other hand, mean-field theory, which was originally dedicated to the analysis of large weakly interacting particle systems in statistical physics and in fluid mechanics (\cites{kac1956foundations,McKean1966}), but has grown recently due to a surge of interest in line with the development of a calculus of variations on the space of probability measures and the study of related optimisation problems from transport theory, mean-field control or mean-field games, see for instance \cite{villani2008optimal} for the former and \cites{LionsVideo,CarmonaDelarue2017book1, CarmonaDelarue2017book2} for the latter two. 
	
	In their most frequent form, as usually described in stochastic analysis (see for instance \cites{Sznitman, meleard1996asymptotic, jourdain1998propagation, kurtz1999particle}), mean-field dynamical models are addressed through a convenient type of distribution dependent stochastic equations in which the pointwise realisation of the unknown interacts with its own statistical time marginal distribution. This specific kind of self-interaction is commonly referred to as a McKean-Vlasov interaction. It captures the asymptotic form of the interactions in the underlying large particle system as the number of particles is taken to infinity. Noticeably, it says that a single tagged particle only sees the others through a deterministic aggregate which is mathematically described by a probability measure. Not only does it provide a robust formulation of the limiting particle system, it also offers a less complex and hence more tractable representation which precisely explains why the mean-field paradigm has shown to be an efficient and scalable solution in the analysis of large population problems. The reduction of complexity is made possible by a statistical averaging phenomenon which is usually called \textit{propagation of chaos} in the existing literature. In brief, when all particles are subjected to independent and identically distributed (and hence chaotic) inputs, chaos is transferred from the inputs to the outputs (as given by the states of the particles) in the limit when the number of those particles becomes infinite. Remarkably, this does not require the particles to be independent when there are only a finite number of them. In other words, independence just manifests asymptotically. Most of the existing theory consider particles subjected to semimartingales drivers, typically Brownian motions for continuous dynamics and possibly Lévy processes for dynamics featuring jumps. 
	
	The adaptation of the rough path theory to the mean-field setting is a very natural and meaningful question which has already stimulated a series of works in the recent years. The analysis of such \textit{rough mean-field models} goes back to the article \cite{CassLyonsEvolving}. However, in the latter work and in most of the following ones on the subject, see for instance \cites{Bailleul2015Flows, deuschel2017enhanced}, the interaction between the particles appears solely in the drift and not in the volatility. Whilst this may appear to be an innocent simplification, this restriction reveals a substantial difficulty that is intrinsic to the mean-field setting. In order to explain this challenge properly, it is worth recalling that rough path theory relies on a relevant form of local expansions (see \cite{gubinelli2004controlling} for a systematic point of view). The solution of a rough differential system driven by the infinitesimal variation of a possibly rough signal is sought within a class of paths that behave locally like the driving signal. This allows one to expand the coefficients of the equation in terms of some local iterated integrals of the signal by means of a suitable version of Taylor's formula. In particular, when the intensity of the driving signal is a nonlinear function of the unknown, it gives a robust sense to the underlying integration provided all the involved iterated integrals are known exogeneously. Importantly, there is no canonical definition of those iterated integrals, at least up to some order $p$, when the underlying signal is below the standard threshold for Young's integration, the value of $p$ being determined by the regularity of the signal (through the H\"older exponent if the path is H\"older continuous). This requires one to choose those iterated integrals that are relevant to the model at hand (for instance, It\^o and Stratonovich integrals lead to two different forms of iterated integrals for the Brownian motion and, depending on the context, one may be more suitable than the other). The collection of iterated integrals forms the \emph{Signature}, which embeds an enriched description of the signal into a single enhanced path.
	
	Thus, a major obstacle to a rough mean-field setting is the provision of a convenient expansion of the coefficients with respect to the distributional argument. The strategy to do so was clarified in a recent contribution by the first author, see \cite{2019arXiv180205882.2B}. At its simplest, the idea is to consider a Taylor expansion of the coefficients using the advances achieved over more than twenty years in the construction of a differential calculus on the space of probability measures, the latter being also referred to as the Wasserstein space. These ideas originate in the work of \cite{Jordan1998variation} on the connection between Fokker-Planck equations (which are partial differential equations for the time marginal law of a McKean-Vlasov equation) and gradient flows on the Wasserstein space. We refer the reader to the monograph \cite{Ambrosio2008Gradient} for a complete overview of the subject. In \cite{2019arXiv180205882.2B}, differential calculus on the Wasserstein space is implemented along Lions' approach, as originally introduced in \cite{LionsVideo} (see also \cite{CarmonaDelarue2017book1}*{Chapter 5}). Although the resulting two notions of gradient coincide, Lions' point of view is very convenient for probabilistic (or Lagrangian) approaches to mean-field equations and models. The main idea is to lift a differentiable function defined on the Wasserstein space to a function on the Hilbert space of square-integrable random variables and take a Fr\'echet derivative. The lifting operation consists of associating with any probability measure $\mu$ a random variable $X$ on a chosen probability space $(\Omega, \cF, \bP)$ which has exactly $\mu = \bP \circ (X)^{-1}$, $X$ being called a representative of $\mu$. Although there may be many choices for $X$, there are in practice canonical representatives for the probability measures under study. Since Lions' derivative is computed as a Fr\'echet derivative on the space of random variables, these canonical representatives appear explicitly in all the related derivatives. This makes it possible to locally expand the distributional dependence in the coefficients of a McKean-Vlasov equation. 
	
	\subsection{Contribution}
	
	The central contribution of our work (not only this paper but also the forthcoming ones) is the exposition of how the theory of regularity structures (applied to rough paths) can be extended to incorporate a Taylor expansion using Lions derivatives and how these microscopic distributional contributions determine the macroscopic dynamics of distribution dependent equations. Further, we will demonstrate how this theory is consistent irregardless of whether it is used to understand systems of interacting equations or their mean-field limit. Motivated by existing knowledge from rough path theory, we would expect the additional necessary information about a driving signal to take the form of iterated integrals. 
	
	To be consistent with the aforementioned standard results on propagation of chaos for weakly interacting particle systems, one would expect the collection of iterated integrals that appear in a rough mean-field setting to correspond to the collection of iterated integrals that appear in a corresponding particle system, where the latter is understood in the rough sense and the number of particles therein tends to infinity. In fact, although this picture is indeed sound, it is not sufficiently accurate. The first observation in this regard is that the mean-field limit may not retain all the iterated integrals from the particle system and that accordingly some of the latter ones may disappear asymptotically because of combinatorial reasons. To rephrase, several iterated integrals from the particle system may contribute to the distributions of iterated integral in the asymptotic regime. 
	In line with this, the second observation is that 
	because of the symmetries that are inherent to the particle system, a given iterated integral can be repeated several times just by modifying, by means of a mere permutation of the particles, the indices that decorate its labelled trees. However, it should be clear that those repetitions do not carry any valuable additional information. Indeed, the Butcher-Connes-Kreimer algebra that is associated to a large but finite system of $n\times d$-dimensional time dependent particles 
	\begin{equation}
		\label{eq:part:sys}
		\big( W_t^{1, N}, ..., W_t^{N, N} \big)_{t\in [0,1]}
	\end{equation}
	naturally leads to consider trees that are indexed by pairs comprising both a label $i \in \{1, ..., N\}$ for the number of the particle and a label $j \in \{1, ..., d\}$ for the component of this particle. There are $d \times N$ labels, as expected if the whole system is regarded as a single equation. \textit{Stricto sensu} and in full generality, an iterated integral that involves particles $i_1, ..., i_k \in \{1, ..., N\}$ differs from the iterated integral, obtained by using the same order of integration but by replacing $i_1$ by $\sigma(i_1)$, $i_2$ by $\sigma(i_2)$, and so forth, for $\sigma$ a permutation on $\{1, ..., N\}$. However, it is very likely that those iterated integrals should be identified, in some manner, in the mean-field setting, because of the intrinsic statistical exchangeability of the system, say for instance if $(W_t^{1, N})_{t\in [0,1]}, ..., (W_t^{N, N})_{t\in [0,1]}$ are independent and identically distributed. 
	
	In the end, these preliminary remarks strongly suggest that the collection of trees that appear in the mean-field limit cannot (and should not) be the mere enumeration of all the possible labelled trees that would appear if the iterated integrals were just constructed on larger and larger particle systems. Somehow, the resulting collection would be much too wide and would miss some of the deep properties of symmetry of a mean-field system. 
	
	\vspace{5mm}
	Our first contribution in our series of works is to identify a smaller class of trees that are well-fitted to the mean-field paradigm. Critically, these trees do not carry any explicit labels for the particles. Similar to a McKean-Vlasov equation, which accounts for the state of one typical tagged particle within an infinite population, the trees that we consider below index the iterated integrals on the top of typical tagged small groups of particles (hence regardless of the precise labels carried by those particles). This construction is made possibly by retaining the sole clusters formed by the labels of the particles in a standard labelled tree: in brief, two nodes in the same labelled tree are in the same cluster if they are labelled by the same particle index. Once the labels of the particles have been erased, we are left with trees equipped with \textit{hyperedges}. This idea is at the core of our research.  
	
	In particular, it is one of our main result to prove that only a certain type of system of hyperedges is admissible. Of course, the fact that some other systems of hyperedges must be ruled out is consistent with our preliminary remark on the need to forget some of the labelled trees when passing from a finite to an infinite population. We refer to trees equipped with these admissible hyperedges as \textit{Lions trees}. Our approach is very advantageous as it fits effectively with the symmetric framework of the problem. However, the main drawback is the fracturing of the standard Butcher-Connes-Kreimer co-algebraic structure. We provide in Section \ref{section:ProbabilisticRoughPaths} below a comprehensive introduction to such Lions trees. 
	
	In order to illustrate one of the main difficulties, the reader must understand that the situation is more difficult in the probabilistic rough path setting than in the usual McKean-Vlasov setting, although this theory shares a similar background. Propagation of chaos says that particles become independent in the limit, which suffices when the noises are Brownian motions. When iterated integrals are introduced, the story becomes somewhat different: iterated integrals are associated with tagged groups of particles and when several of these integrals have to be manipulated at the same time, it is no longer clear which particle is which. Statistically, the latter makes a strong difference: Groups equipped with implicitly distinct particles must be independent (in some manner), whilst groups sharing a common particle must be correlated. As a result, we require a richer algebraic operation that clarifies the possible ``correlations'' between the hyperedges of the two trees. This problem is the same as the \textit{coupling} problem in probability: given two probability measures on two spaces, there may be plenty of ways to reconstruct a probability measure on the product space with the two original probability measures as marginal measures. This gives rise to the \emph{coupled Hopf algebra}, which we also present in Section \ref{section:ProbabilisticRoughPaths}. 
	
	Once a convenient graded algebraic structure satisfying a commutativity relationship similar to those of Hopf algebras has been introduced, probabilistic rough paths can be defined as paths on the associated characters, see Section \ref{section:Models} for the main principles. This turns out to be equivalent to the definition provided in \cite{2018arXiv180205882B}, \cite{2019arXiv180205882.2B} although our approach is far more general and provides greater insight into the underlying algebraic and analytic relationships that arise from the probabilistic framework. However, we additionally observe that the group of characters is typically too large to capture all the necessary properties and we instead focus on subgroups that capture whether the associated equation describes either a particle system or a continuum. Although not immediately obvious, this is actually very natural: the characters of the model space of a regularity structure describe the structure group since the character properties describe those of a Taylor expansion. Should that Taylor expansion have additional properties, we would expect these to describe a subgroup of the characters that contains our structure group. 
	
	\subsubsection*{Organization of the paper}
	
	We stress again that the paper has been designed as an introduction to the theory of \textit{probabilistic rough paths}. As such, we have decided to focus below on the main aspects and, if needed, to skip some of the proofs and just to refer to some of our other works on the subject for the most technical details. Also, it must be clear to the reader that the scope of the article just covers the mean-field analogue of a regularity structure, or equivalently of a rough path and its signature. In particular, there is no rough McKean-Vlasov equation handled in this paper. This would require other developments of the theory, which will be eventually addressed in forthcoming works, a very sketchy presentation of which is given in the next paragraph. Whilst this limitation to the basics of a probabilistic rough path, without any further analysis of the equations driven by such a rough path, may be rather frustrating for the reader, the material that is presented here is in fact sufficiently demanding to deserve a paper on its own. To rephrase, our objective in these notes is to provide a theory that accommodates both the high-dimensional rough signal that drives a particle system and the infinite-dimensional structure that is expected to drive in the end the corresponding rough mean-field equation. In this regard, we feel very important to stress that, although it looks very formal, our construction has empirical roots. It is indeed one of our main concern below to have a theory that is consistent with the standard rough path paradigm when the two approaches are used to interpret the same particle system. 
	
	Section \ref{section:ProbabilisticRoughPaths} is devoted to the introduction of the algebraic structure supporting probabilistic rough paths. This is certainly the main ingredient in our theory and we must admit that some effort is needed to capture all the ingredients. We hope that the examples that are given next will be useful to the reader. In particular, we spend some time   introducing the definition of a Lions tree (see Definition \ref{definition:Forests}), which is the basic block in our construction. As we already explained, a Lions tree is a special kind of decorated tree, whose structure is connected with the higher-dimensional Taylor formula for smooth functionals defined on the space of probability measures. For this reason, we provide in the next Subsection \ref{section:TaylorExpansions} a summary of some important results about the latter. The next step is to associate with such trees a module that replaces the aforementioned Butcher-Connes-Kreimer algrebra. The very feature of it is that elements of this module expand as sums of random variables constructed on tensorised probability spaces, with each tensorisation being dictated by the structure of the underlying Lions tree. We refer to Subsection \ref{subse:module:lions:forests}. The definition a coupled tensor product, which is a key operation on elements of the module, is explained in Subsection \ref{subsec:CoupledBialgebra}. Some work is then needed to obtain a coupled bialgebra and then to equip it with a grading in Subsection \ref{subse:2.5} and a group structure in Subsection \ref{subse:2.6}. Those two are very important ingredients for the next Section \ref{section:Models}, in which we introduce the definition of a probabilistic rough path as a character, see Definition
	\ref{definition:General-PRP}. The reader will find in Subsection \ref{subsection:Examples} two very natural examples, which will be key in our next works to understand the connection between a rough particle system (with weak interaction) and the corresponding rough mean-field equation. 
	
	\subsubsection*{Previous work}
	The first ideas about probabilistic rough paths can be found in the ArXiv preprint \cite{2018arXiv180205882B} which was subsequently published under a different title \cite{2019arXiv180205882.2B}. 
	
	This paper is based on many of the key ideas of the previous unpublished work \cite{2021Probabilistic}, although it also includes some ideas present in the preprint \cite{salkeld2021Probabilistic2}. However, many ideas have been reformulated and restructured in order to improve the accessibility of this body of research. 
	
	This work also regularly references the companion papers \cite{salkeld2022Lions} which explores Lions-Taylor expansions in much greater detail and \cite{salkeld2022LionsTrees} which covers many of the results from \cite{2021Probabilistic} which are not included in this work, some updated material that finds its roots in \cite{salkeld2021Probabilistic2} as well as several new and important results. 
	
	Compared to the latter ones, the focus of this paper is answering the question ``what is a probabilistic rough path and what do I need to understand in order to work with them?''. 
	
	\subsection{Lions-Taylor expansions and coupled structures}
	\label{section:TaylorExpansions}
	
	We start by considering how we should represent a Taylor expansion for a function of the form
	\begin{equation*}
		f: \bR^e \times \cP_2(\bR^e) \to \bR^d. 
	\end{equation*}
	
	To explain this, we provide a reminder about derivatives constructed on the so-called ``Wasserstein space'' of probability measures with finite second moment. The latter space is a particular case of the following definition. For any $p\geq 1$, let $\cP_p(\bR^e)$ be the set of all probability measures on $\big( \bR^e, \cB(\bR^e) \big)$. This is equipped with the $\bW^{(p)}$-Wasserstein distance defined by:
	\begin{equation}
		\label{eq:WassersteinDistance}
		\bW^{(p)}(\mu,\nu) = \inf_{\Pi \in \cP_{p}(\bR^e \times \bR^e)}
		\biggl( \int_{\bR^e \times \bR^e}
		|x-y|^p
		d \Pi(x,y) \biggr)^{1/p},
	\end{equation}
	the infimum being taken with respect to all the probability measures $\Pi$ on the product space $\bR^e \times \bR^e$ with $\mu$ and $\nu$ as respective $e$-dimensional marginal laws.
	
	\subsubsection*{The Lions derivative}
	\label{subsection:1-Lip_sup_envelope}
	
	For a function $f:\cP_2(\bR^e) \to \bR^d$, we consider the canonical lift $F: L^2(\Omega, \cF, \bP; \bR^e) \to \bR^d$ defined by $F(X) = f( \bP\circ X^{-1})$. We say that $f$ is $L$-differentiable at $\mu$ if $F$ is Fr\'echet differentiable at some point $X$ such that $\mu = \bP\circ X^{-1}$. Denoting the Fr\'echet derivative by $DF$, it is now well known (see for instance \cite{GangboDifferentiability2019} that $DF$ is a $\sigma(X)$-measurable random variable of the form $DF(\mu, \cdot):\bR^e \to \lin(\bR^e, \bR^d)$ depending on the law of $X$ and satisfying $DF(\mu, \cdot) \in L^2\big( \bR^e, \cB(\bR^e), \mu; \lin(\bR^e, \bR^d) \big)$. We denote the $L$-derivative of $f$ at $\mu$ by the mapping $\partial_\mu f(\mu)(\cdot): \bR^e \ni x \to \partial_\mu f(\mu, x) \in \lin(\bR^e, \bR^d)$ satisfying $DF(\mu, X) = \partial_\mu f(\mu, X)$. This derivative is known to coincide with the so-called Wasserstein derivative, as defined in for instance \cite{Ambrosio2008Gradient}, \cite{CarmonaDelarue2017book1} and \cite{GangboDifferentiability2019}. As we explained in the introduction, Lions' approach is well-fitted to probabilistic approaches for mean-field models since, very frequently, we have  a ``canonical'' random variable $X$ for representing the law of a given probability measure $\mu$. 
	
	Just as rough path theory makes an intense use of Taylor expansions of smooth functions, our program relies on a form of Taylor expansion based on the derivative $\partial_\mu$. Of course, this requires first extending Lions-Taylor expansions to higher orders, which requires some care. To wit, the second order derivatives are obtained by differentiating $\partial_{\mu} f$ with respect to $x$ (in the standard Euclidean sense) and $\mu$ (in the same Lions' sense). The two derivatives $\nabla_{x} \partial_{\mu} f$ and $\partial_{\mu} \partial_{\mu} f$ are thus very different functions: The first one is defined on ${\cP}_{2}({\bR}^e) \times {\bR}^e$ and writes $(\mu,x) \mapsto \nabla_{x} \partial_{\mu} f(\mu,x)$ whilst the second one is defined on ${\cP}_{2}({\bR}^e) \times {\bR}^e \times {\bR}^e$ and writes $(\mu,x,x') \mapsto \partial_{\mu} \partial_{\mu} f(\mu,x,x')$. The $e$-dimensional entries of $\nabla_{x} \partial_{\mu} f$ and $\partial_{\mu} \partial_{\mu} f$ are called here the \textit{free} variables, since they are integrated with respect to the measure $\mu$ itself. 
	
	In practice, we will be working with measure functionals that are dependent on a multi-variable $(x_0, \mu) \in \bR^e \times \cP_2(\bR^e)$ so we need to consider a Taylor expansion involving both Lions and spacial derivatives. Unlike with the free variables generated by the application of iterative Lions derivatives where the $m^{th}$ free variable $x_m$ can only appear if the function has been differentiated at least $m$ times with respect to $\mu$, the variable $x_0$ can appear in any derivatives of $f$ whenever $f$ is a function of the form $f(x_0, \mu)$. 
	
	This leads us to Definition \ref{def:a} below, the principle of which can be stated as follows for the first and second order Lions derivatives: The derivative symbol $\partial_{\mu}$ can be denoted by $\partial_{1}$ and then the two derivative symbols $\nabla_{x_1} \partial_{\mu}$ and $\partial_{\mu} \partial_{\mu}$ can be respectively denoted by $\partial_{(1,1)}$ and $\partial_{(1,2)}$. Derivatives with respect to the $x_0$-component are encoded through the inclusion of a ‘0’. Repeated 0’s thus account for repeated derivatives in the direction of $x_0$.
	
	\begin{definition}
		\label{def:a}
		The sup-envelope of an integer-valued sequence $(a_{i})_{i=1,...,n}$ of length $n$ is the non-decreasing sequence $(\max_{i=1,...,k} a_{i})_{k=1,...,n}$. The sup-envelope is said to be $1$-Lipschitz (or just $1$-Lip) if, for any $k \in \{2,...,n\}$, 
		$$
		\max_{i=1,...,k} a_{i} \leq 1+ \max_{i=1,...,k-1} a_{i}.
		$$
		We call $A_{n}$ the collection of all ${\bN}$-valued sequences of length $n$, with $a_{1}=1$ as initial value and with a 1-Lip sup-envelope. Thus $A_n$ is the collection of all sequences $(a_k)_{k=1, ..., n} \in A_n$ taking values on $\{1, ..., n\}$ such that 
		$$
		a_1 = 1, \quad a_k \in \Big\{1, ..., 1+ \max_{i=1, ..., k-1} a_i \Big\}. 
		$$
		We refer to $A_n$ as the collection of \emph{partition sequences}. 
		
		Let $k, n\in \bN_0$ and let $A_{k,n}[0]$ be the collection of all sequences 
		$a' = (a_i')_{i=1, ..., k+n} = \sigma( (0) \cdot a)$ where $(0) = (0)_{i=1, ..., k}$ is the sequence of length $k$ with all entries $0$, $a\in A_n$ and $\sigma$ is a $(k,n)$-shuffle, i.e., a permutation of $\{1, ..., n+k\}$ such that $\sigma(1) < ... < \sigma(k)$ and $\sigma(k+1)<... < \sigma(n+k)$. 
		
		For a given $n \in \bN$, we let 
		$$
		A_{n}[0]= \bigcup_{k=0}^n A_{k,n-k}[0], \quad A^n[0] = \bigcup_{i=0}^n A_i[0] \quad \mbox{and}\quad A[0] = \bigcup_{n\in \bN_0} \A{n}.
		$$
		Given $a\in \A{n}$, we denote 
		$$
		|a| = n \quad \mbox{and} \quad m[a] = \max_{i=1, ..., n} a_i.
		$$ 
	\end{definition}
	
	Within this framework, we can define the derivative $\partial_a$ for some $a \in A_n[0]$ with $n\in \bN$. By induction, we let:
	\begin{align*}
		\partial_{(0)} =& \nabla_{x_0}, \quad \partial_{(1)} = \partial_\mu \quad \mbox{and}
		\\
		\partial_{(a_1, ..., a_{k-1}, a_k)} =&
		\begin{cases} 
			\nabla_{x_0} \cdot \partial_{(a_1, ..., a_{k-1})} & \quad a_k =0, 
			\\
			\nabla_{x_{a_k}} \cdot \partial_{(a_1, ..., a_{k-1})} & \quad 0 < a_k \leq \max \{ a_1, ..., a_{k-1}\}, 
			\\
			\partial_\mu \cdot \partial_{(a_1, ..., a_{k-1})} & \quad a_k > \max \{a_1, ..., a_{k-1}\}.  
		\end{cases}
	\end{align*}
	\textbf{Further notations.} In line with Definition \ref{def:a} we introduce the following notations, but the reader can skip this paragraph over the first reading. With in the same framework as above, we denote $\llbracket a\rrbracket_c$ to be the equivalence class of all sequences such that
	\begin{align*}
		(b_i)_{i=1, ..., n} \in \llbracket a\rrbracket_c \quad \iff \quad &\Big\{ b^{-1} \Big\} = \Big\{ a^{-1}[j]: j=0, ..., m[a] \Big\} \in \scP\big( \{1, ..., n\} \big)
		\\ 
		& \mbox{and}\quad b^{-1}[c] = a^{-1}[0]. 
	\end{align*}
	For $a, a'\in A_n[0]$, we say that $a \subseteq a'$ if and only if
	\begin{equation}
		\label{eq:partial-ordering}
		\begin{aligned}
			&a^{-1}[0] \subseteq (a')^{-1}[0] \quad \mbox{and}
			\\
			&\forall j=1, ..., m[a]\quad \exists j' \in \{0, 1, ..., m[a']\} \quad \mbox{such that} \quad  a^{-1}[j] \subseteq (a')^{-1}[j'].  
		\end{aligned}
	\end{equation}
	For a sequence $\boldsymbol{b} = (b_i)_{i=1, ..., n}$ and $a\in A_n[0]$ such that $a \subseteq \llbracket b\rrbracket_c$ we define
	\begin{equation}
		\label{eq:K-sequence2}
		\boldsymbol{b}\circ(a) = \Big( b_{a^{-1}[j]} \Big)_{j=1, ..., m[a]}. 
	\end{equation}	
	Let $\alpha, \beta>0$. We define $\scG_{\alpha, \beta}: A[0] \to \bR^{+}$ by
	$$
	\scG_{\alpha, \beta}[a] = \alpha \cdot \Big| a^{-1}\big[ 0 \big] \Big| + \beta \cdot \Big| \bigcup_{j=1}^{m[a]} a^{-1}\big[ j\big] \Big|. 
	$$
	Then for $\gamma>\alpha \wedge \beta$, we define
	$$
	A^{\gamma, \alpha, \beta}[0]:= \Big\{ a\in A[0]: \scG_{\alpha, \beta}[a] \leq \gamma \Big\}. 
	$$
	
	\subsubsection*{A Lions-Taylor expansion}
	
	We now provide our Lions-Taylor expansion. Let $k, n \in \bN_0$ and let $a\in A_{k, n}[0]$. Let $x_0, y_0\in \bR^e$ and let $\Pi^{\mu, \nu}$ be a measure on $(\bR^e)^{\oplus 2}$ with marginal distribution $\mu, \nu \in \cP_{n+1}(\bR^e)$ (in other words, a coupling of $\mu$ and $\nu$). Then, for $f$ a function defined on $\bR^e \times \cP_2(\bR^e)$ (with values in some Euclidean space), we define the operator
	\begin{align}
		\nonumber
		\rD^a& f(x_0, \mu)[ y_0-x_0, \Pi^{\mu, \nu}]
		\\
		\label{eq:rDa}
		&= \underbrace{\int_{(\bR^e)^{\oplus 2}} ... \int_{(\bR^e)^{\oplus 2}} }_{\times m[a]} \partial_a f \Big( x_0, \mu, \boldsymbol{x}_{m\{a\}} \Big) \cdot \bigotimes_{i=1}^{|a|} ( y_{a_i} - x_{a_i})  \cdot d\big( \Pi^{\mu, \nu}\big)^{\times m[a]} \Big( (\boldsymbol{x}, \boldsymbol{y})_{m\{a\}} \Big). 
	\end{align}
	Here, for compact notation we have denoted
	\begin{align*}
		&d\big( \Pi^{\mu, \nu}\big)^{\times m[a]} \Big( (\boldsymbol{x}, \boldsymbol{y})_{m\{a\}}\Big) = d\Pi^{\mu, \nu}(x_1, y_1) \times ... \times d\Pi^{\mu, \nu}(x_{m[a]}, y_{m[a]} ) \quad \mbox{and}
		\\
		& \boldsymbol{x}_{m\{a\}} = (x_1, ..., x_{m[a]}). 
	\end{align*}

	The following statement makes a connection between the operator $\rD^a$ and (standard) derivatives in Euclidean spaces. This connection generalizes the results reviewed in \cite{CarmonaDelarue2017book1}*{Chapter 5} and provides a dictionary that will enable us to pass (typically) from expansions for particle systems to expansions for mean-field equations.
	
	\begin{proposition}
		\label{proposition:classicTay<=>LionsTay*}
		Let $n, N\in \bN$. Let $f: \bR^e \times \cP_2(\bR^e) \to \bR^d$ and suppose that $f\in C_b\big[A^{n}[0]\big]\big( \bR^e \times \cP_2(\bR^e); \bR^d\big)$. Let $x^{1, N}, ..., x^{N, N} \in \bR^e$ and denote $\boldsymbol{x} = (x^{1, N}, ..., x^{N, N})$. 
		
		We define the empirical projection $\boldsymbol{f}: (\bR^e)^{\oplus N} \to (\bR^d)^{\oplus N}$ of the function $f$ by 
		\begin{equation}
			\label{eq:bf-field}
			\boldsymbol{f}\big( \boldsymbol{x} \big) = \boldsymbol{f}\big( (x^{1, N}, ..., x^{N, N}) \big) = \bigoplus_{i=1}^N \overline{f}_i \big( \boldsymbol{x} \big)
		\end{equation}
		where $\overline{f}_i: (\bR^e)^{\oplus N} \to \bR^d$ is the projection pinned at particle $i$ defined by
		\begin{equation*}
			\overline{f}_i \big( \boldsymbol{x} \big) = \overline{f}_i \big( (x^{1, N}, ..., x^{N, N}) \big) = f\Big( x^{i, N}, \bar{\mu}_N\big[ \boldsymbol{x}^N \big] \Big), \quad \bar{\mu}_N\big[ \boldsymbol{x} \big] = \sum_{j=1}^N \delta_{x^{j, N}}. 
		\end{equation*}
		Let $\boldsymbol{i} = (i_1, ..., i_n)$ where each $i_j \in \{1, ..., e\}$. Then 
		\begin{equation}
			\label{eq:example:Sums-finer-partitions}
			\nabla_{\boldsymbol{i}} \overline{f}_i \big( \boldsymbol{x} \big) = \sum_{\substack{a\in A_n[0] \\ a\subseteq \llbracket \boldsymbol{i} \rrbracket_i }} \frac{1}{N^{\times m[a]}} \cdot \partial_a f\Big( x^{i, N}, \bar{\mu}_N\big[\boldsymbol{x}\big] , \boldsymbol{x}^{\boldsymbol{i}\circ (a), N} \Big)  
		\end{equation}
		where we used the notation 
		\begin{equation*}
			\boldsymbol{x}^{\boldsymbol{i} \circ a, N} = \Big( x^{i_{a^{-1}[1]}, N}, ..., x^{i_{a^{-1}[m[a]]}, N} \Big). 
		\end{equation*}
		
		Additionally, if for $\boldsymbol{x}, \boldsymbol{y} \in (\bR^e)^{\oplus N}$ we define $\Pi^{\boldsymbol{x}, \boldsymbol{y}} \in \cP_2(\bR^e \oplus \bR^e)$ as the trivial coupling between the empirical measures associated with $\boldsymbol{x}$ and $\boldsymbol{y}$, namely 
		\begin{equation}
			\label{eq:proposition:classicTay<=>LionsTay}
			\Pi^{\boldsymbol{x}, \boldsymbol{y}} = \frac{1}{N} \sum_{j=1}^N \delta_{(x^{j, N}, y^{j, N})}. 
		\end{equation}
		then finite-dimensional Taylor's formula yields
		\begin{equation}
			\label{eq:proposition:classicTay<=>LionsTay*}
			\boldsymbol{f}\Big( \boldsymbol{y} \Big) = \sum_{a\in A^{n}[0]} \frac{1}{|a|!} \bigoplus_{i=1}^N \rD^a f \Big( x^{i, N}, \bar{\mu}_N\big[ \boldsymbol{x} \big] \Big)\Big[ y^{i, N} - x^{i. N}, \Pi^{\boldsymbol{x}, \boldsymbol{y}} \Big] + \boldsymbol{R}_n^{\boldsymbol{x}, \boldsymbol{y}}\big( \boldsymbol{f} \big)
		\end{equation}
		where the remainder term satisfies
		\begin{equation*}
			\boldsymbol{R}_n^{\boldsymbol{x}, \boldsymbol{y}}\big( \boldsymbol{f} \big) \leq O\Big(  |\boldsymbol{y} - \boldsymbol{x}|^{n+1} \Big). 
		\end{equation*}
	\end{proposition}
	
	A mean-field version of Proposition \ref{proposition:classicTay<=>LionsTay*} is the purpose of the following Theorem:	
	\begin{theorem}
		\label{theorem:LionsTaylor2}
		Let $f: \bR^e \times \cP_2(\bR^e) \to \bR^d$ and suppose that $f$ has Lions and spacial derivatives of all orders. Then for $n \in \bN$, any $\mu, \nu \in \cP_{n+1} (\bR^e)$ with joint distribution $\Pi^{\mu, \nu}$ and any $x_0, y_0\in \bR^e$, we have that
		\begin{align}
			\label{eq:FullTaylorExpansion}
			f(y_0, \nu) 
			= \sum_{a\in A[0]}^{\gamma, \alpha, \beta} \frac{1}{|a|!} \cdot \rD^a f\Big( x_0, \mu \Big)\Big[ y_0 - x_0, \Pi^{\mu, \nu} \Big]  
			+ R_{\gamma, \alpha, \beta}^{(x_0, y_0), \Pi^{\mu, \nu}} \big( f \big). 
		\end{align}
		where the remainder term is of order $n+1$. 
	\end{theorem}
	
	For a full proof along with an explicit statement for $R_{\gamma, \alpha, \beta}^{(x_0, y_0), \Pi^{\mu, \nu}} \big( f \big)$ along with some analysis of its asymptotic nature as $\mu \to \nu$ and $x_0 \to y_0$ (which is beyond the scope of this work), we direct the reader to \cite{salkeld2022Lions}. This result is included to give context to the ideas that we will be exploring in great detail in the coming paper. 
	
	Another natural observation is that the Fr\'echet derivatives are symmetric multi-linear operators. This leads us to the following
	
	\begin{theorem}[Schwarz Theorem for Lions Derivatives]
		\label{thm:Schwarz-Lions}
		Let $a[0]$ and let $f: \bR^e \times \cP_2(\bR^e) \to \bR^d$ such that $\partial_a f$ exists. Let $\Shuf(|a|)$ be the set of permutations on the set $\{1, ..., |a| \}$. 
		
		Then $\forall \sigma\in \Shuf(|a|)$, $v_1, ..., v_{|a|} \in \bR^e$, $\mu\in \cP_2(\bR^e)$, $x_0, x_1, ..., x_{|a|} \in \bR^e$, 
		\begin{equation}
			\label{eq:proposition:Schwartz-2}
			\partial_{a} f\Big( x_0, \mu, \boldsymbol{x}_{m\{a\}} \Big) \cdot \bigotimes_{i=1}^{|a|} v_i 
			=
			\partial_{\llbracket \sigma[a]\rrbracket_0} f\Big( x_0, \mu, \boldsymbol{x}_{\llbracket \sigma(a)\rrbracket \circ (a)} \Big) \cdot \bigotimes_{i=1}^{|a|} v_{\sigma(i)} . 
		\end{equation}
		where
		\begin{equation*}
			\boldsymbol{x}_{m\{a\}} = \Big( x_1, ..., x_{m[a]} \Big)
			\quad \mbox{and}\quad 
			\boldsymbol{x}_{\llbracket \sigma(a)\rrbracket \circ (a)} = \Big( x_{\sigma(a)_{a^{-1}[1]}}, ..., x_{\sigma(a)_{a^{-1}[m[a]]}} \Big). 
		\end{equation*}
	\end{theorem}
	
	\begin{remark}
		Restricting ourselves to the case $|a| = 2$, we remark that Equation \eqref{eq:proposition:Schwartz-2} implies
		$$
		\partial_\mu \partial_\mu f \Big( x_0, \mu, x_1, x_2\Big) \cdot v_1 \otimes v_2 = \partial_\mu \partial_\mu f \Big( x_0, \mu, x_2, x_1\Big) \cdot v_2 \otimes v_1.  
		$$
		Hence, a key insight of Theorem \ref{thm:Schwarz-Lions} is that the order of the free variables generated by the application of Lions derivatives is not key. Thus, we typically represent the collection of free variables using an unordered sequence (a set) rather than an ordered sequence. 
	\end{remark}
	
	\section{Coupled algebraic structures}
	\label{section:ProbabilisticRoughPaths}
	
	A rough differential equation has embedded in the heart of its solution structure a Taylor expansion which includes a polynomial of tensor products plus a remainder term that does not make meaningful contributions at the microscopic level. These polynomials provide a structure that the driving rough path must necessarily satisfy. In the case of weak geometric rough paths, this is encoded via the shuffle product while for branched rough paths this takes the form of the tree product. 
	
	In the mean-field setting, things are a little different: the tensor product has an integration operator associated to it determined by the particular Lions derivative with which this tensor product is paired. To emphasise how important this integral operator is in determining whether terms contribute or not to the Taylor expansion, let $W$ be a Brownian motion on a canonical probability space and consider the difference between
	\begin{equation*}
		\bE^{1, 2} \Big[ W_{s, t}(\omega_1) \otimes W_{s, t} (\omega_2) \Big] 
		\quad \mbox{and} \quad 
		\bE^{1} \Big[ W_{s, t}(\omega_1) \otimes W_{s, t} (\omega_1) \Big],
	\end{equation*}
	where the first expectation is understood as being taken over $(\omega_{1},\omega_{2})$ with respect to the product of the Wiener measure and the second one is understood as being taken over the sole $\omega_{1}$. The first matrix is $0$ while the second one is a non-zero variance matrix for the increment $W_{s,t}$. In connection with the Lions-Taylor expansion stated in Theorem \ref{theorem:LionsTaylor2}, the above increments should be understood as terms of the form $y_{a_{i}} - x_{a_{i}}$ in \eqref{eq:rDa}. In the statement of Theorem \ref{theorem:LionsTaylor2}, the integrals are set over the Euclidean state space, which comes from an implicit application of the transfer Theorem in measure theory. Those integrals originate from expectations defined on the path space and the above two examples are very simple instances of this. In this framework, the labels $\omega_{1}$ and $\omega_{2}$ are 
	somewhat the analogues of $(x_{1},y_{1})$ and $(x_{2},y_{2})$ in 
	\eqref{eq:rDa}.

	This emphasises the fact that $\omega_{1}$ and $\omega_{2}$ are indeed attached to free variables occurring in our Lions derivatives. That said, the reader may find it useful to consider another point of view on the practical interpretation of $\omega_{1}$ and $\omega_{2}$: Focusing on the particle system \eqref{eq:part:sys}, $\omega_{1}$ and $\omega_{2}$ should be understood as continuous labels and, consistently, we should think of $W^1$ and $W^2$ in \eqref{eq:part:sys} as $W^1_{s,t}(\omega)=W_{s,t}(\omega_1)$ and $W^2_{s,t}(\omega)=W_{s,t}(\omega_2)$, with $\omega=(\omega_1, ..., \omega_n)$ belonging to the Wiener space $C([0,1];\bR^d)$ and $W$ being the canonical random variable. 
	Thus we specify the ``statistical correlations'' between those ``labels'' when performing products. In turn, this yields several different contexts for performing the product of two trees (depending on the choice of these correlations).
	
	Hence, there is no longer a ``single product'' from which mean-field Taylor expansions can be iteratively defined. There are going to be 5 different bilinear forms, each one associated with a different element of $A_2[0]$. Further, there are 15 trilinear forms, each one associated with a different element of $\A{3}$. These \emph{couplings} provide a different kind of polynomial expansion that is more closely aligned with the Lions-Taylor expansions developed in Section \ref{section:TaylorExpansions}. The goal of this section is to document how these couplings create a new theory on which we will build our rough path for mean-field equations. 
	
	\subsection{From rough particle systems to rough McKean-Vlasov dynamics}
	
	\subsubsection*{A heuristic introduction to approximating systems of interacting equations}
	
	\emph{Branched rough paths} and the Butcher-Connes-Kreimer Hopf algebra provide a formal setting within which one can consider interacting systems of equations of the form
	\begin{equation}
		\label{eq:particle:system}
		dX_{t}^{i, N} = f\Big( X_{t}^{i, N}, \bar{\mu}_N\big[ \boldsymbol{X}_t \big] \Big) \cdot dW_{t}^{i, N} \quad \bar{\mu}_N\big[ \boldsymbol{X}_t \big] = \tfrac{1}{N} \sum_{j=1}^N \delta_{X_t^{j, N}}\quad i \in \{1, ..., N \}. 
	\end{equation}
	where we denoted $\boldsymbol{X}_t = \big( X_t^{1, N}, ..., X_t^{N, N} \big)$. Thus, the dynamics of the entire system can be expressed in terms of the rough differential equation
	\begin{equation*}
		d\boldsymbol{X}_t = \boldsymbol{f}\big( \boldsymbol{X}_t \big) \cdot d\boldsymbol{W}_t, 
		\quad
		\boldsymbol{f}\big( \boldsymbol{x} \big) = \bigoplus_{i=1}^N f\Big( x^{i, N} , \bar{\mu}_N\big[\boldsymbol{x}\big] \Big)
		\quad \mbox{and}\quad
		\boldsymbol{x} = \big( x^{1, N}, ..., x^{N, N} \big). 
	\end{equation*}
	
	Suppose for the moment that the function $f:\bR^e \times \cP_2(\bR^e) \to \lin(\bR^d, \bR^e)$ is infinitely differentiable and each $W^i:[0,1] \to \bR^d$ is a smooth paths so there are no additional considerations to the underlying calculus. By expressing the system in terms of a differentiable vector field over the entire system and fixing $s, t\in [0,1]$, we are able to use classical Taylor expansion techniques to find an approximation of the form
	\begin{align*}
		\boldsymbol{X}_{s, t} =& \int_s^t \boldsymbol{f}\Big( \boldsymbol{X}_r \Big) \cdot d\boldsymbol{W}_r
		\approx
		\sum_{\boldsymbol{j}} \frac{1}{|\boldsymbol{j}|} \int_s^t \nabla_{\boldsymbol{j}} \boldsymbol{f}\Big( \boldsymbol{X}_r \Big) \cdot \bigotimes_{\iota=1}^{|\boldsymbol{j}|} X_{s, r}^{j_{\iota}, N} \cdot d\boldsymbol{W}_r
	\end{align*}
	where the multi-index $\boldsymbol{j}$ is understood as $\boldsymbol{j} = (j_\iota)_{\iota=1, ..., |j|}$ and where, here and throughout, $\approx$ indicates that equality holds but at the price of an additional remainder term whose contributions is negligible at this stage of the presentation. 
	
	Let us assume that the coefficient $f$ is infinitely differentiable and that the driving signal $W:[0,1] \to \bR^d$ is a smooth path. By approximating $X_{s, t}^{j, N} = f(X_s^{j, N}) \cdot W_{s, t}^{j, N}$ plus a remainder term, we get
	\begin{align}
		\label{eq:heuristic-Approx-+}
		X_{s, t}^{i, N}
		\approx&
		\sum_{\boldsymbol{j} } \frac{1}{|\boldsymbol{j}|!} \nabla_{\boldsymbol{j}} f \Big( X_s^{i, N}, \bar{\mu}_N\big[ \boldsymbol{X}_s \big] \Big)
		\cdot 
		\bigotimes_{k=1}^{|\boldsymbol{j}|} f\Big( X_s^{j_k, N}, \bar{\mu}_N\big[ \boldsymbol{X}_s \big] \Big) 
		\cdot 
		\int_{s}^{t}\Big( \bigotimes_{k=1}^{|\boldsymbol{j}|} W_{s,r}^{j_k, N} \Big) \otimes dW_r^{i, N}. 
	\end{align}
	
	On the one hand, this setting is well documented and we are able to directly identify the products of derivatives and iterated integrals in terms of the labelled, directed trees using the branched rough path techniques described in \cite{gubinelli2010ramification}. On the other hand, thanks to Proposition \ref{proposition:classicTay<=>LionsTay*} we have that Equation \eqref{eq:heuristic-Approx-+} is equivalent to
	\begin{align}
		\nonumber
		X_{s, t}^{i, N} \approx&\sum_{a\in A[0]} \frac{1}{|a|!} \frac{1}{N^{m[a]}} \sum_{j_1, ..., j_{m[a]} = 1}^N \partial_a f\Big( X_s^{i, N}, \bar{\mu}_N\big[ \boldsymbol{X}_s \big], \boldsymbol{X}_s^{\boldsymbol{j}\circ m\{a\}, N} \Big) \cdot \bigotimes_{k=1}^{|a|} f\Big( X_{s}^{j_{a_k}, N}, \bar{\mu}_N\big[ \boldsymbol{X}_s \big] \Big)
		\\
		\label{eq:heuristic-Approx-Particle}
		&\qquad \cdot \int_s^t \Big( \bigotimes_{k=1}^{|a|} W_{s, r}^{j_{a_k}, N} \Big) \otimes dW_r^{i, N}, 
	\end{align}
	where 
	\begin{equation*}
		\boldsymbol{X}_s^{\boldsymbol{j} \circ m\{a\}, N} = \big( X_s^{j_1, N}, ..., X_s^{j_{m[a]}, N} \big)
		\quad \mbox{and} \quad 
		X_s^{j_0, N} = X_s^{i, N} . 
	\end{equation*}
	
	An equivalent formulation for Equation \eqref{eq:particle:system} is to consider the canonical lift of the functional $f:\bR^e \times \cP_2(\bR^e) \to \lin(\bR^d, \bR^e)$ to the function $F:\bR^e \times L^2(\Omega', \bP'; \bR^e) \to \lin(\bR^d, \bR^e)$ and rewrite Equation \eqref{eq:particle:system} as 
	\begin{equation}
		dX_t^{i, N} = F\Big( X_t^{i, N}, X_t^{u, N} \Big) \cdot dW_t^{i, N}, 
		\quad 
		u:\Omega' \to \{1, ..., N\} 
		\quad \mbox{and} \quad 
		i \in \{1, ..., N\}
	\end{equation}
	where $\bP' \circ (u)^{-1}$ is the uniform distribution over the set $\{1,..., N\}$. The canonical lift of the distribution $\bar{\mu}_N\big[ \boldsymbol{X}_t \big]$ is the random variable $X_t^{u(\omega')}(\omega_0)$ on the product probability space $\Omega \times \Omega'$. Of course, $\Omega$ and $\Omega'$ could be identified, but we feel it improves the clarity of the expansions to emphasise that $\Omega$ (the probability space on which the solution exists) and $\Omega'$ (the probability space on which the lift of the solution law exists) are distinguished. 

	Using the same approximation techniques to obtain Equation \eqref{eq:heuristic-Approx-Particle}, we can replace the various empirical means by an expectation of a sequence of independent uniformly distributed random variable on the set $\{1, ..., N\}$ and obtain
	\begin{align}
		\nonumber
		X_{s, t}^{i, N}&(\omega_0) \approx \sum_{a\in A[0]} \frac{1}{|a|!}  (\bE')^{m\{a\}}\bigg[ \partial_a f\Big( X_s^{i, N}(\omega_0), \bar{\mu}_N\big[ \boldsymbol{X}_s(\omega_0) \big], X_s^{u(\omega'), N}(\omega_0)_{m\{a\}} \Big) 
		\\
		\label{eq:heuristic-Approx-Particle2}
		&\cdot \bigotimes_{k=1}^{|a|} f\Big( X_{s}^{u(\omega_{a_k}'), N}(\omega_0), \bar{\mu}_N\big[ \boldsymbol{X}_s(\omega_0) \big] \Big)
		\cdot \int_s^t \Big( \bigotimes_{k=1}^{|a|} W_{s, r}^{u(\omega_{a_k}'), N}(\omega_0) \Big) \otimes dW_r^{i, N}(\omega_0) \bigg], 
	\end{align}
	where
	\begin{equation*}
		X_s^{u(\omega'), N}(\omega_0)_{m\{a\}} = \Big( X_s^{u(\omega_1'), N}(\omega_0), ..., X_s^{u(\omega_{m[a]}'), N}(\omega_0) \Big)
		\quad \mbox{and} \quad
		X_s^{u(\omega'_0), N}(\omega_0) = X_s^{i, N}(\omega_0). 
	\end{equation*} 
	
	In Equation \eqref{eq:heuristic-Approx-Particle2}, the expectation $(\bE')^{m\{a\}}$ runs over the product probability space 
	\begin{equation*}
		\big( (\Omega')_1, \cF_1, \bP_1 \big) \times ... \times \big( (\Omega')_{m[a]}, \cF_{m[a]}, \bP_{m[a]} \big),
	\end{equation*} 
	a generic element of which is denoted $\omega_{m\{a\}} = (\omega_1, ..., \omega_{m[a]})$. 
	
	
	\subsubsection*{A comparison with the associated McKean-Vlasov mean-field limit}
	\label{subsubsec:Heuristic-McKean}
	
	On the other hand, let us consider the distribution dependent equation
	\begin{equation}
		\label{eq:meanfield:equation}
		dX_t = f\Big( X_t, \cL_t^X \Big) \cdot dW_t, \quad X_0 = \xi, \quad \cL_t^X = \bP \circ (X_t)^{-1}. 
	\end{equation}
	This is often referred to as a McKean-Vlasov equation as the coefficients are dependent on the position and law of the solution. These equations describe the limiting statistical behaviour of dynamics such as Equation \eqref{eq:particle:system} as $N\to \infty$. 
	
	Suppose as before that the coefficient $f:\bR^e \times \cP_2(\bR^e) \to \lin(\bR^d, \bR^e)$ is infinitely differentiable and $W:[0,1] \to \bR^d$ is a smooth path so there are no additional considerations to the underlying calculus. Using Theorem \ref{theorem:LionsTaylor2}, we have that the dynamics of Equation \eqref{eq:meanfield:equation} can be approximated for $s, t\in [0,1]$ by
	\begin{align}
		\nonumber
		X_{s, t}(\omega_0) =& \int_s^t f \Big( X_r(\omega_0), \cL_r^X \Big) \cdot dW_r(\omega_0) = \int_s^t F \Big( X_r(\omega_0), X_r \Big) \cdot dW_r(\omega_0)
		\\
		\label{eq:heuristic-Approx}
		\approx& \sum_{a\in A[0]} \frac{1}{|a|!} \int_{s}^{t} \rD^a f\Big( X_s(\omega_0), \cL_s^X \Big)\Big[ X_{s, r}(\omega_0), \Pi_{s, r}^X \Big] \cdot dW_r, 
	\end{align}
	where again $\approx$ indicates that equality holds but at the price of an additional remainder whose contribution is negligible for the purposes of this presentation. Next, by approximating 
	\begin{equation*}
		X_{s, r} \approx f(X_s, \cL_s^X) \cdot W_{s, r} 
		\quad \mbox{and}\quad
		\Pi_{s, r}^X \approx \cL\big( X_s, X_s + f(X_s, \cL_s^X) \cdot W_{s, r} \big), 
	\end{equation*} 
	we obtain the Lions elementary differentials and iterated integrals of the driving noise $W$ that exists on some product probability space: 
	\begin{align}
		\nonumber
		X_{s, t}(\omega_0)
		\approx&
		\sum_{a\in A[0] } \frac{1}{|a|!} \bE^{m\{a\}}\bigg[  \partial_a f\Big( X_s(\omega_0), \cL_s^X, X_s(\omega)_{m\{a\}} \Big)
		\cdot 
		\bigotimes_{k=1}^{|a|} f\Big(X_s(\omega_{a_k}), \cL_s^X\Big) 
		\\
		\label{eq:heuristic-Approx-}
		&\qquad 
		\cdot 
		\int_{s}^{t}\Big( \bigotimes_{k=1}^{|a|} W_{s,r}(\omega_{a_k}) \Big) \otimes dW_r(\omega_0) \bigg], 
	\end{align}
	where we have denoted the sequence of random variables
	\begin{equation*}
		X_s(\omega)_{m\{a\}} = \big( X_s(\omega_1), ..., X_s(\omega_{m[a]}) \big). 
	\end{equation*}
	
	Thus, when each of the driving signals of some interacting rough differential are taken to be statistically exchangeable and the coefficients of an equation only interact with other equations via an empirical measure, we hope that we can reduce the amount of information about the driving signal required to solve the McKean-Vlasov rough differential equation by replacing the lattice of iterated integrals of the driving signals by some distribution which captures the statistical nature of that lattice of iterated integrals. 
	
	If the world were a reasonable place, we would expect that as the number of equations within the collective system described by Equation \eqref{eq:particle:system} were taken to be large while the collection of driving signals were statistically exchangeable, the dynamics could be approximated by the dynamics of the associated McKean-Vlasov equation \eqref{eq:meanfield:equation}. Hence, we would also hope that the infinitesimal dynamics would also match up (so that there is an inter-relationship between the symbols used for the index set of the abstract Lions-Taylor expansion described in Equation \eqref{eq:heuristic-Approx-} and the labelled directed trees used as the index set for the abstract Taylor expansion in \eqref{eq:heuristic-Approx-+}). 

	\subsubsection*{Representation of abstract Lions-Taylor expansions}
	
	We fix $N\in \bN$ and $i\in \{1, ..., N\}$. Dropping remainder terms for the moment, let us use the symbol $\sim$ in order to associate below the key components of the expansions from Equations \eqref{eq:heuristic-Approx-Particle2} and \eqref{eq:heuristic-Approx-} with abstract symbols (with each symbol, we associate two quantities, one for the particle
	system and one for the McKean-Vlasov equation): for $\iota \in \{1, ..., d\}$
	\begin{equation*}
		\begin{tikzpicture}
			\node[vertex, label=right:{\footnotesize $\iota$}] at (0,0) {};
		\end{tikzpicture}
		\sim 
		\left\{ 
		\begin{aligned}
			&f^\iota\Big( X_s^{i, N}(\omega_0), \bar{\mu}_N\big[ \boldsymbol{X}_s(\omega_0) \big] \Big) \cdot W_{s, t}^{\iota, (i, N)}(\omega_0)
			\\
			&f^\iota\Big( X_s(\omega_0), \cL_s^X \Big) \cdot W_{s, t}^{\iota}(\omega_0)
		\end{aligned}
		\right.
	\end{equation*}
	and
	\begin{align}
		\label{eq:partitioning1}
		&
		\begin{tikzpicture}
			\node[vertex, label=right:{\footnotesize $\iota$}] at (0,0) {};
			\node[vertex, label=right:{\footnotesize $\iota_1$}] at (-2,1) {};
			\node[vertex, label=right:{\footnotesize $\iota_2$}] at (-1,1) {};
			\node[vertex, label=right:{\footnotesize $\iota_{|a|-1}$}] at (0.75,1) {}; 
			\node[vertex, label=right:{\footnotesize $\iota_{|a|}$}] at (2,1) {}; 
			\node at (0,1) {...};
			\node at (-2.5,1) {$\cE^a\Big[ $};
			\node at (2.75,1) {$\Big]$};
			\draw[edge] (0,0) -- (-1,1);
			\draw[edge] (0,0) -- (-2,1);
			\draw[edge] (0,0) -- (0.75,1);
			\draw[edge] (0,0) -- (2,1);
		\end{tikzpicture}
		\\
		&\sim \left\{
		\begin{aligned}
			(\bE')^{m\{a\}}&\bigg[ \partial_a f^{\iota}\Big( X_s^{i, N}(\omega_0), \bar{\mu}_N\big[ \boldsymbol{X}_s(\omega_0) \big], X_s^{u(\omega'), N}(\omega_0)_{m\{a\}} \Big) 
			\\
			\centerdot &\prod_{k=1}^{|a|} f^{\iota_k}\Big( X_{s}^{u(\omega_{a_k}'), N}(\omega_0), \bar{\mu}_N\big[ \boldsymbol{X}_s(\omega_0) \big] \Big)
			\cdot \int_s^t \Big( \prod_{k=1}^{|a|} W_{s, r}^{\iota_k, (u(\omega_{a_k}'), N)}(\omega_0) \Big) \times dW_r^{\iota, (i, N)}(\omega_0) \bigg] 
			\\
			\bE^{m\{a\}}&\bigg[ \partial_a f^{\iota} \Big( X_s(\omega_0), \cL_s^X, X_s(\omega)_{m\{a\}} \Big) 
			\\
			\centerdot &\prod_{k=1}^{|a|} f^{\iota_k} \Big( X_s(\omega_{a_k}), \cL_s^X \Big) \cdot \int_s^t \Big( \prod_{k=1}^{|a|} W_{s, r}^{\iota_k}(\omega_{a_k}) \Big) \times dW^{\iota}_r(\omega_0) \bigg]
		\end{aligned}
		\right.
	\end{align}
	where $a\in A[0]$ is a partition sequence and $\cE^a$ is some operator that for the moment we only think of as ``partitioning'' the elements of the second level. Thanks to Theorem \ref{thm:Schwarz-Lions} and commutativity of scalar products, we should additionally emphasise that this tree structure is non-planar in the sense that we can commute the directional derivatives via a shuffling of the partition sequence. 
	
	While the graph \eqref{eq:partitioning1} is used to represent a product of integrals (together with several additional integration operators), we can easily guess that similar symbols should be necessary to represent higher order iterated integrals. In fact, by continuously repeating this technique, we will obtain more and more diffeomorphisms that are dependent only on the pair $(X_s, \cL_s^X)$ (respectively $(X_s^{i, N}, \bar{\mu}_N[\boldsymbol{X}_s])$) in $\bR^e \times \cP_2(\bR^e)$ and the iterated integrals of the driving signal. Further, we can (informally) see that this process leads to terms that we identify with tree structures that carry an additional partition structure. A natural question is therefore what does the abstract symbol
	\begin{equation}
		\label{eq:partitioning2}
		\begin{aligned}
			\begin{tikzpicture}
				\node[vertex, label=right:{\footnotesize $i$}] at (0,0) {};
				\node[vertex, label=right:{\footnotesize $i_1$}] at (-2,1) {};
				\node[vertex, label=right:{\footnotesize $i_2$}] at (-1,1) {};
				\node[vertex, label=right:{\footnotesize $i_{|a|-1}$}] at (0.75,1) {};
				\node[vertex, label=right:{\footnotesize $i_{|a|}$}] at (2,1) {};
				\node[vertex, label=right:{\footnotesize $i_1'$}] at (-1,2) {};
				\node[vertex, label=right:{\footnotesize $i_2'$}] at (0,2) {};
				\node[vertex, label=right:{\footnotesize $i_{|a'|-1}'$}] at (1.75,2) {};
				\node[vertex, label=right:{\footnotesize $i_{|a'|}'$}] at (3,2) {};
				\node at (0,1) {...};
				\node at (1,2) {...};
				\node at (-2.5,1) {$\cE^a\Big[ $};
				\node at (2.75,1) {$\Big]$};
				\node at (-1.5,2) {$\cE^{a'}\Big[ $};
				\node at (4,2) {$\Big]$};
				\draw[edge] (0,0) -- (-1,1);
				\draw[edge] (0,0) -- (-2,1);
				\draw[edge] (0,0) -- (0.75,1);
				\draw[edge] (0,0) -- (2,1);
				\draw[edge] (0.75,1) -- (1.75,2);
				\draw[edge] (0.75,1) -- (3,2);
				\draw[edge] (0.75,1) -- (0,2);
				\draw[edge] (0.75,1) -- (-1,2);
				\node at (5,1) {$=...?$};
			\end{tikzpicture}
		\end{aligned}
	\end{equation}
	correspond to and how should we keep track of the associated probability spaces needed to give this meaning. 
	
	\subsection{Lions trees and their associated algebra}
	\label{subsection:2.1LionsTrees}
	
	Motivated by the theory of rough paths and the Butcher-Connes-Kreimer Hopf algebra (see for instance \cite{gubinelli2010ramification}), we start this section by introducing notation and concepts that are central to this work. 
	
	By a directed graph, we mean a pair $(\scN, \scE)$ where $\scN$ is a set of distinct elements which we call nodes (sometimes referred to as vertices in the literature) and $\scE\subseteq \scN\times \scN$ which we call edges. A directed tree is a graph that is connected, acyclic and all directed edges point towards a single node which we call the root. If a directed graph has a finite number of connected components, each of which is a directed tree, we refer to it as a directed forest. When referring to more than one directed trees, it will always be assumed that any two sets of nodes are disjoint. 
	
	For a directed forest $(\scN, \scE)$ with set of roots $\fr(\scN)$, every element $y\in \scN$ has a unique sequence $(y_i)_{i=1, ..., n} \in \scN$ such that $y_1=y$, $y_n \in \fr(\scN)$ and $(y_i, y_{i+1}) \in \scE$ for $i=1, ..., n-1$.
	
	There is a Reflexive, Transitive binary relationship $\leq$ over the nodes of a directed forest $\scN$ determined by the number of edges along the unique path from each node to a root. 
	
	That is $x \leq y \iff$
	\begin{align*}
		&\exists! (x_i)_{i=1, ..., m}: x_1 = x, x_m \in \fr(\scN) \mbox{ and } \forall i=1, ..., m-1, (x_i, x_{i+1}) \in \scE, 
		\\
		&\exists! (y_i)_{i=1, ..., n}: y_1 = y, y_n \in \fr(\scN) \mbox{ and } \forall i=1, ..., n-1, (y_i, y_{i+1}) \in \scE,
		\\
		&\quad \mbox{and} \quad m \leq n.  
	\end{align*}
	When $x\leq y$ and not $y \leq x$, we additionally denote the Transitive binary relationship $<$. Thus $(\scN, \leq)$ is a preorder and the set of all nodes such that
	$$
	x \leq \geq y \quad \mbox{or equivalently} \quad x\leq y \mbox{ and } y \leq x
	$$
	form equivalence classes of the nodes. 
	
	We denote the set of directed forests by $\fF$, the set of labelled, directed forests $(\scN, \scE, \scL)$ where $\scL: \scN \to \{1, ..., d\}$ by $\fF_{d}$,  the set of directed trees by $\fT$ and the set of labelled, directed trees by $\fT_{d}$. Further, we denote by $\rId = (\emptyset, \emptyset, \emptyset)$ the empty forest and $\fF_{0, d} = \fF_d \cup\{\rId\}$. 
	
	For $x \in \scN$, a label $\scL(x)$ should be understood as the dimensional coordinate of a particle attached to the node $x$ but this is in no way a label of the particle itself. In our construction, particles have no distinctive labellings. 
	
	We start with the following notion which will be used intensively for the rest of the paper:
	\begin{definition}
		Let $\scN$ be a non-empty set containing a finite number of elements and let $H \subseteq 2^{\scN}\backslash \emptyset$ be a collection of subsets of $\scN$. Then we say that the pair $(\scN, H)$ is a hypergraph. The elements $h\in H$ are referred to as hyperedges. 
	\end{definition}
	
	A hypergraph $(\scN, H)$ is a generalisation of a graph in which "edges" $h\in H$ can contain any positive number of nodes, rather than specifically two. Hypergraphs are sometimes referred to as \emph{range spaces} in computational geometry, \emph{simple games} in cooperative game theory and in some literature hyperedges are referred to as \emph{hyperlinks} or \emph{connectors}. 	
	A hyperedge is said to be \emph{$d$-regular} if every node is contained in exactly $d$ hyperedges. 
	
	In this work, we make no assumptions about the users background in hypergraph theory. However, the curious reader may choose to refer to \cite{Bretto2014Hypergraph} for a more general introduction to the theory of hypergraphs. 
	\begin{example}
		\label{example:3:2}
		The following are all visual representations of 1-regular hypergraphs:
		$$
		
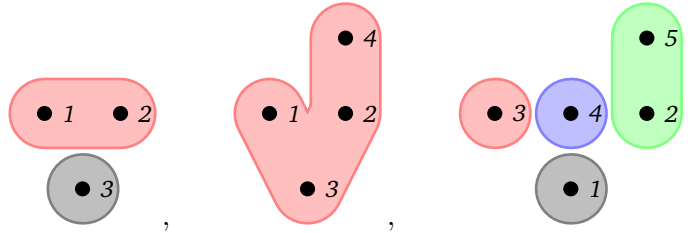
\begin{tikzpicture}
			\node[vertex, label=right:{\footnotesize 3}] at (0,0) {}; 
			\node[vertex, label=right:{\footnotesize 2}] at (0.5,1) {}; 
			\node[vertex, label=right:{\footnotesize 1}] at (-0.5,1) {}; 
			\begin{pgfonlayer}{background}
				\draw[zhyedge1] (0,0) -- (-0,0);
				\draw[zhyedge2] (0,0) -- (-0,0);
				\draw[hyedge1, color=red] (-0.5,1) -- (0.5,1);
				\draw[hyedge2] (-0.5,1) -- (0.5,1);
			\end{pgfonlayer}
		\end{tikzpicture}
		,\qquad
		\begin{tikzpicture}
			\node[vertex, label=right:{\footnotesize 3}] at (0,0) {}; 
			\node[vertex, label=right:{\footnotesize 2}] at (0.5,1) {}; 
			\node[vertex, label=right:{\footnotesize 1}] at (-0.5,1) {}; 
			\node[vertex, label=right:{\footnotesize 4}] at (0.5,2) {}; 
			\begin{pgfonlayer}{background}
				\draw[hyedge1, color=red] (0.5,2) -- (0.5,1) -- (0,0) -- (-0.5,1);
				\draw[hyedge2] (0.5,2) -- (0.5,1) -- (0,0) -- (-0.5,1);
			\end{pgfonlayer}
		\end{tikzpicture}
		, \qquad
		\begin{tikzpicture}
			\node[vertex, label=right:{\footnotesize 1}] at (0,0) {}; 
			\node[vertex, label=right:{\footnotesize 2}] at (1,1) {}; 
			\node[vertex, label=right:{\footnotesize 3}] at (-1,1) {}; 
			\node[vertex, label=right:{\footnotesize 4}] at (0,1) {}; 
			\node[vertex, label=right:{\footnotesize 5}] at (1,2) {}; 
			\begin{pgfonlayer}{background}
				\draw[zhyedge1] (0,0) -- (0,0);
				\draw[zhyedge2] (0,0) -- (0,0);
				\draw[hyedge1, color=red] (-1,1) -- (-1,1);
				\draw[hyedge2] (-1,1) -- (-1,1);
				\draw[hyedge1, color=blue] (0,1) -- (0,1);
				\draw[hyedge2] (0,1) -- (0,1);
				\draw[hyedge1, color=green] (1,1) -- (1,2);
				\draw[hyedge2] (1,1) -- (1,2);
			\end{pgfonlayer}
		\end{tikzpicture}. 
		$$
		These correspond to the hypergraphs 
		$$
		\Big( \{1,2,3\}, \big\{ \{1,2\}, \{3\}\big\} \Big), \quad
		\Big( \{1,2,3,4\}, \big\{ \{1,2,3,4\}\big\}\Big), \quad
		\Big( \{1,2,3,4,5\}, \big\{ \{1\}, \{2,5\}, \{3\}, \{4\}\big\}\Big).
		$$ 
		We emphasise that despite the aethetics of these examples, these are not trees (or even graphs). 
	\end{example}
	As we stated earlier, particles have no distinctive labellings. All particles are tagged through a hypergraph, in the sense that two different hyperedges are implicitly associated with two distinct particles. In short, a tree (or a forest) is associated with a finite collection of distinct particles (the number of particles being equal to the number of hyperedges), no matter which particles these are precisely. This is in stark contrast with the trees from the Butcher-Connes-Kreimer algebra associated with the finite-dimensional particle system where any node carries an additional label, which is specific to a given particle.
	
	The following definition is at the centre of this paper:
	\begin{definition}
		\label{definition:Forests}
		Let $\scN$ be a non-empty set containing a finite number of elements. Let $\scE\subset \scN\times \scN$ such that $(x,y)\in \scE \implies (y,x)\notin \scE$. Let $h_0 \subseteq \scN$ and $H\in \scP(\scN\backslash h_0)$ is a partition of $\scN\backslash h_0$. We denote $H' = (H \cup \{h_0\}) \backslash \{ \emptyset\}$. Let $\scL:\scN \to \{1, ..., d\}$. We refer to $T=(\scN, \scE, h_0, H)$ as a \emph{Lions forest} if
		\begin{enumerate}
			\item $(\scN, \scE)$ is a directed forest with preordering $\leq$. 
			\item $(\scN, H')$ is a 1-regular hypergraph where the hyperedges $h_0, h \in H'$ satisfy:
			\begin{enumerate}[label=(2.\arabic*)]
				\item 
				\label{definition:Forests:2.1}
				If $h_0 \neq \emptyset$, then $\exists x\in h_0$ such that $\forall y\in \scN$, $x\leq y$. 
				\item 
				\label{definition:Forests:2.2}
				For $h_i \in H'$, suppose $x, y\in h_i$ and $x<y$. Then $\exists z \in h_i$ such that $(y, z)\in \scE$. 
				\item 
				\label{definition:Forests:2.3}
				For $h_i \in H'$, suppose $x_1, y_1 \in h_i$, $x_1 \leq \geq y_1$, $(x_1, x_2), (y_1, y_2) \in \scE$ and $x_2 \neq y_2$. Then $x_2, y_2 \in h_i$. 
			\end{enumerate}
		\end{enumerate}
		The collection of all such Lions forests is denoted by $\scF$. When a Lions forest $T=(\scN,\scE,h_0, H)$ satisfies that $(\scN, \scE)$ is a directed tree, it is referred to as a Lions tree and the collection of all lions trees is denoted $\scT$. We define $\scF_0= \scF\cup\{\rId\}$ where $\rId=(\emptyset, \emptyset, \emptyset, \emptyset, \scL)$  is the empty tree. On a separate note, we define $\scT_0 = \{T \in \scT: h_0^T \neq \emptyset\}$. 
		
		Equivalently, we refer to $T=(\scN, \scE, h_0, H, \scL)$ as a \emph{labelled Lions forest} if $(\scN, \scE, h_0, H)$ is a Lions forest and $\scL$ is a labelling taking values in $\{1, ..., d\}$. We denote the collection of labelled Lions forests $\scF_d$, and the set of labelled Lions trees by $\scT_d$. Finally, we denote $\scF_{0,d}= \scF_d \cup \{\rId\}$ and $\scT_{0, d} = \{T \in \scT_d: h_0^T \neq \emptyset\}$. 
	\end{definition}
	
	\begin{remark}
		\label{remark:labellings-hypergraphs}
		Labels $\scL$ here could be understood as labels of the particles in the particle system, bearing in mind that such particles carry in fact two types of labels: One label for the dimensional coordinate and a second label that identifies the particle. The above labels are implicitly understood as labels for the dimensional coordinate. 
		
		It is our choice not to put labels for the numbering of the particles, except for one tagged reference particle, which we mark by distinguishing a hyperedge $h_0$ from all others. Hyperedges should be thought as clusters of the same particle (or equivalently of particles equipped with the same number). However, hyperedges should satisfy additional properties, as formalised by Definition \ref{definition:Forests}
		
		This ``tagged'' particle has the same role as the tagged particle for standard McKean-Vlasov equations. The others are regarded as forming an infinite cloud, which we sometimes call a ``continuum''.
	\end{remark}
	
	Motivated by Section \ref{subsubsec:Heuristic-McKean}, our goal is to show that the hypergraphic structure described in Definition \ref{definition:Forests} arises as a result of the properties of the Lions derivatives. 
	
	\begin{theorem}
		\label{theorem:EquivalenceTrees}
		Let $(\scN, \scE, \scL)$ be a labelled directed forest with roots $\fr(\scN) \subseteq \scN$. Let $\hat{h}_0$ be a (possibly empty) subset of $\fr(\scN)$ and let $\hat{H} \in \scP\big( \fr(\scN)\backslash \hat{h}_0 \big)$. 
		
		For $x \in \scN$, we denote the set of children of $x$ by
		\begin{equation*}
			\scN_x := \big\{ y \in \scN: (y, x) \in \scE\big\}. 
		\end{equation*}
		Let $\fh: \scN \to 2^{2^\scN}$ such that
		\begin{equation*}
			\fh[x] \in \scP\Big( \scN_x \cup \{x\} \Big). 
		\end{equation*}
		We say that $\big( \scN, \scE, (\hat{h}_0, \hat{H}), \fh, \scL \big)$ is a \emph{Lions partition forest}. We denote $\hat{\scF}$ to be the set of all Lions partition forests. 
		
		Then there is an isomorphism between the set of all Lions forests $\scF$ and the set of all Lions partition forests $\hat{\scF}$. As such, we view these two definitions as equivalent. 
	\end{theorem}

	Theorem \ref{theorem:EquivalenceTrees}, which is proved in \cite{salkeld2022LionsTrees}, is critical because it captures the equivalence between the Lions hypergraphic structure that arises from Conditions \ref{definition:Forests:2.1}, \ref{definition:Forests:2.2} and \ref{definition:Forests:2.3} and the separate hypergraphic structure that arises from partitioning the tree ``locally''. This is the fundamental result that links the partions of Lions trees of Definition \ref{definition:Forests} and the partition structure of the elementary differentials described in Equations \eqref{eq:partitioning1} and \eqref{eq:partitioning2} previously. 
	
	\subsubsection{Operations on Lions trees}
	
	The set of all directed forests $\fF_0$ is described via a pair of operations which are referred to as the forest product and the rooting operation. For instance, the forest product $\odot:\fF_{0, d} \times \fF_{0, d} \to \fF_{0, d}$ is defined so that given $\tau_1=(\scN^{\tau_1}, \scE^{\tau_1}, \scL^{\tau_1})$ and $\tau_2 = (\scN^{\tau_2}, \scE^{\tau_2}, \scL^{\tau_2})$, we have 
	\begin{equation}
		\label{eq:Tree_product}
		\begin{aligned}
			&\tau_1 \odot \tau_2 = (\tilde{\scN}, \tilde{\scE}, \tilde{\scL}) 
			\quad \mbox{such that} \quad
			\\
			&\tilde{\scN} = \scN^{\tau_1} \cup \scN^{\tau_2}, 
			\quad
			\tilde{\scE} = \scE^{\tau_1} \cup \scE^{\tau_2},
			\quad 
			\tilde{\scL}: \tilde{\scN} \to \{1, ..., d\}
			\quad \mbox{such that} \quad 
			\tilde{\scL}|_{\scN^{\tau_i}} = \scL^{\tau_i}. 
		\end{aligned}
	\end{equation}
	Similarly, for $i \in \{1, ..., d\}$ the forest rooting operator $\lfloor \cdot \rfloor_i: \fF_{0, d} \to \fT_{d}$ is defined so that given $\tau=(\scN^{\tau}, \scE^{\tau}, \scL^{\tau})$, we have $\lfloor \tau \rfloor_i = (\tilde{\scN}, \tilde{\scE}, \tilde{\scL})$ where
	\begin{align*}
		&\tilde{\scN} = \scN \cup \{x_0\}, 
		\quad
		\tilde{\scE} = \scE \cup \big\{ (x, x_0): x \in \scN \quad \mbox{root} \}
		\\
		&\tilde{\scL}:\tilde{\scN} \to \{1, ..., d\} \quad \mbox{such that} \quad \tilde{\scL}|_\scN = \scL, \tilde{\scL}[x_0] = i. 
	\end{align*}

	
	When multiplying two Lions forests, the $h_{0}$-hyperedges are merged; all the others are kept separated. In other words, $h_{0}$-hyperedges of the two forests are implicitly understood as carrying the same particle (that is somehow tagged by a $0$) and all the other hyperedges are implicitly understood as carrying statistically independent particles (a non-$h_{0}$-hyperedge of the first forest carries an independent particle from any other hyperedge of the second forest). The rationale for this definition is that particles in non-$h_{0}$-hyperedges ($h\in H$) are sampled out from the cloud and each sample is independent of other samples. 
	
	In order to state the definition properly, we introduce the notation:
	\begin{definition}
		\label{definition:Lionsproduct}
		We define $\circledast: \scF_{0,d} \times \scF_{0,d} \to \scF_{0,d}$ so that for two Lions forests 
		$$
		T_1=(\scN^1, \scE^1, h_0^1, H^1, \scL^1)
		\quad \mbox{and} \quad 
		T_2 = (\scN^2, \scE^2, h_0^2, H^2, \scL^2), 
		$$
		we have $T_1 \circledast T_2 = (\tilde{\scN}, \tilde{\scE}, \tilde{h_0}, \tilde{H}, \tilde{\scL})$ such that 
		\begin{align*}
			&\tilde{\scN} = \scN^1 \cup \scN^2,
			\quad
			\tilde{\scE} = \scE^1 \cup \scE^2, 
			\quad
			\tilde{h}_0 = h_0^1 \cup h_0^2 
			\quad
			\tilde{H} = H^1 \cup H^2,  
			\\
			&\tilde{\scL}: \tilde{\scN} \to \{1, ..., d\} \quad \mbox{ such that } \quad \tilde{\scL}|_{\scN^i} = \scL^i. 
		\end{align*}
		The $\circledast$ operation is associative and commutative product on the space of forests with unit $\rId$. 
		
		Let $\cE: \scF_{0,d} \to \scF_{0,d}$ be the operator defined for $T=(\scN, \scE, h_0, H, \scL) \in \scF_{0,d}$ by
		$$
		\cE[T] = (\scN, \scE, \emptyset, H', \scL).
		$$ 
		We refer to $\cE$ as the decoupling operation. 
		
		For $i\in\{1, ..., d\}$, the rooting operator $\lfloor \cdot \rfloor_i: \scF_{0, d} \to \scT_{0, d}$ is defined so that 
		\begin{align*}
			\mbox{For} \quad& T=(\scN, \scE, h_0, H, \scL) \quad \lfloor T\rfloor_i = (\tilde{\scN}, \tilde{\scE}, \tilde{h}_0, H, \tilde{\scL}) \quad \mbox{where:}
			\\
			& x_0 \notin \scN,
			\quad
			\tilde{\scN} = \scN \cup \{ x_0\}, 
			\quad
			\tilde{h}_0 = h_0 \cup\{x_0\}.  
			\\
			& \{x_1, ..., x_n\}=\fr[T] \subseteq \scN, 
			\quad \tilde{\scE} = \scE \cup \big\{ (x_1, x_0), ..., (x_n, x_0) \big\}.
		\end{align*}
	\end{definition}
	The decoupling operation takes the $0$-hyperedge $h_0$ (if non-empty) and turns it into a hyperedge carrying a particle sampled uniformly from the cloud. The symbol $\cE$ is chosen because this operation corresponds (in a very specific sense that will be elucidated later) to integrating. The $0$-particle is usually not expected to be integrated, since only the statistical law for the other particles matters. Thus the nodes of the Lions tree contained in $h_0$ are mapped from nodes that are not summed over to nodes that are collectively summed over by mapping $h_0$ into $H$. 

	Let $a \in \A{n}$. We define the operator $\cE^a: (\scF_{0, d})^{\times n} \to \scF_{0, d}$ by
	\begin{equation}
		\label{eq:definition:E^a-notation}
		\cE^a\Big[ T_1, ..., T_n\Big] = \Big[ \underset{\substack{i:\\ a_i=0}}{\circledast} T_i \Big] \circledast \cE\Big[ \underset{\substack{i:\\ a_i=1}}{\circledast} T_i \Big] \circledast ... \circledast \cE\Big[ \underset{\substack{i:\\ a_i=m[a]}}{\circledast} T_i \Big] . 
	\end{equation}
	Intuitively, the decoupling by the partition sequence $\cE^a$ groups the sequence of trees into $m[a]+1$ partitions with a single group tagged and all others detagged. 
	
	There is a fundamental link between Lions forests and the three operators $\circledast$, $\cE$ and $\lfloor \cdot \rfloor$:
	\begin{proposition}
		\label{proposition:CompletionOfTrees}
		The completion of the set $\{ \rId \}$ with respect to the operations $\circledast$, $\cE$ and $\lfloor \cdot\rfloor_i$ such that $i\in\{1, ..., d\}$ is the complete set of forests $\scF_{0, d}$. 
	\end{proposition}
	The proof of Proposition \ref{proposition:CompletionOfTrees} can be found in \cite{salkeld2022LionsTrees}. 
	
	\subsubsection{Visualisations}
	
	To better explain the structures of these forests, we provide the reader with a collection of diagramatic representations. 
	
	\begin{example}
		The following are all trees that satisfy Definition \ref{definition:Forests}
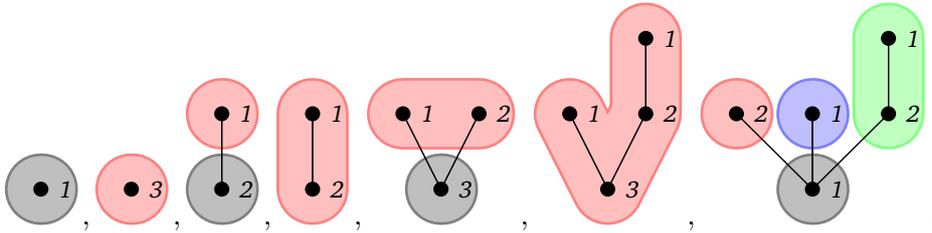
		$$
		\begin{tikzpicture}
			\node[vertex, label=right:{\footnotesize 1}] at (0,0) {}; 
			\begin{pgfonlayer}{background}
				\draw[zhyedge1] (0,0) -- (-0,0);
				\draw[zhyedge2] (0,0) -- (-0,0);
			\end{pgfonlayer}
		\end{tikzpicture}
		, 
		\begin{tikzpicture}
			\node[vertex, label=right:{\footnotesize 3}] at (0,0) {}; 
			\begin{pgfonlayer}{background}
				\draw[hyedge1, color=red] (0,0) -- (-0,0);
				\draw[hyedge2] (0,0) -- (-0,0);
			\end{pgfonlayer}
		\end{tikzpicture}
		, 
		\begin{tikzpicture}
			\node[vertex, label=right:{\footnotesize 2}] at (0,0) {}; 
			\node[vertex, label=right:{\footnotesize 1}] at (0,1) {}; 
			\draw[edge] (0,0) -- (0,1);
			\begin{pgfonlayer}{background}
				\draw[zhyedge1] (0,0) -- (-0,0);
				\draw[zhyedge2] (0,0) -- (-0,0);
				\draw[hyedge1, color=red] (0,1) -- (0,1);
				\draw[hyedge2] (0,1) -- (0,1);
			\end{pgfonlayer}
		\end{tikzpicture}
		,
		\begin{tikzpicture}
			\node[vertex, label=right:{\footnotesize 2}] at (0,0) {}; 
			\node[vertex, label=right:{\footnotesize 1}] at (0,1) {}; 
			\draw[edge] (0,0) -- (0,1);
			\begin{pgfonlayer}{background}
				\draw[hyedge1, color=red] (0,0) -- (0,1);
				\draw[hyedge2] (0,0) -- (0,1);
			\end{pgfonlayer}
		\end{tikzpicture}
		,
		\begin{tikzpicture}
			\node[vertex, label=right:{\footnotesize 3}] at (0,0) {}; 
			\node[vertex, label=right:{\footnotesize 2}] at (0.5,1) {}; 
			\node[vertex, label=right:{\footnotesize 1}] at (-0.5,1) {}; 
			\draw[edge] (0,0) -- (0.5,1);
			\draw[edge] (0,0) -- (-0.5,1);
			\begin{pgfonlayer}{background}
				\draw[zhyedge1] (0,0) -- (-0,0);
				\draw[zhyedge2] (0,0) -- (-0,0);
				\draw[hyedge1, color=red] (-0.5,1) -- (0.5,1);
				\draw[hyedge2] (-0.5,1) -- (0.5,1);
			\end{pgfonlayer}
		\end{tikzpicture}
		,
		\begin{tikzpicture}
			\node[vertex, label=right:{\footnotesize 3}] at (0,0) {}; 
			\node[vertex, label=right:{\footnotesize 2}] at (0.5,1) {}; 
			\node[vertex, label=right:{\footnotesize 1}] at (-0.5,1) {}; 
			\node[vertex, label=right:{\footnotesize 1}] at (0.5,2) {}; 
			\draw[edge] (0,0) -- (0.5,1);
			\draw[edge] (0,0) -- (-0.5,1);
			\draw[edge] (0.5,1) -- (0.5,2);
			\begin{pgfonlayer}{background}
				\draw[hyedge1, color=red] (0.5,2) -- (0.5,1) -- (0,0) -- (-0.5,1);
				\draw[hyedge2] (0.5,2) -- (0.5,1) -- (0,0) -- (-0.5,1);
			\end{pgfonlayer}
		\end{tikzpicture}
		,
		\begin{tikzpicture}
			\node[vertex, label=right:{\footnotesize 1}] at (0,0) {}; 
			\node[vertex, label=right:{\footnotesize 2}] at (1,1) {}; 
			\node[vertex, label=right:{\footnotesize 2}] at (-1,1) {}; 
			\node[vertex, label=right:{\footnotesize 1}] at (0,1) {}; 
			\node[vertex, label=right:{\footnotesize 1}] at (1,2) {}; 
			\draw[edge] (0,0) -- (0,1);
			\draw[edge] (0,0) -- (-1,1);
			\draw[edge] (0,0) -- (1,1);
			\draw[edge] (1,1) -- (1,2);
			\begin{pgfonlayer}{background}
				\draw[zhyedge1] (0,0) -- (0,0);
				\draw[zhyedge2] (0,0) -- (0,0);
				\draw[hyedge1, color=red] (-1,1) -- (-1,1);
				\draw[hyedge2] (-1,1) -- (-1,1);
				\draw[hyedge1, color=blue] (0,1) -- (0,1);
				\draw[hyedge2] (0,1) -- (0,1);
				\draw[hyedge1, color=green] (1,1) -- (1,2);
				\draw[hyedge2] (1,1) -- (1,2);
			\end{pgfonlayer}
		\end{tikzpicture}. 
		$$
		Here, the symbols \begin{tikzpicture} \node[vertex, label=right:{\footnotesize i}] at (0,0) {}; \end{tikzpicture} represent nodes, elements of the set $\scN$ which are mapped by $\scL$ to the element $i$. Black lines between nodes are the elements of the set $\scE$. The partial ordering $<$ obtained from $\scE$ is graphically represented with the minimal element at the bottom of the diagram. 
		
		As with Example \ref{example:3:2}, the highlighted regions identify the hyperedges of the tree. The hyperedge $h_0$ is the subset that is colored gray. Other hyperedges are coloured with distinct colours representing differing hyperedges, so that $H$ is the collection of all coloured hyperedges. 
	\end{example}
	
	\begin{example}
		To visualise the operation $\circledast$, we have
		$$
		
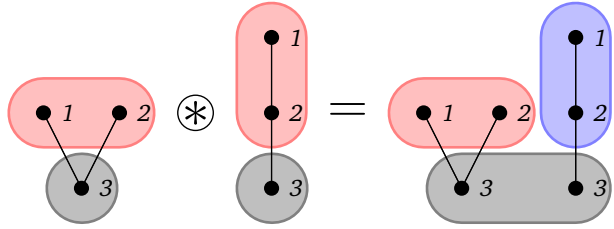
\begin{tikzpicture}
			\node[vertex, label=right:{\footnotesize 3}] at (0,0) {}; 
			\node[vertex, label=right:{\footnotesize 2}] at (0.5,1) {}; 
			\node[vertex, label=right:{\footnotesize 1}] at (-0.5,1) {}; 
			\node at (1.5, 1) {{\huge $\circledast$}};
			\node[vertex, label=right:{\footnotesize 3}] at (2.5,0) {}; 
			\node[vertex, label=right:{\footnotesize 2}] at (2.5,1) {}; 
			\node[vertex, label=right:{\footnotesize 1}] at (2.5,2) {}; 
			\node at (3.5, 1) {{\huge $=$}};
			\node[vertex, label=right:{\footnotesize 3}] at (5,0) {}; 
			\node[vertex, label=right:{\footnotesize 1}] at (4.5,1) {}; 
			\node[vertex, label=right:{\footnotesize 2}] at (5.5,1) {}; 
			\node[vertex, label=right:{\footnotesize 3}] at (6.5,0) {}; 
			\node[vertex, label=right:{\footnotesize 2}] at (6.5,1) {}; 
			\node[vertex, label=right:{\footnotesize 1}] at (6.5,2) {}; 
			\draw[edge] (0,0) -- (0.5,1);
			\draw[edge] (0,0) -- (-0.5,1);
			\draw[edge] (2.5,0) -- (2.5,1);
			\draw[edge] (2.5,1) -- (2.5,2);
			\draw[edge] (5,0) -- (4.5,1);
			\draw[edge] (5,0) -- (5.5,1);
			\draw[edge] (6.5,0) -- (6.5,1);
			\draw[edge] (6.5,1) -- (6.5,2);
			\begin{pgfonlayer}{background}
				\draw[zhyedge1] (0,0) -- (-0,0);
				\draw[zhyedge2] (0,0) -- (-0,0);
				\draw[hyedge1, color=red] (-0.5,1) -- (0.5,1);
				\draw[hyedge2] (-0.5,1) -- (0.5,1);
				\draw[zhyedge1] (2.5,0) -- (2.5,0);
				\draw[zhyedge2] (2.5,0) -- (2.5,0);
				\draw[hyedge1, color=red] (2.5,1) -- (2.5,2);
				\draw[hyedge2] (2.5,1) -- (2.5,2);
				\draw[zhyedge1] (5,0) -- (6.5,0);
				\draw[zhyedge2] (5,0) -- (6.5,0);
				\draw[hyedge1, color=red] (4.5,1) -- (5.5,1);
				\draw[hyedge2] (4.5,1) -- (5.5,1);
				\draw[hyedge1, color=blue] (6.5,1) -- (6.5,2);
				\draw[hyedge2] (6.5,1) -- (6.5,2);
			\end{pgfonlayer}
		\end{tikzpicture}
		$$
		Thus we see that $\circledast$ takes the union of sets $h_0^{T_1}$ and $h_0^{T_2}$, but the sets $H^{T_1}$ and $H^{T_2}$ are unioned although their elements remain distinct. 
	\end{example}
	
	\begin{example}
		The operation $\cE$ maps the hyperedge $h_0$ into the set $H$. Thus $\cE$ makes any gray hyperedges coloured. Therefore
		$$
		
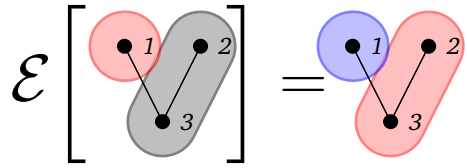
\begin{tikzpicture}
			\node[vertex, label=right:{\footnotesize 3}] at (0,-0.5) {}; 
			\node[vertex, label=right:{\footnotesize 2}] at (0.5,0.5) {}; 
			\node[vertex, label=right:{\footnotesize 1}] at (-0.5,0.5) {}; 
			\node at (-1.5,0) {{\Huge$\cE\bigg[$}}; 
			\node at (1.5,0) {{\Huge$\bigg]=$}}; 
			\node[vertex, label=right:{\footnotesize 3}] at (3,-0.5) {}; 
			\node[vertex, label=right:{\footnotesize 2}] at (3.5,0.5) {}; 
			\node[vertex, label=right:{\footnotesize 1}] at (2.5,0.5) {}; 
			\draw[edge] (0,-0.5) -- (0.5,0.5);
			\draw[edge] (0,-0.5) -- (-0.5,0.5);
			\draw[edge] (3,-0.5) -- (3.5,0.5);
			\draw[edge] (3,-0.5) -- (2.5,0.5);
			\begin{pgfonlayer}{background}
				\draw[zhyedge1] (0,-0.5) -- (0.5,0.5);
				\draw[zhyedge2] (0,-0.5) -- (0.5,0.5);
				\draw[hyedge1, color=red] (-0.5,0.5) -- (-0.5,0.5);
				\draw[hyedge2] (-0.5,0.5) -- (-0.5,0.5);
				\draw[hyedge1, color=red] (3,-0.5) -- (3.5,0.5);
				\draw[hyedge2] (3,-0.5) -- (3.5,0.5);
				\draw[hyedge1, color=blue] (2.5,0.5) -- (2.5,0.5);
				\draw[hyedge2] (2.5,0.5) -- (2.5,0.5);
			\end{pgfonlayer}
		\end{tikzpicture}
		$$
		However, when there are no gray hyperedges, the operator $\cE$ is just the identity so that
		$$
		
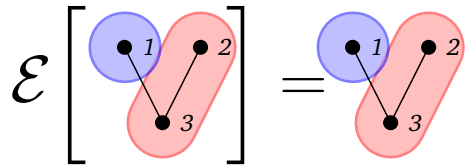
\begin{tikzpicture}
			\node[vertex, label=right:{\footnotesize 3}] at (0,-0.5) {}; 
			\node[vertex, label=right:{\footnotesize 2}] at (0.5,0.5) {}; 
			\node[vertex, label=right:{\footnotesize 1}] at (-0.5,0.5) {}; 
			\node at (-1.5,0) {{\Huge$\cE\bigg[$}}; 
			\node at (1.5,0) {{\Huge$\bigg]=$}}; 
			\node[vertex, label=right:{\footnotesize 3}] at (3,-0.5) {}; 
			\node[vertex, label=right:{\footnotesize 2}] at (3.5,0.5) {}; 
			\node[vertex, label=right:{\footnotesize 1}] at (2.5,0.5) {}; 
			\draw[edge] (0,-0.5) -- (0.5,0.5);
			\draw[edge] (0,-0.5) -- (-0.5,0.5);
			\draw[edge] (3,-0.5) -- (3.5,0.5);
			\draw[edge] (3,-0.5) -- (2.5,0.5);
			\begin{pgfonlayer}{background}
				\draw[hyedge1, color=red] (0,-0.5) -- (0.5,0.5);
				\draw[hyedge2] (0,-0.5) -- (0.5,0.5);
				\draw[hyedge1, color=blue] (-0.5,0.5) -- (-0.5,0.5);
				\draw[hyedge2] (-0.5,0.5) -- (-0.5,0.5);
				\draw[hyedge1, color=red] (3,-0.5) -- (3.5,0.5);
				\draw[hyedge2] (3,-0.5) -- (3.5,0.5);
				\draw[hyedge1, color=blue] (2.5,0.5) -- (2.5,0.5);
				\draw[hyedge2] (2.5,0.5) -- (2.5,0.5);
			\end{pgfonlayer}
		\end{tikzpicture}
		$$
	\end{example}
	
	\begin{example}
		The operation $\lfloor \cdot \rfloor$ acts in the same way as the rooting operation for the collection of forests in the Butcher-Connes-Kreimer setting, with the new node being added to the hyperedge $h_0$. Thus
		$$
		
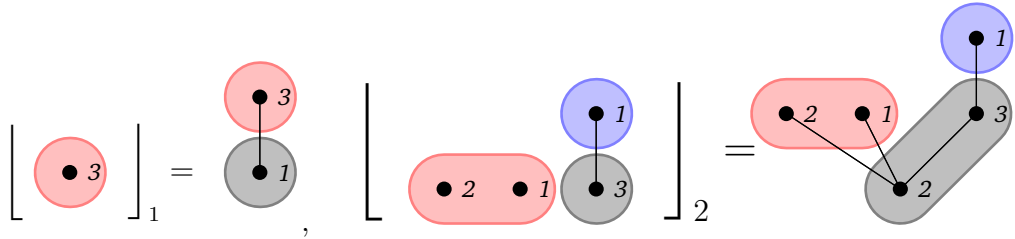
\begin{tikzpicture}
			\node[vertex, label=right:{\footnotesize 3}] at (0,0) {}; 
			\node at (-0.75,0) {{\large$\Bigg\lfloor$}}; 
			\node at (1.2,0) {{\large$\Bigg\rfloor_1 =$}}; 
			\node[vertex, label=right:{\footnotesize 1}] at (2.5,0) {}; 
			\node[vertex, label=right:{\footnotesize 3}] at (2.5,1) {}; 
			\draw[edge] (2.5,0) -- (2.5,1);
			\begin{pgfonlayer}{background}
				\draw[hyedge1, color=red] (0,0) -- (-0,0);
				\draw[hyedge2] (0,0) -- (-0,0);
				\draw[hyedge1, color=red] (2.5,1) -- (2.5,1);
				\draw[hyedge2] (2.5,1) -- (2.5,1);
				\draw[zhyedge1] (2.5,0) -- (2.5,0);
				\draw[zhyedge2] (2.5,0) -- (2.5,0);
			\end{pgfonlayer}
		\end{tikzpicture}
		, \quad
		\begin{tikzpicture}
			\node[vertex, label=right:{\footnotesize 2}] at (0,0) {}; 
			\node[vertex, label=right:{\footnotesize 1}] at (1,0) {}; 
			\node[vertex, label=right:{\footnotesize 3}] at (2,0) {}; 
			\node[vertex, label=right:{\footnotesize 1}] at (2,1) {}; 
			\node at (-1, 0.5) {{\LARGE $\Bigg\lfloor$}};
			\node at (3.5, 0.5) {{\LARGE $\Bigg\rfloor_{2} =$}};
			\node[vertex, label=right:{\footnotesize 2}] at (6,0) {}; 
			\node[vertex, label=right:{\footnotesize 2}] at (4.5,1) {}; 
			\node[vertex, label=right:{\footnotesize 1}] at (5.5,1) {}; 
			\node[vertex, label=right:{\footnotesize 3}] at (7,1) {}; 
			\node[vertex, label=right:{\footnotesize 1}] at (7,2) {}; 
			\draw[edge] (2,0) -- (2,1);
			\draw[edge] (6,0) -- (4.5,1);
			\draw[edge] (6,0) -- (5.5,1);
			\draw[edge] (6,0) -- (7,1);
			\draw[edge] (7,2) -- (7,1);
			\begin{pgfonlayer}{background}
				\draw[zhyedge1] (2,0) -- (2,0);
				\draw[zhyedge2] (2,0) -- (2,0);
				\draw[hyedge1, color=red] (0,0) -- (1,0);
				\draw[hyedge2] (0,0) -- (1,0);
				\draw[hyedge1, color=blue] (2,1) -- (2,1);
				\draw[hyedge2] (2,1) -- (2,1);
				\draw[zhyedge1] (6,0) -- (7,1);
				\draw[zhyedge2] (6,0) -- (7,1);
				\draw[hyedge1, color=red] (4.5,1) -- (5.5,1);
				\draw[hyedge2] (4.5,1) -- (5.5,1);
				\draw[hyedge1, color=blue] (7,2) -- (7,2);
				\draw[hyedge2] (7,2) -- (7,2);
			\end{pgfonlayer}
		\end{tikzpicture}
		$$
	\end{example}

	\subsection{Module of Lions forests}
	\label{subse:module:lions:forests}
	
	In the standard theory of rough paths, one defines the Butcher-Connes-Kreimer algebra over the ring $(\bR^e, +, \centerdot)$ (where $\centerdot: \bR^e \times \bR^e \to \bR^e$ is coordinatewise product) as
	\begin{equation*}
		\cH_d:= \bigoplus_{\tau \in \fF_{0, d}} \bR^e. 
	\end{equation*}
	This can equivalently be written as
	\begin{equation}
		\label{eq:Connes-Kreimer}
		\cH_d:= \bigoplus_{\tau \in \fF_0} \lin\big( (\bR^d)^{\otimes |\scN^\tau|}, \bR^e\big)
	\end{equation}
	and throughout we will typically use the form of Equation \eqref{eq:Connes-Kreimer}. The operation $\odot$ defined in Equation \eqref{eq:Tree_product} can be extended linearly to the space $\cH_d$ to obtain $\odot:\cH_d \times \cH_d \to \cH_d$ with unit $\rId$ (the empty forest) making $(\cH_d, \odot, \rId)$ an associative, commutative unital algebra over the ring $(\bR^e, + \centerdot)$. 
	
	One would like to have (up to some extent) a similar algebraic structure, but for Lions trees. The main difficulty in this regard is that, in the end, hyperedges are equipped with random labels. For this reason, we consider next two probability spaces $(\Omega, \cF, \bP)$ and $(\Omega', \cF', \bP')$. As explained in the Appendix \ref{subsection:Notation}, for any $m,n\in \bN$ and $\rp = (p_1, ..., p_m)$ where $p_i\in (1, \infty)$ we have that
	\begin{equation*}
		L^0\bigg(\Omega, \bP; L^{\rp}\Big( (\Omega')^{\times m}, (\bP')^{\times m}; \lin\big( (\bR^d)^{\otimes n}, \bR^e\big) \Big) \bigg)
	\end{equation*}
	is a module with respect to the ring $\big( L^0(\Omega, \bP; \bR^e), +. \centerdot \big)$. When $m=0$, we have
	\begin{equation*}
		L^0\bigg(\Omega, \bP; L^{\rp}\Big( (\Omega')^{\times 0}, (\bP')^{\times 0}; \lin\big( (\bR^d)^{\otimes n}, \bR^e\big) \Big) \bigg) = L^0\Big( \Omega, \bP; \lin\big( (\bR^d)^{\otimes n}, \bR^e\big) \Big). 
	\end{equation*}
	The intuition for distinguishing between the tagged and detagged probability spaces is important for the theory this work exposes but not specifically for this work: In many distributional dependent dynamics (but critically not McKean-Vlasov type dynamics) the probability space on which the random variable representing the lift of the distribution exists need not be the same as the probability space on which the solution exists. As a simple motivational example, consider McKean-Vlasov dynamics with a common noise. Thus we write $(\Omega, \cF, \bP)$ to represent the tagged probability space (on which the solution exists) and $(\Omega', \cF', \bP')$ the detagged probability space on which the lift of the law exists. 
	
	In the next definition, we let this module act on Lions trees equipped with $n$ nodes and $m$ hyperedges (excluding the $h_{0}$-hyperedge). Intuitively, this action may be explained as follows: If we had to follow the construction of the Butcher-Connes-Kreimer algebra (see for instance Remark \ref{rem:CK:algebra} below), we should equip all the hyperedges of the trees with indices. The $h_{0}$-hyperedge, which carries the tagged particle, should be equipped with $\omega_{0} \in \Omega$ to emphasize the realisation of its underlying idiosyncratic noise. 
	
	Similarly, the $m$ remaining hyperedges should be equipped with labels 
	$(\omega_{1}, ..., \omega_{m}) \in \Omega^{\times m}$. However, since trees are not equipped with such kinds of labels, we must restore the occurrences of $\omega_{0},\omega_{1}, ..., \omega_{m}$ in the module structure of the mean-field analogue of the Butcher-Connes-Kreimer algebra. This is the rationale for the following definition:
	
	\begin{definition}
		\label{definition:AlgebraRVs}
		Let $(\Omega, \cF, \bP)$ and $(\Omega', \cF', \bP')$ be probability spaces. For each $T\in \scF$, let
		\begin{equation}
			\label{eq:definition:AlgebraRVs-int}
			p[T] = (p[T]_h)_{h\in H^T} \in (1, \infty)^{\times |H^T|} 
			\quad \mbox{and} \quad
			\rp = \big( p[T] \big)_{T\in \scF} \in \bigsqcup_{T\in \scF} (1, \infty)^{\times |H^T|}. 
		\end{equation}
		We say that $\rp = \big( (p_h)_{h\in H^T} \big)_{T\in \scF}$ is an integrability functional. 
		
		We define $\scH_{\rp}'(\Omega')$ to be the $\bR^e$-module 
		\begin{equation*}
			\scH_{\rp}'(\Omega') = \bigoplus_{T\in \scF_0} L^{p[T]} \Big( (\Omega')^{\times |H^T|}, (\bP')^{\times |H^T|}; \lin \big( (\bR^d)^{\otimes |\scN^T|}, \bR^e\big) \Big). 
		\end{equation*}
		
		We also define $\scH_{\rp}(\Omega, \Omega')$ to be the $L^0(\Omega, \bP; \bR^e)$-module 
		\begin{align}
			\nonumber
			\scH_\rp(\Omega, \Omega') &= L^0\big( \Omega, \bP; \scH_{\rp}'(\Omega') \big)
			\\
			\label{eq:definition:AlgebraRVs}
			&= L^0\bigg( \Omega, \bP; \bigoplus_{T\in \scF_0} L^{p[T]} \Big( (\Omega')^{\times |H^T|}, (\bP')^{\times |H^T|}; \lin \big( (\bR^d)^{\otimes |\scN^T|}, \bR^e\big) \Big) \bigg).  
		\end{align}
		In the same spirit as Equation \eqref{eq:Connes-Kreimer}, this can be written as
		\begin{equation*}
			\scH_{\rp}(\Omega, \Omega') = L^0\bigg( \Omega, \bP; \bigoplus_{T\in \scF_{0,d}} L^{p[T]} \Big( (\Omega')^{\times |H^T|}, (\bP')^{\times |H^T|}; \bR^e \Big) \bigg) 
		\end{equation*}
		so we see that the inclusion of labellings corresponds coordinate evaluation in \eqref{eq:definition:AlgebraRVs}. 
	\end{definition}

	When there is no ambiguity over the choice of probability space, we will often simply write $\scH_{\rp}$. We will use the convention that $\scH_p(\Omega, \Omega')$ has a representation in terms of the set $\scF_0$: for $X\in \scH(\Omega, \Omega')$, we write
	\begin{equation*}
		X(\omega_0, \cdot) = \sum_{T\in \scF_0} \Big\langle X, T \Big\rangle(\omega_0, \cdot) 
	\end{equation*}
	where
	\begin{equation}
		\label{eq:definition:AlgebraRVs:<X,T>}
		\begin{aligned}
			&\Big\langle X, T \Big\rangle \in L^{0}\bigg( \Omega, \bP; L^{p[T]} \Big( (\Omega')^{\times |H^T|}, (\bP')^{\times |H|^T}; \lin\big( (\bR^d)^{\otimes |\scN^T|}, \bR^e\big)  \Big) \bigg) 
			\quad\mbox{or}
			\\
			&\Big\langle X, T \Big\rangle(\omega_0, \cdot) \in L^{p[T]} \Big( (\Omega')^{\times |H^T|}, (\bP')^{\times |H^T|}; \lin\big( (\bR^d)^{\otimes |\scN^T|}, \bR^e \big) \Big). 
		\end{aligned}
	\end{equation}
	
	Now suppose that $\rp$ additionally satisfies that 
	\begin{equation}
		\label{eq:definition:AlgebraRVs-int+}
		\forall T_1, T_2 \in \scF
		\quad \mbox{we have that} \quad
		p\big[ T_1 \circledast T_2 \big] = \big( p[T_1], p[T_2] \big)
	\end{equation}

	\begin{definition}
		Let $(\Omega, \cF, \bP)$ and $(\Omega', \cF', \bP')$ be probability spaces and let $\rp = \big( (p_h)_{h\in H^T} \big)_{T\in \scF}$ be an integrability functional. 
		
		By extending $\circledast$ to be bilinear, the pair of triples $\big( \scH_{\rp}'(\Omega'), \circledast, \rId \big)$ and $\big( \scH_{\rp}(\Omega, \Omega'), \circledast, \rId \big)$ becomes commutative algebras over the respective rings $\bR^e$ and $L^0(\Omega, \bP; \bR^e)$. In the former, this is because
		\begin{equation}
			\label{eq:definition:AlgebraRVs:Product-}
			X_1 \circledast X_2 
			= 
			\sum_{T\in \scF_0} \Big\langle X_1 \circledast X_2, T \Big\rangle
			= 
			\sum_{T\in \scF_0} \bigg( \sum_{\substack{T_1, T_2\in \scF_0\\ T_1\circledast T_2 = T }} \Big\langle X_1, T_1 \Big\rangle(\cdot) \otimes \Big\langle X_2, T_2 \Big\rangle(\cdot) \bigg) 
		\end{equation}
		where
		\begin{align*}
			\Big\langle X_1&, T_1\Big\rangle (\omega'_{H^{T_1}}) \otimes \Big\langle X_2, T_2\Big\rangle ( \omega'_{H^{T_2}}) 
			\\
			\in& L^{p[T_1]}\Big( (\Omega')^{\times |H^{T_1}|}, (\bP')^{\times |H^{T_1}|}; \lin\big( (\bR^d)^{\otimes |\scN^{T_1}|}, \bR^e\big) \Big)
			\\
			&\quad \otimes L^{p[T_2]}\Big( (\Omega')^{\times |H^{T_2}|}, (\bP')^{\times |H^{T_2}|}; \lin\big( (\bR^d)^{\otimes |\scN^{T_2}|}, \bR^e\big) \Big) 
			\\
			=& L^{p[T_1], p[T_2]}\Big( (\Omega')^{\times |H^{T_1}|} \times (\Omega')^{\times |H^{T_2}|}, (\bP')^{\times |H^{T_1}|} \times (\bP')^{\times |H^{T_2}|}; \lin\big( (\bR^d)^{\otimes (|\scN^{T_1}|+ |\scN^{T_2}|)}, \bR^e\big) \Big) 
			\\
			=& L^{p[T_1\circledast T_2]}\Big( (\Omega')^{\times |H^{T_1\circledast T_2}|}, (\bP')^{\times |H^{T_1\circledast T_2}|} ; \lin\big( (\bR^d)^{\otimes |\scN^{T_1\circledast T_2}|}, \bR^e\big) \Big) 
			\\
			=& L^{p[T]}\Big( (\Omega')^{\times |H^{T}|}, (\bP')^{\times |H^{T}|} ; \lin\big( (\bR^d)^{\otimes |\scN^{T}|}, \bR^e\big) \Big). 
		\end{align*}
		In the latter, 
		\begin{align}
			\nonumber
			X_1 \circledast X_2(\omega_0, \cdot) &= \sum_{T\in \scF_0} \Big\langle X_1 \circledast X_2, T \Big\rangle(\omega_0, \cdot) 
			\\
			\label{eq:definition:AlgebraRVs:Product}
			&= \sum_{T\in \scF_0} \bigg( \sum_{\substack{T_1, T_2\in \scF_0\\ T_1\circledast T_2 = T }} \Big\langle X_1, T_1 \Big\rangle(\omega_0, \cdot) \otimes \Big\langle X_2, T_2 \Big\rangle(\omega_0, \cdot) \bigg) 
		\end{align}
		where
		\begin{align*}
			\Big\langle X_1&, T_1\Big\rangle (\omega_0, \omega'_{H^{T_1}}) \otimes \Big\langle X_2, T_2\Big\rangle (\omega_0, \omega'_{H^{T_2}}) 
			\\
			\in& L^0\bigg( \Omega, \bP; L^{p[T_1]}\Big( (\Omega')^{\times |H^{T_1}|}, (\bP')^{\times |H^{T_1}|}; \lin\big( (\bR^d)^{\otimes |\scN^{T_1}|}, \bR^e\big) \Big)
			\\
			&\quad \otimes L^{p[T_2]}\Big( (\Omega')^{\times |H^{T_2}|}, (\bP')^{\times |H^{T_2}|}; \lin\big( (\bR^d)^{\otimes |\scN^{T_2}|}, \bR^e\big) \Big) \bigg)
			\\
			=& L^0\bigg( \Omega, \bP; L^{p[T]}\Big( (\Omega')^{\times |H^{T}|}, (\bP')^{\times |H^{T}|} ; \lin\big( (\bR^d)^{\otimes |\scN^{T}|}, \bR^e\big) \Big) \bigg). 
		\end{align*}
	\end{definition}
	
	We use the convention that $(\omega'_{H^{T}}) \in (\Omega')^{\times |H^{T}|}$ where $H$ is the set of hyperedges of the tree $T$ written as $H^T = \{h^T, ...\}$. In this work, it will be important to distinguish free variables but not to order them as is common convention. 
	
	Pay close attention in \eqref{eq:definition:AlgebraRVs:<X,T>}, even when the set $h_0^T$ is empty, the random variable $\langle X, T\rangle$ will still be a map from $\Omega \times (\Omega')^{\times |H^T|}$. 
	\begin{remark}
		\label{rem:CK:algebra}
		The Butcher-Connes-Kreimer Hopf algebra $\cH_d$ (see Equation \eqref{eq:Connes-Kreimer}) shares some aesthetic similarities with Equation \eqref{eq:definition:AlgebraRVs}. The Butcher-Connes-Kreimer is a module over the ring $\bR^e$ much like $\scH'(\Omega')$ whereas $\scH(\Omega, \Omega')$ is a module over the ring $L^\infty(\Omega, \bP; \bR^e)$. In the Butcher-Connes-Kreimer algebra, each tree has a labelling which corresponds to the components of the driving noise associated to each term and the module in front of $\tau$ is implicitly required to depend on the labels that $\tau$ carries. 
		
		Therefore, one might naively guess that when considering probabilistic dynamics, one would obtain some sort of space with labellings taking values on $\bN$ where each independent driving signal has its own labelling. This is not the case. The hypergraphic structure of Lions trees is analogous to this, although the restrictions imposed on the hyperedges in Definition \ref{definition:Forests} do emphasise that there is more going on here. Each hyperedge identifies the collection of nodes that correspond to increments of the driving signals that share a common source of noise. The $0$-hyperedge corresponds to the \emph{tagged} probability space within which the solution is identified (which can be identified with $L^\infty(\Omega, \bP; \bR^e)$). 
	\end{remark}

	\begin{proposition}
		\label{proposition:associativity}
		Let $(\Omega, \cF,\bP)$ and $(\Omega', \cF', \bP')$ be probability spaces and let $\rp$ be an integrability functional that satisfies Equation \eqref{eq:definition:AlgebraRVs-int+}. Let $(\bR^e, +, \centerdot)$ and $\big( L^0(\Omega, \bP; \bR^e), +, \centerdot \big)$ be rings and let $\scH_\rp'(\Omega')$, $\scH_\rp(\Omega, \Omega')$ and $\circledast$ be as defined in Definition \ref{definition:AlgebraRVs}. 
		
		Then $\big(\scH_\rp'(\Omega'), \circledast , \rId \big)$ and $\big( \scH_\rp(\Omega, \Omega'), \circledast , \rId \big)$ are associative and commutative algebras over the rings $(\bR^e, +, \centerdot)$ and $\big( L^0(\Omega, \bP; \bR^e), +, \centerdot \big)$ respectively. 
	\end{proposition}
	The proof of Proposition \ref{proposition:associativity} is found in \cite{salkeld2022LionsTrees}. 
	
	\subsubsection{Quotient}
	
	To understand the link between labelled directed trees and Lions trees, we first direct the reader to Proposition \ref{proposition:classicTay<=>LionsTay*}: When we take iterative spacial derivatives of a function of the form 
	\begin{equation*}
		f\Big( x^{i, N}, \tfrac{1}{N} \sum_{j=1}^N \delta_{x^{j, N}} \Big)	
	\end{equation*}
	with respect to a multi-index (for a detailed proof see \cite{salkeld2022Lions}) we get a collection of terms described by all the partitions that are finer than the partition obtained from the multi-index: 
	\begin{example}
		\label{example:Sums-finer-partitions}
		We highlight Equation \eqref{eq:example:Sums-finer-partitions} here to emphasise the following point: when taking a directional derivative of some function of an empirical measure, one obtains a summation of Lions derivatives indexed by all the partition sequences \emph{that are finer} than the partition sequence obtained from the pre-image of the multi-index from the directional derivative. If we can expect similar summations, we would like to exploit the natural symmetry that arises from terms having this form. 
		
		For instance, let $f:\bR^e \times \cP_2(\bR^e) \to \lin(\bR^d, \bR^e)$ and if we consider the Butcher-Connes-Kreimer expansion of a function $\boldsymbol{f}: (\bR^e)^{\oplus N} \to \big( \lin(\bR^d, \bR^e) \big)^{\oplus N}$ of the form
		\begin{equation*}
			\boldsymbol{f}\big( \boldsymbol{x} \big) = \boldsymbol{f}\Big( (x^{1, N}, ..., x^{N, N}) \Big) = \bigoplus_{i=1}^N f\Big( x^{i, N}, \mu_N\big[ \boldsymbol{x} \big] \Big)
		\end{equation*}
		and for $i, j, k \in \{1, ..., N\}$ look at the elementary differentials associated to the tree
		\begin{equation*}
			\begin{tikzpicture}
				\node[vertex, label=right:{\footnotesize i}] at (0,0) {};
				\node[vertex, label=right:{\footnotesize j}] at (0,1) {};
				\node[vertex, label=right:{\footnotesize k}] at (0,2) {};
				\draw[edge] (0,0) -- (0,1);
				\draw[edge] (0,2) -- (0,1);
				\node at (1,1) {$\sim$};
				\node at (6,1) {$\nabla_{j} f\Big( x^{i, N}, \mu_N\big[ \boldsymbol{x}\big] \Big) \cdot \nabla_{k} f\Big( x^{j, N}, \mu_N\big[ \boldsymbol{x} \big] \Big) \cdot f\Big( x^{k, N}, \mu_N\big[\boldsymbol{x} \big] \Big)$. };
			\end{tikzpicture}
		\end{equation*}
		Expanding this out using Equation \eqref{eq:example:Sums-finer-partitions} and summing over $(j, k) \in \{1, ..., N \}^{\times 2}$, we get that
		\begin{equation*}
			\begin{tikzpicture}
				\node at (-1, 1) {{\Huge$\sum$}};
				\node at (-1, 0) {$j, k=1$};
				\node at (-1, 2) {$N$};
				\node[vertex, label=right:{\footnotesize i}] at (0,0) {};
				\node[vertex, label=right:{\footnotesize j}] at (0,1) {};
				\node[vertex, label=right:{\footnotesize k}] at (0,2) {};
				\draw[edge] (0,0) -- (0,1);
				\draw[edge] (0,2) -- (0,1);
				\node at (1,1) {$\sim$};
				\node at (6,1) {$\partial_0 f\Big( x^{i, N}, \mu_N\big[ \boldsymbol{x}\big] \Big) \cdot \partial_0 f\Big( x^{i, N}, \mu_N\big[ \boldsymbol{x} \big] \Big) \cdot f\Big( x^{i, N}, \mu_N\big[\boldsymbol{x} \big] \Big)$};
				\node at (7,0) {$+\tfrac{1}{N}\sum_{j} \partial_0 f\Big( x^{i, N}, \mu_N\big[ \boldsymbol{x}\big] \Big) \cdot \partial_1 f\Big( x^{i, N}, \mu_N\big[ \boldsymbol{x} \big] \Big) \cdot f\Big( x^{j, N}, \mu_N\big[\boldsymbol{x} \big] \Big)$};
				\node at (7,-1) {$+\tfrac{1}{N}\sum_{j} \partial_1 f\Big( x^{i, N}, \mu_N\big[ \boldsymbol{x}\big] \Big) \cdot \partial_0 f\Big( x^{j, N}, \mu_N\big[ \boldsymbol{x} \big] \Big) \cdot f\Big( x^{j, N}, \mu_N\big[\boldsymbol{x} \big] \Big)$};
				\node at (7,-2) {$+\tfrac{1}{N^2}\sum_{j, k} \partial_1 f\Big( x^{i, N}, \mu_N\big[ \boldsymbol{x}\big] \Big) \cdot \partial_1 f\Big( x^{j, N}, \mu_N\big[ \boldsymbol{x} \big] \Big) \cdot f\Big( x^{k, N}, \mu_N\big[\boldsymbol{x} \big] \Big)$. };
			\end{tikzpicture}
		\end{equation*}
		Our strategy is to represent these summations in terms of the four hypergraphic structures that give rise to the Lions trees
		\begin{equation}
			\label{eq:example:Sums-finer-partitions(all)}
			\begin{tikzpicture}
				\node at (-1, 1) {{\Huge$\sum$}};
				\node at (-1, 0) {$j, k=1$};
				\node at (-1, 2) {$N$};
				\node[vertex, label=right:{\footnotesize i}] at (0,0) {};
				\node[vertex, label=right:{\footnotesize j}] at (0,1) {};
				\node[vertex, label=right:{\footnotesize k}] at (0,2) {};
				\draw[edge] (0,0) -- (0,1);
				\draw[edge] (0,2) -- (0,1);
				\node at (1,1) {$\Rightarrow$};
				\node[vertex] at (2,0) {};
				\node[vertex] at (2,1) {};
				\node[vertex] at (2,2) {};
				\draw[edge] (2,0) -- (2,1);
				\draw[edge] (2,2) -- (2,1);
				\node at (3,1) {$+$};
				\node[vertex] at (4,0) {};
				\node[vertex] at (4,1) {};
				\node[vertex] at (4,2) {};
				\draw[edge] (4,0) -- (4,1);
				\draw[edge] (4,2) -- (4,1);
				\node at (5,1) {$+$};
				\node[vertex] at (6,0) {};
				\node[vertex] at (6,1) {};
				\node[vertex] at (6,2) {};
				\draw[edge] (6,0) -- (6,1);
				\draw[edge] (6,2) -- (6,1);
				\node at (7,1) {$+$};
				\node[vertex] at (8,0) {};
				\node[vertex] at (8,1) {};
				\node[vertex] at (8,2) {};
				\draw[edge] (8,0) -- (8,1);
				\draw[edge] (8,2) -- (8,1);
				\begin{pgfonlayer}{background}
					\draw[zhyedge1] (2,0) -- (2, 2);
					\draw[zhyedge2] (2,0) -- (2, 2);
					\draw[zhyedge1] (4,0) -- (4, 1);
					\draw[zhyedge2] (4,0) -- (4, 1);
					\draw[zhyedge1] (6,0) -- (6, 0);
					\draw[zhyedge2] (6,0) -- (6, 0);
					\draw[zhyedge1] (8,0) -- (8, 0);
					\draw[zhyedge2] (8,0) -- (8, 0);
					\draw[hyedge1,color=red] (4,2) -- (4, 2);
					\draw[hyedge2] (4,2) -- (4, 2);
					\draw[hyedge1,color=red] (6,1) -- (6, 2);
					\draw[hyedge2] (6,1) -- (6, 2);
					\draw[hyedge1,color=red] (8,1) -- (8, 1);
					\draw[hyedge2] (8,1) -- (8, 1);
					\draw[hyedge1,color=blue] (8,2) -- (8, 2);
					\draw[hyedge2] (8,2) -- (8, 2);
				\end{pgfonlayer}
			\end{tikzpicture}
		\end{equation}
		where the \textcolor{red}{Red} hyperedges represent a summation over $j$ and the \textcolor{blue}{Blue} hyperedge represents a summation over $k$. 
		
		Suppose first that we consider the case when $j\neq i$ and $k \neq i, j$ so that one might make the naive observation
		\begin{equation*}
			\begin{tikzpicture}
				\node[vertex, label=right:{\footnotesize i}] at (0,0) {};
				\node[vertex, label=right:{\footnotesize j}] at (0,1) {};
				\node[vertex, label=right:{\footnotesize k}] at (0,2) {};
				\draw[edge] (0,0) -- (0,1);
				\draw[edge] (0,2) -- (0,1);
				\node at (1,1) {$\sim$};
				\node at (6.5,1) {$\frac{1}{N^2} \partial_{1} f\Big( x^{i, N}, \mu_N\big[ \boldsymbol{x}\big] \Big) \cdot \partial_1 f\Big( x^{j, N}, \mu_N\big[ \boldsymbol{x} \big] \Big) \cdot f\Big( x^{k, N}, \mu_N\big[\boldsymbol{x} \big] \Big)$};
				\node at (12,1) {$\in$};
				\node[vertex] at (13,0) {};
				\node[vertex] at (13,1) {};
				\node[vertex] at (13,2) {};
				\draw[edge] (13,0) -- (13,1);
				\draw[edge] (13,2) -- (13,1);
				\begin{pgfonlayer}{background}
					\draw[zhyedge1] (13,0) -- (13, 0);
					\draw[zhyedge2] (13,0) -- (13, 0);
					\draw[hyedge1,color=red] (13,1) -- (13, 1);
					\draw[hyedge2] (13,1) -- (13, 1);
					\draw[hyedge1,color=blue] (13,2) -- (13, 2);
					\draw[hyedge2] (13,2) -- (13, 2);
				\end{pgfonlayer}
			\end{tikzpicture}
		\end{equation*}
		that is the partition obtained by inverting the labelling yields the hypergraph (with the addition that the $0$-hyperedge is the pre-image of $i$). However this is not quite the case. Recalling Theorem \ref{theorem:EquivalenceTrees}, the Lions hypergraphic structure is equivalent to partitioning the tree locally around a node and all nodes directly above it. This means that we also have
		\begin{equation*}
			\begin{tikzpicture}
				\node[vertex, label=right:{\footnotesize i}] at (0,0) {};
				\node[vertex, label=right:{\footnotesize j}] at (0,1) {};
				\node[vertex, label=right:{\footnotesize i}] at (0,2) {};
				\draw[edge] (0,0) -- (0,1);
				\draw[edge] (0,2) -- (0,1);
				\node at (1,1) {$\sim$};
				\node at (6.5,1) {$\frac{1}{N^2} \partial_{1} f\Big( x^{i, N}, \mu_N\big[ \boldsymbol{x}\big] \Big) \cdot \partial_1 f\Big( x^{j, N}, \mu_N\big[ \boldsymbol{x} \big] \Big) \cdot f\Big( x^{i, N}, \mu_N\big[\boldsymbol{x} \big] \Big)$};
				\node at (12,1) {$\in$};
				\node[vertex] at (13,0) {};
				\node[vertex] at (13,1) {};
				\node[vertex] at (13,2) {};
				\draw[edge] (13,0) -- (13,1);
				\draw[edge] (13,2) -- (13,1);
				\begin{pgfonlayer}{background}
					\draw[zhyedge1] (13,0) -- (13, 0);
					\draw[zhyedge2] (13,0) -- (13, 0);
					\draw[hyedge1,color=red] (13,1) -- (13, 1);
					\draw[hyedge2] (13,1) -- (13, 1);
					\draw[hyedge1,color=blue] (13,2) -- (13, 2);
					\draw[hyedge2] (13,2) -- (13, 2);
				\end{pgfonlayer}
			\end{tikzpicture}
		\end{equation*}
		since $i\neq j$ and $j\neq i$ so each node must be in a different hyperedge from the node directly above it. This is contrast to the (incorrect) hypothesis that 
		\begin{equation*}
			\begin{tikzpicture}
				\node[vertex, label=right:{\footnotesize i}] at (0,0) {};
				\node[vertex, label=right:{\footnotesize j}] at (0,1) {};
				\node[vertex, label=right:{\footnotesize i}] at (0,2) {};
				\draw[edge] (0,0) -- (0,1);
				\draw[edge] (0,2) -- (0,1);
				\node at (1,1) {$\in$};
				\node[vertex] at (2,0) {};
				\node[vertex] at (2,1) {};
				\node[vertex] at (2,2) {};
				\draw[edge] (2,0) -- (2,1);
				\draw[edge] (2,2) -- (2,1);
				\node at (5,1) {which is \underline{False}.};
				\begin{pgfonlayer}{background}
					\draw[zhyedge1] (2,0) -- (3, 0) -- (3,2) -- (2, 2);
					\draw[zhyedge2] (2,0) -- (3, 0) -- (3,2) -- (2, 2);
					\draw[hyedge1,color=red] (2,1) -- (2, 1);
					\draw[hyedge2] (2,1) -- (2, 1);
				\end{pgfonlayer}
			\end{tikzpicture}
		\end{equation*}
		Note also that this hypergraph is not admissible under Definition \ref{definition:Forests}. Indeed, Equation \eqref{eq:example:Sums-finer-partitions(all)} represents a summation over all possible Lions forests that have the same underlying directed forest. 
		
		A key term to obtaining this representation is when $i=j=k$ in which case
		\begin{equation*}
			\begin{tikzpicture}
				\node[vertex, label=right:{\footnotesize i}] at (0,0) {};
				\node[vertex, label=right:{\footnotesize i}] at (0,1) {};
				\node[vertex, label=right:{\footnotesize i}] at (0,2) {};
				\draw[edge] (0,0) -- (0,1);
				\draw[edge] (0,2) -- (0,1);
				\node at (1,1) {$\sim$};
				\node at (6,1) {$\partial_0 f\Big( x^{i, N}, \mu_N\big[ \boldsymbol{x}\big] \Big) \cdot \partial_0 f\Big( x^{i, N}, \mu_N\big[ \boldsymbol{x} \big] \Big) \cdot f\Big( x^{i, N}, \mu_N\big[\boldsymbol{x} \big] \Big)$};
				\node at (6.5,0) {$+\tfrac{1}{N}\partial_0 f\Big( x^{i, N}, \mu_N\big[ \boldsymbol{x}\big] \Big) \cdot \partial_1 f\Big( x^{i, N}, \mu_N\big[ \boldsymbol{x} \big] \Big) \cdot f\Big( x^{i, N}, \mu_N\big[\boldsymbol{x} \big] \Big)$};
				\node at (6.5,-1) {$+\tfrac{1}{N}\partial_1 f\Big( x^{i, N}, \mu_N\big[ \boldsymbol{x}\big] \Big) \cdot \partial_0 f\Big( x^{i, N}, \mu_N\big[ \boldsymbol{x} \big] \Big) \cdot f\Big( x^{i, N}, \mu_N\big[\boldsymbol{x} \big] \Big)$};
				\node at (6.5,-2) {$+\tfrac{1}{N^2}\partial_1 f\Big( x^{i, N}, \mu_N\big[ \boldsymbol{x}\big] \Big) \cdot \partial_1 f\Big( x^{i, N}, \mu_N\big[ \boldsymbol{x} \big] \Big) \cdot f\Big( x^{i, N}, \mu_N\big[\boldsymbol{x} \big] \Big)$};
			\end{tikzpicture}
		\end{equation*}
		which is making critical contributions to each of the summations associated to each Lions tree. If we simplify and  think of each of the relevant Lions forests as being identified with a pair of partition sequences of length 1 ($A_1[0] = \{0, 1\}$), we observe that this is a summation over pairs of partition sequences that are finer (with respect to the partial ordering defined in Equation \eqref{eq:partial-ordering}) than the pair of partition sequences obtained from the local pre-image of the labelling. 
	\end{example}
	
	With this is mind, we introduce the following definition:
	\begin{definition}
		\label{definition:Quotient-Partition}
		Let $\tau=(\scN, \scE)$ be a directed forest with preordering $\leq$. We define $\scQ^{\scE}(\scN)[0]$ to be the set of all pairs $(q_0, Q)$ such that $q_0 \subseteq \scN$, $Q\in \scP(\scN\backslash q_0)$ and
		\begin{enumerate}[label=(2.\arabic*)]
			\item If $q_0 \neq \emptyset$, then $q_0 \cap \fr[\tau] \neq \emptyset$. 
			\item $\forall q \in Q$ such that $x, y \in q$ and  $x<y$, $\exists z \in q$ such that $(y,z) \in \scE$. 
			\item $\forall q \in Q$ such that $x_1, y_1 \in q$ and $x_1 \leq \geq y_1$, $(x_1, x_2), (y_1, y_2) \in \scE$ and $x_2 \neq y_2$, then $x_2, y_2 \in q$. 
		\end{enumerate}
		
		For $(p_0, P), (q_0, Q) \in \scQ^E(\scN)[0]$, we say
		\begin{align*}
			(p_0, P) \subseteq (q_0, Q) \quad \iff \quad p_0 \subseteq q_0 \quad \mbox{and} \quad \forall p \in P, \exists q \in \big(Q \cup q_0 \big)\backslash\{\emptyset\} \quad \mbox{such that} \quad p \subseteq q. 
		\end{align*}
		We refer to $\scQ^{\scE}(\scN)[0]$ as the collection of tagged Lions partitions. 
	\end{definition}
	
	Although the exact meaning of this is explored in more detail in \cite{salkeld2022LionsTrees}, Definition \ref{definition:Quotient-Partition} provides a partial ordering for the set of all possible tagged hypergraphs of a directed forest $(\scN, \scE)$ that satisfy \ref{definition:Forests:2.1}, \ref{definition:Forests:2.2} and \ref{definition:Forests:2.3}. In words, the elements of $\scQ^{\scE}(\scN)[0]$ satisfy that $(p_0, P) \subseteq (q_0, Q)$ if the partition $\big( \{p_0\} \cup O \big)\backslash\{\emptyset\}$ is finer that the partition $\big( \{q_0\} \cup Q \big) \backslash \{\emptyset\}$. Further, for any $(h_0, H) \in \scQ^{\scE}(\scN)[0]$ we have that $(\scN, \scE, h_0, H) \in \scF$. 
	
	The reader should also keep in mind the link between partial orderings, Lions derivatives and Equation \eqref{eq:example:Sums-finer-partitions} which is also discussed in \cite{salkeld2022Lions}:
	
	\begin{lemma}
		\label{lemma:Link-infin_Syst-Part}
		Let $N\in \bN$, let $\Omega' = \{1, ..., N\}$, $\cF' = 2^{\{1, ..., N\}}$ and $\bP'$ be the uniform measure on $(\Omega', \cF')$. Let $i\in \Omega'$. For $\tau=(\scN, \scE, \scL, \hat{\scL}) \in \fF_{0, d\times N}$ where
		\begin{equation*}
			\scL:\scN \to \{1, ..., d\} \quad \mbox{and} \quad \hat{\scL}:\scN \to \{1, ..., N\} = \Omega', 
		\end{equation*}
		we denote
		\begin{equation}
			\label{eq:Label-to-TaggedPartition}
			\begin{split}
				\big( \hat{h}_0^{\tau}, \hat{H}^{\tau} \big) &= \bigg( \big( \hat{\scL}\big|_{\fr[\tau]} \big)^{-1}[i], \quad \Big\{ \big( \hat{\scL}\big|_{\fr[\tau]} \big)^{-1}[j]: j\in \Omega'\backslash \{i\} \Big\}\backslash \{\emptyset\} \bigg)\quad \mbox{and}  
				\\
				\fh[x] &= \bigg\{ \Big( \hat{\scL}\Big|_{\scN_x \cup \{x\} } \Big)^{-1} \big[ j \big]: j \in \{1, ..., N \} \bigg\}\backslash \{\emptyset\}. 
			\end{split} 
		\end{equation}
		Let $(h_0^{\tau}, H^{\tau} ) \in \scQ^\scE(\scN)[0]$ be the tagged partition such that
		\begin{equation}
			\label{eq:Label-to-TaggedPartition-}
			\big(\scN, \scE, h_0^{\tau}, H^{\tau}, \scL \big) = \Phi\Big[ \big(\scN, \scE, (\hat{h}_0^{\tau}, \hat{H}^{\tau}), \fh, \scL \big) \Big]
		\end{equation}
		where $\Phi: \hat{\scF}_0 \to \scF_0$ is the isomorphism from Theorem \ref{theorem:EquivalenceTrees}. 
		
		Let $\rp$ be an integrability functional. We define the $\bR^e$ module
		\begin{align}
			\nonumber
			\overline{\cH}_{d}:=&\bigg\{ X \in \cH_{d \times N}: \quad \forall \tau =( \scN, \scE, \scL, ) \in \fF_{0, d}, \quad \exists \Psi^{\tau}: \bigsqcup_{T\in \scF_{0, d}} (\Omega')^{\times |H^T|} \to \bR^e. 
			\\
			\label{eq:proposition:Link-infin_Syst-Part}
			&\Big\langle X, (\scN, \scE, \scL, \hat{\scL}) \Big\rangle = \sum_{\substack{(p_0, P)\in \scQ^{\scE}(\scN)[0] \\ (p_0, P) \subseteq (h_0^{\tau}, H^{\tau}) }} \tfrac{1}{N^{\times |P|}} \cdot \Psi^{\tau} \Big[ \big( \scN, \scE, p_0, P, \scL \big), \big( \hat{\scL}[p] \big)_{p\in P} \Big] \bigg\}. 
		\end{align}
		Then $\overline{\cH}_d \equiv \scH_{\rp}'(\Omega')$.  
	\end{lemma}
	The choice of $i\in \Omega'$ corresponds to the label of a tagged particle within a system of interacting particles and the mapping $\Psi$ maps all of the possible labelled trees with indices $\hat{\scL}$ to the state space. 
	
	Additionally, fix a probability space $(\Omega, \cF, \bP)$. We note that Lemma \ref{lemma:Link-infin_Syst-Part} can be extended to provide an embedding 
	\begin{equation*}
		L^0\Big( \Omega, \bP; \overline{\cH}_d \Big) \equiv \scH_{\rp}(\Omega, \Omega'). 
	\end{equation*}
	
	\begin{proof}
		First, we recall that for any $(\scN, \scE, \scL, \hat{\scL} ) \in \scF_{0, d\times N}$ and $\overline{X} \in \overline{\cH}_d$, the element
		\begin{equation*}
			\Big\langle \overline{X}, (\scN, \scE, \scL, \hat{\scL}) \Big\rangle \in \bR^e
		\end{equation*}
		so that the summation in Equation \eqref{eq:proposition:Link-infin_Syst-Part} runs over ring elements. Fix $\tau = (\scN, \scE, \scL) \in \fF_{0, d}$ and $\overline{X} \in \overline{\cH}_d$. Consider the summation
		\begin{equation*}
			\sum_{\hat{\scL}:\scN \to \Omega'} \Big\langle \overline{X}, (\scN, \scE, \scL, \hat{\scL}) \Big\rangle
			= 
			\sum_{\hat{\scL}:\scN \to \Omega'} \sum_{\substack{(p_0, P)\in \scQ^\scE(\scN)[0] \\ (p_0, P) \subseteq (h_0^{\tau}, H^{\tau}) }} \tfrac{1}{N^{\times |P|}} \cdot \Psi^{\tau}\Big[ \big( \scN, \scE, p_0, P, \scL \big), \big(\hat{\scL}[p]\big)_{p\in P} \Big]
		\end{equation*}
		where $(h_0^{\tau}, H^{\tau})$ was defined as Equation \eqref{eq:Label-to-TaggedPartition-} in terms of the choice of $\hat{\scL}$. Changing the order of the two summations, we get
		\begin{equation*}
			\sum_{\hat{\scL}:\scN \to \Omega'} \Big\langle \overline{X}, (\scN, \scE, \scL, \hat{\scL}) \Big\rangle= \sum_{(h_0, H) \in \scQ^\scE(\scN)[0] } \tfrac{1}{N^{\times |H|}} \sum_{i_{H} \in (\Omega')^{\times |H|}} \Psi^{\tau}\Big[ \big( \scN, \scE, h_0, H, \scL \big), i_{H} \Big]. 
		\end{equation*}
		Therefore, 
		\begin{align*}
			\overline{X} =& \sum_{(\scN, \scE, \scL) \in \fF_{0, d}} \sum_{\hat{\scL}:\scN \to \Omega'} \Big\langle \overline{X}, (\scN, \scE, \scL, \hat{\scL} ) \Big\rangle \in \bigoplus_{\tau \in \fF_{0, d\times N}} \bR^e
			\\
			&= \sum_{\substack{(\scN, \scE, \scL) \in \fF_{0, d} \\ (h_0, H) \in \scQ^{\scE}(\scN)[0]}} \tfrac{1}{N^{\times |H|}} \sum_{i_H \in (\Omega')^{\times |H|}} \Psi^{\tau}\Big[ \big( \scN , \scE, h_0, H, \scL_1\big), i_H \Big]
		\end{align*}
		so that for any choice of $T=(\scN, \scE, h_0, H, \scL) \in \scF_{0,d}$ we can define the measurable function
		\begin{equation*}
			\Big\langle X, T \Big\rangle: (\Omega')^{\times |H|} \to \bR^e 
			\quad
			\Big\langle X, T \Big\rangle(i_H) = \tfrac{1}{N^{\times |H|}} \cdot \Psi^{(\scN, \scE, \scL)}\Big[ \big( \scN , \scE, h_0, H, \scL \big), i_H \Big]. 
		\end{equation*}
		Then for any choice of integrability functional $\rp$ we have that
		\begin{align*}
			&\Big\langle X, T \Big\rangle \in L^{p[T]}\Big( (\Omega')^{\times |H^T|} , (\bP')^{\times |H^T|}; \bR^e \Big) \quad \mbox{so that}
			\\
			&X = \sum_{T\in \scF_{0, d}} \Big\langle X, T \Big\rangle \in \bigoplus_{T\in \scF_{0, d}} L^{p[T]}\Big( (\Omega')^{\times |H^T|} , (\bP')^{\times |H^T|}; \bR^e \Big)
		\end{align*}
		so that we have $\overline{\cH}_d \subseteq \scH_{\rp}'(\Omega')$. 
		
		On the other hand, suppose that $X \in \scH_{\rp}'(\Omega')$ so that for each $T = (\scN, \scE, h_0, H, \scL) \in \scF_{0,d}$ we have that
		\begin{equation*}
			\Big\langle X, T \Big\rangle \in L^{p[T]}\Big( (\Omega')^{\times |H^T|}, (\bP)^{\times |H^T|}; \bR^e \Big). 
		\end{equation*}
		Fix $T = (\scN, \scE, h_0, H, \scL) \in \scF_{0, d}$ and $\hat{\scL}: \scN \to \Omega'$. Then we define
		\begin{equation*}
			\Big\langle \overline{X}, (\scN, \scE, \scL, \hat{\scL}) \Big\rangle = \sum_{\substack{(p_0, P)\in \scQ^{\scE}(\scN)[0] \\ (p_0, P) \subseteq (h_0^{\tau}, H^{\tau}) }} \tfrac{1}{N^{\times |P|}} \cdot \Big\langle X, (\scN, \scE, p_0, P, \scL) \Big\rangle\big( (\hat{\scL}[p])_{p\in P} \big)
		\end{equation*}
		Then, by construction we have that
		\begin{equation*}
			\overline{X} = \sum_{(\scN, \scE, \scL, \hat{\scL}) \in \fF_{0, d}} \Big\langle \overline{X}, (\scN, \scE, \scL, \hat{\scL} ) \Big\rangle
		\end{equation*}
		satisfies that $X \in \overline{\cH}_d$ so that $\scH_{\rp}'(\Omega') \subseteq \overline{\cH}_d$. 
	\end{proof}
	
	\subsection{Products of distributions and coupled bialgebras}
	\label{subsec:CoupledBialgebra}
	
	Recall from \cite{connes1999hopf} that the Butcher-Connes-Kreimer coproduct $\triangle: \cH_d \to \cH_d \otimes \cH_d$ is the unique linear algebra homomorphism that satisfies the recursive relationship
	\begin{equation}
		\label{eq:francois:label:CK:co-product}
		\triangle \big[ \rId\big] =  \rId \otimes \rId
		\quad\mbox{and} \quad
		\triangle\Big[ \big\lfloor \tau \big\rfloor_i \Big] = \big\lfloor \tau \big\rfloor_i \times \rId + \Big( I \times \big\lfloor \cdot \big\rfloor_i \Big) \circ \triangle \Big[ \tau \Big]. 
	\end{equation}
	This is paired with the counit $\epsilon: \cH_d \to \bR^e$ defined by $\epsilon[X] = \langle X, \rId\rangle$, making $(\cH_d, \triangle, \epsilon)$ a co-associative coalgebra over the ring $(\bR^e, +, \centerdot)$.
	
	In this section, given two integrability functionals $\rp_1$ and $\rp_2$ we want to consider products of random of variables from
	\begin{equation*}
		\begin{split}
			&L^0 \bigg( \Omega, \bP; L^{\rp_1}\Big( (\Omega')^{\times |H^{1}|}, (\bP')^{\times |H^{1}|}; \lin \big( (\bR^d)^{\times |\scN^{1}|}, \bR^e\big) \Big) \bigg)
			\\
			\textrm{\rm and}
			\quad 
			&L^0 \bigg( \Omega, \bP; L^{\rp_2}\Big( (\Omega')^{\times |H^{2}|}, (\bP')^{\times |H^{2}|}; \lin \big( (\bR^d)^{\times |\scN^{2}|}, \bR^e\big) \Big) \bigg). 
		\end{split}
	\end{equation*}
	Recall that the former two spaces arise in the expansion \eqref{eq:definition:AlgebraRVs}. From a purely probabilistic point of view, the variables $(\omega'_{H^{1}})$ and $(\omega'_{H^{2}})$ are thought as free, which is the same as saying that we have two canonical representatives of two conditional laws given the realization $\omega_0$ and that we want to make their product. Of course, there are plenty of ways to do so, depending on the conditional joint law that should be above the two conditional laws. Equivalently, when defining the product of two elements of the above two spaces, we must recreate the correlations between the two corresponding conditional marginal laws. It turns out that, in our formalism, this amounts to specify which hyperedges in $T_{1}$ and $T_{2}$ are associated to the same free variable. For this, we found that the most parsimonious way was using a collection of operators over the space of partitions of a finite set.
	
	\subsubsection{Couplings}
	\label{subsubsection:Couplings}
	
	A coupling between two partitions $P$ and $Q$ is a way of describing a bijective mapping between a subset of $P$ and a subset of $Q$. Here and throughout, the elements of $P$ and $Q$ are an abstraction of probability spaces, and by coupling partitions we are analogously choosing a joint measure on the product of these probability spaces that is either the product measure (when there is no coupling) or the deterministic transport. 
	
	\begin{definition}
		Let $\scM$ and $\scN$ be two non-empty, disjoint, finite sets and let $P\in \scP(\scM)$ and $Q\in \scP(\scN)$ be partitions. We define
		$$
		P \tilde{\cup} Q = \Big\{ G \in \scP(\scM\cup \scN): \{ g\cap \scM: g\in G\}\backslash \emptyset = P, \{ g\cap \scN: g\in G\}\backslash \emptyset = Q\Big\}. 
		$$
		We refer to $P \tilde{\cup} Q$ as the set of coupled partitions. 
	\end{definition}
	
	$P \tilde{\cup} Q $ is the collection of partitions of the set $\scM\cup \scN$ that agree with the partition $P$ when restricting to $\scM$ and agrees with the partition $Q$ when restricting to $\scN$. As such, we have the symmetry relation
	$$
	P \tilde{\cup} Q = Q \tilde{\cup} P. 
	$$
	
	\begin{example}
		Suppose that we have the set $\scM=\{1, 2\}$ and $\scN=\{3, 4\}$ with partitions 
		$$
		P=\Big\{ \{1\}, \{2\}\Big\}\quad \mbox{ and }\quad Q=\Big\{ \{3,4\}\Big\}.
		$$ 
		Hence (see also Figure \ref{fig:coupling}) 
		$$
		P\tilde{\cup}Q = \Big\{ \big\{ \{1,3,4\}, \{2\} \big\}, \big\{ \{1\}, \{2,3,4\}\big\}, \big\{ \{1\}, \{2\}, \{3,4\} \big\} \Big\}.
		$$ 
		\begin{figure}[htb]
			\centering
			\begin{tikzpicture}
				\node[vertex, label=right:{\footnotesize 1}] at (0,1) {}; 
				\node[vertex, label=right:{\footnotesize 2}] at (0,0) {}; 
				\node[vertex, label=right:{\footnotesize 3}] at (1,1) {}; 
				\node[vertex, label=right:{\footnotesize 4}] at (1,0) {}; 
				\node at (2,0) {,}; 
				\node[vertex, label=right:{\footnotesize 1}] at (3,1) {}; 
				\node[vertex, label=right:{\footnotesize 2}] at (3,0) {}; 
				\node[vertex, label=right:{\footnotesize 3}] at (4,1) {}; 
				\node[vertex, label=right:{\footnotesize 4}] at (4,0) {}; 
				\node at (5,0) {,}; 
				\node[vertex, label=right:{\footnotesize 1}] at (6,1) {}; 
				\node[vertex, label=right:{\footnotesize 2}] at (6,0) {}; 
				\node[vertex, label=right:{\footnotesize 3}] at (7,1) {}; 
				\node[vertex, label=right:{\footnotesize 4}] at (7,0) {}; 
				\node at (-1,0.5) {{\Huge$\Bigg\{$}}; 
				\node at (8,0.5) {{\Huge$\Bigg\}$}}; 
				\node at (3.5,2) {{\huge$\Downarrow$}}; 
				\node[vertex, label=right:{\footnotesize 1}] at (2.5,4) {}; 
				\node[vertex, label=right:{\footnotesize 2}] at (2.5,3) {}; 
				\node[vertex, label=right:{\footnotesize 3}] at (4.5,4) {}; 
				\node[vertex, label=right:{\footnotesize 4}] at (4.5,3) {}; 
				\node at (3.5,3) {,}; 
				\begin{pgfonlayer}{background}
					\draw[hyedge1, color=blue] (0,0) -- (0,0);
					\draw[hyedge2] (0,0) -- (0,0);
					\draw[hyedge1, color=green] (0,1) -- (1,1) -- (1,0);
					\draw[hyedge2] (0,1) -- (1,1) -- (1,0);
					\draw[hyedge1, color=red] (3,1) -- (3,1);
					\draw[hyedge2] (3,1) -- (3,1);
					\draw[hyedge1, color=green] (3,0) -- (4,0) -- (4,1);
					\draw[hyedge2] (3,0) -- (4,0) -- (4,1);
					\draw[hyedge1, color=red] (6,1) -- (6,1);
					\draw[hyedge2] (6,1) -- (6,1);
					\draw[hyedge1, color=blue] (6,0) -- (6,0);
					\draw[hyedge2] (6,0) -- (6,0);
					\draw[hyedge1, color=green] (7,0) -- (7,1);
					\draw[hyedge2] (7,0) -- (7,1);
					\draw[hyedge1, color=red] (2.5,4) -- (2.5,4);
					\draw[hyedge2] (2.5,4) -- (2.5,4);
					\draw[hyedge1, color=blue] (2.5,3) -- (2.5,3);
					\draw[hyedge2] (2.5,3) -- (2.5,3);
					\draw[hyedge1, color=green] (4.5,3) -- (4.5,4);
					\draw[hyedge2] (4.5,3) -- (4.5,4);
				\end{pgfonlayer}
			\end{tikzpicture}
			\caption{A visualisation of the set $P \tilde{\cup} Q$}
			\label{fig:coupling}
		\end{figure}
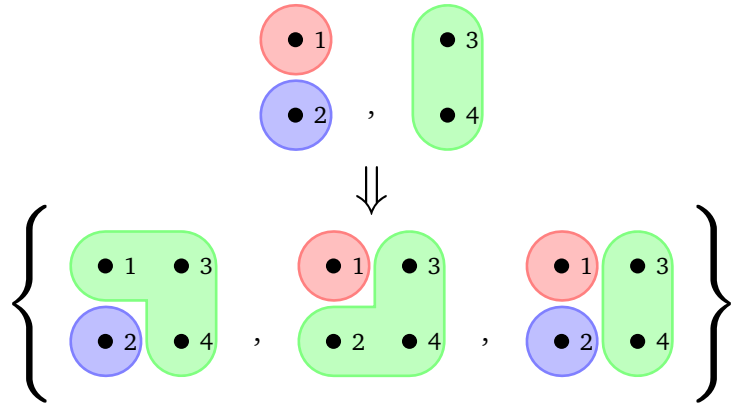
	\end{example}
	
	\begin{definition}
		\label{definition:psi-Mappings}
		Let $\scM, \scN$ be disjoint sets and $P \in \scP(\scM)$, $Q\in \scP(\scN)$. 
		For every $G\in P\tilde{\cup} Q$, we denote the two injective mappings $\psi^{P,G}: P \to G$ and $\psi^{Q,G}: Q \to G$ such that $\forall p\in P$ and $\forall q\in Q$, 
		\begin{equation}
			\label{eq:definition:psi}
			\psi^{P,G}[ p] \cap \scM = p, \quad \psi^{Q,G}[ q] \cap \scN = q. 
		\end{equation}
		Further, we define $\phi^{G, P,Q}: G \to P \cup Q$ and $\varphi^{G, P, Q}:P \to P\cup Q$ as follows:
		\begin{align}
			\label{eq:definition:phi}
			\phi^{G, P,Q} [g] =& 
			\begin{cases} 
				g & \quad \in P \quad \mbox{if } g\cap Q = \emptyset, 
				\\
				g\cap Q \qquad  & \quad \in Q  \quad \mbox{if } g\cap Q \neq \emptyset. 
			\end{cases}
			\\
			\label{eq:definition:varphi}
			\varphi^{G, P, Q}[p] =& 
			\begin{cases}
				p & \quad \in P \quad \mbox{if } \psi^{P, G}[p]\cap Q  = \emptyset, 
				\\
				\psi^{P, G}[p]\cap Q & \quad \in Q \quad \mbox{if } \psi^{P, G}[p]\cap Q  \neq \emptyset. 
			\end{cases}
		\end{align}
	\end{definition}
	The operator $\psi^{P, G}$ (see Equation \eqref{eq:definition:psi}) maps the set $p$, an element of the partition $P$, to the unique element of the partition $G$ that contains the set $p$. 
	
	The operator $\phi^{G, P, Q}$ (see Equation \eqref{eq:definition:phi}) maps the set $g$, an element of the partition $G$ (a coupled partition between $P$ and $Q$) to either the unique element of the partition $Q$ that is contained within $g$ (should that element exist) or the unique of the partition $P$ that is equal to $g$. These two scenarios are mutually exclusive and both elements exist as appropriate. 
	
	On the other hand, the operator $\varphi^{G, P, Q}$ (see Equation \eqref{eq:definition:varphi}) maps the set $p$, an element of the partition $P$ to either the element $q\in Q$ such that $p\cup q \in G$ or to itself if no such element exists in $G$. Again, these two scenarios are mutually exclusive and both elements exist as appropriate. 
	
	Critically, we can now extend the set of couplings into an ordered pairing operation on forests:
	\begin{definition}
		\label{definition:coupled_trees}		
		Let $\Upsilon, Y\in \scF_0$ and let $G \in H^{\Upsilon} \tilde{\cup} H^Y$. We say that $G$ is a \emph{Lions coupling} if
		\begin{equation}
			\label{eq:definition:LionsCoupling}
			\forall h\in H^{\Upsilon} \quad \mbox{such that}\quad h\cap \fr[\Upsilon] = \emptyset \quad \implies \varphi^{G, H^{\Upsilon}, H^Y}[h] = h. 
		\end{equation}
		We denote the set of Lions couplings between two Lions forests $\Upsilon$ and $Y$ by $\lion(\Upsilon, Y)\subseteq H^{\Upsilon} \tilde{\cup} H^Y$. Finally, we define the set
		\begin{equation}
			\label{eq:definition:coupled_trees2}
			\scF_0 \tilde{\times} \scF_0:=\Big\{ \Upsilon \times^G Y: \quad \Upsilon, Y\in \scF_0,\quad G\in \lion(\Upsilon, Y) \Big\} = \bigsqcup_{\Upsilon \times Y \in \scF_0 \times \scF_0} \lion(\Upsilon, Y). 
		\end{equation}
	\end{definition}
	
	In words, for any two Lions forests $\Upsilon, Y\in \scF_0$ we find that the additional hypergraphic structures mean that there are fewer ways to couple the hyperedges in a way that is combinatorially meaningful. It turns out that the set of possible couplings is too large and actually the set of Lions couplings is a better collection of couplings to work with. 
	
	\begin{remark}
		\label{remark:Ttimes-Gtimes}
		For notational purposes, when a coupling between two Lions trees is explicitly known, it will be stated as a superscript of the Cartesian product $\times^G$ as in Definition \ref{definition:coupled_trees}. However, if the coupling is only implicitly known (as will often be the case for operations on the collection of coupled pairs), we will use the notation $\tilde{\times}$ to indicate that there is a coupling but the situation is the same regardless of the choice. 
	\end{remark}
	
	Next, we use the coupled pairs to parametrise a new graded module in the same fashion as Definition \ref{definition:AlgebraRVs}. While the Cartesian product of two index sets $I\times J$ can be used to index the tensor product of the two modules indexed by $I$ and $J$, the set of coupled pairs will be used to index the coupled tensor product. Having defined a coupled pair in Definition \ref{definition:coupled_trees} which corresponds to a richer collection of elements than a Cartesian product, we extend the same ideas to define a coupled tensor product 
	between modules indexed by Lions forests.  
	
	Further, we need to introduce a functional for measuring the integrability of various products of random variables obtained by a tensorisation operation that couples the different probability spaces based on the choice of Lions coupling. 
	\begin{definition}
		\label{definition:XTensorProduct2}
		Let $(\Omega, \cF, \bP)$ and $(\Omega', \cF', \bP')$ be probability spaces. For each $\Upsilon \times^G Y\in \scF$, let
		\begin{align*}
			&\hat{p}[Y] = \big( \hat{p}[Y]_h \big)_{h\in H^Y} \in (1, \infty)^{\times |H^Y|} 
			\quad \mbox{and} \quad
			\hat{p}[\Upsilon, G] = \big( \hat{p}[\Upsilon, G]_h \big)_{h\in H^{\Upsilon}: \varphi^{G, \Upsilon, Y}[h] = h}
			\\
			&\hat{\rp} \in  \bigsqcup_{Y\in \scF} \Bigg( (1, \infty)^{\times |H^Y|} \times \bigg( \bigsqcup_{\substack{\Upsilon \in \scF_0 \\ G\in \lion(\Upsilon, Y)}} (1, \infty)^{\times (|G| - |H^Y|)} \bigg) \Bigg),
			\qquad \hat{\rp} = \Big( \hat{p}[Y], \hat{p}[\Upsilon,G] \Big)_{\Upsilon \times^G Y \in \scF_0 \tilde{\times} \scF_0}
		\end{align*}
		We say that $\hat{\rp} = \big( \hat{p}[Y], \hat{p}[\Upsilon,G] \big)_{\Upsilon \times^G Y \in \scF_0 \tilde{\times} \scF_0}$ is a 2-integrability functional. 
		
		We define the $\bR^e$-module
		\begin{align*}
			&\scH'(\Omega') \tilde{\otimes}_{\rp} \scH'(\Omega') = \bigoplus_{Y\in \scF_0} L^{\hat{p}[Y]}\Bigg( (\Omega')^{\times |H^{Y}|}, (\bP')^{\times |H^{Y}|};
			\\
			&\quad \bigoplus_{\substack{\Upsilon \in \scF_0 \\ G\in \lion(\Upsilon, Y) }} L^{\hat{p}[\Upsilon, G]}\bigg( (\Omega')^{\times (|G| - |H^{Y}|)}, (\bP')^{\times (|G|-|H^{Y}|)}; \lin\Big( (\bR^d)^{\otimes |\scN^{\Upsilon}|}, \lin\big( (\bR^d)^{\otimes |\scN^{Y}|}, \bR^e \big) \Big) \bigg) \Bigg).  
		\end{align*}
		We also define the $L^0(\Omega, \bP; \bR^e)$-module
		\begin{align*}
			&\scH \tilde{\otimes}_{\hat{\rp}} \scH(\Omega, \Omega') = L^0\Big( \Omega, \bP; \scH'(\Omega') \tilde{\otimes}_{\hat{\rp}} \scH'(\Omega') \Big)
			\\
			&= L^0\Bigg( \Omega, \bP; \bigoplus_{Y\in \scF_0} L^{\hat{p}[Y]}\Bigg( (\Omega')^{\times |H^{Y}|}, (\bP')^{\times |H^{Y}|};
			\\
			& \bigoplus_{\substack{\Upsilon \in \scF_0 \\ G\in \lion(\Upsilon, Y) }} L^{\hat{p}[\Upsilon, G]} \bigg( (\Omega')^{\times (|G| - |H^{Y}|)}, (\bP')^{\times (|G|-|H^{Y}|)}; \lin\Big( (\bR^d)^{\otimes |\scN^{\Upsilon}|}, \lin\big( (\bR^d)^{\otimes |\scN^{Y}|}, \bR^e \big) \Big) \bigg) \Bigg) \Bigg).  
		\end{align*}
	\end{definition}

	We use the convention that $\scH\tilde{\otimes}_{\hat{\rp}} \scH(\Omega, \Omega')$ has a representation in terms of the set $\scF_0 \tilde{\times} \scF_0$ (recall Equation \eqref{eq:definition:coupled_trees2}): for $X\in \scH \tilde{\otimes}_{\hat{\rp}} \scH(\Omega, \Omega')$, we write
	\begin{equation*}
		X(\omega_0) = \sum_{Y\in \scF_0} \sum_{\substack{\Upsilon \in \scF_0 \\ G\in \lion(\Upsilon, Y)}} \Big\langle X, 	\Upsilon \times^G Y \Big\rangle(\omega_0, \cdot, \cdot) 
	\end{equation*}
	where for each $Y \in \scF_0$, we have
	\begin{align*}
		&\Big\langle X, \Upsilon \times^G Y \Big\rangle(\omega_0, \omega_{H^{Y}}, \cdot ) 
		\\
		&\in L^{\hat{p}[\Upsilon, G]}\bigg( (\Omega')^{\times(|G| - |H^{Y}|)}, \bP^{{\times(|G| - |H^{Y}|)}} ; \lin\Big( (\bR^d)^{\otimes |\scN^{\Upsilon}|}, \lin\big((\bR^d)^{\otimes |\scN^{Y}|}, \bR^e\big) \Big) \bigg) \Bigg)
	\end{align*}
	so that
	\begin{align*}
		&\sum_{\substack{\Upsilon \in \scF_0 \\ G \in \lion(\Upsilon, Y)}} \Big\langle X, \Upsilon \times^G Y \Big\rangle(\omega_0, \omega_{H^{Y}}, \cdot ) 
		\\
		&\in \bigoplus_{\substack{\Upsilon \in \scF_0 \\ G \in \lion(\Upsilon, Y)}} L^{p[\Upsilon, G]}\bigg( (\Omega')^{\times(|G| - |H^{Y}|)}, \bP^{{\times(|G| - |H^{Y}|)}} ; \lin\Big( (\bR^d)^{\otimes |\scN^{\Upsilon}|}, \lin\big((\bR^d)^{\otimes |\scN^{Y}|}, \bR^e\big) \Big) \bigg) \Bigg)
	\end{align*}
	and
	\begin{align*}
		&\sum_{Y \in \scF_0} \sum_{\substack{\Upsilon \in \scF_0 \\ G\in \lion( \Upsilon,Y) }} \Big\langle X, \Upsilon \times^G Y \Big\rangle (\omega_0, \cdot, \cdot)
		\\
		&\in \bigoplus_{Y\in \scF_0} L^{\hat{p}[Y]} \Bigg( (\Omega')^{\times |H^{Y}|}, (\bP')^{\times |H^Y|}; 
		\\
		&\bigoplus_{\substack{\Upsilon \in \scF_0 \\ G\in \lion(\Upsilon, Y) }} L^{\hat{p}[\Upsilon, G]} \bigg( (\Omega')^{\times (|G| - |H^Y|)}, (\bP')^{\times (|G| - |H^Y|)}; \lin\Big( (\bR^d)^{\otimes |\scN^{\Upsilon}|}, \lin\big((\bR^d)^{\otimes |\scN^{Y}|}, \bR^e\big) \Big) \bigg) \Bigg)
	\end{align*}
	
	For a more detailed discussion of what a coupled tensor product as part of a wider discussion with regard to topological algebras, we refer the reader to \cite{salkeld2022LionsTrees}. 
	
	\begin{remark}
		Following on from Definition \ref{definition:psi-Mappings}, when there is no ambiguity we will often equivalently write
		\begin{align*}
			\phi^{T, \Upsilon, Y} = \phi^{H^T, H^{\Upsilon}, H^{Y}} 
			\quad \mbox{for} \quad T, \Upsilon, Y \in \scF, 
			\\
			\varphi^{T, \Upsilon, Y} = \varphi^{H^T, H^{\Upsilon}, H^{Y}} 
			\quad \mbox{for} \quad T, \Upsilon, Y \in \scF. 
		\end{align*}
		Similarly, we have that
		\begin{align*}
			\phi^{G, \Upsilon, Y} = \phi^{G, H^{\Upsilon}, H^{Y}} 
			&\quad \mbox{and} \quad
			\varphi^{G, \Upsilon, Y} = \varphi^{G, H^{\Upsilon}, H^{Y}}
			\\
			&\mbox{for} \quad \Upsilon, Y \in \scF_0 \quad \mbox{and} \quad G\in \lion(\Upsilon, Y). 
		\end{align*}
	\end{remark}
		
	\subsubsection{Cuts for Lions trees}
	
	An admissible cut is a way of dividing a directed tree into two subtrees, one a directed tree (referred to as the root) and one a directed forest (referred to as the prune), see for instance \cite{connes1999hopf}. In the context of Lions forests, a cut is a subset of the edges of a Lions forest that can be removed for combinatorial purposes but this operation should not alter the hypergraphic structure. Thus, a cut takes a Lions tree to a coupled pair of Lions forests contained in $\scF_0 \tilde{\times} \scF_0$ rather than simply $\scF_0 \times \scF_0$ as one might naively expect in the Butcher-Connes-Kreimer setting. 
	
	\begin{definition}
		\label{definition:CutsCoproduct}
		Let $T=(\scN, \scE, h_0, H, \scL) \in \scT_d$ be a non-empty labelled Lions tree. A subset $c\subseteq \scE$ is called an admissible cut if $\forall y\in \scN$, the unique path $(e_i)_{i=1, ..., n}$ from $y$ to the root $x$ satisfies that if $e_i\in c$ then $\forall j\neq i, e_j\notin c$. The set of admissible cuts for the Lions tree $T$ is denoted $C(T)$. 
		
		We call the tuple $T_c^R$ the root of the cut $c$ of the labelled Lions tree $T$ 
		\begin{align*}
			T_c^R &= (\scN_c^R, \scE_c^R, h_0 \cap \scN_c^R, H_c^R, \scL_c^R), \quad \mbox{where}
			\\
			\scN_c^R:&= \big\{ y\in \scN: \exists (y_i)_{i=1, ..., n}\in \scN, y_1 = y, y_n=x, (y_i, y_{i+1})\in \scE \backslash c \big\}, 
			\\
			\scE_c^R:&= \big\{ (y,z)\in \scE: y,z\in \scN_c^R \big\}, 
			\\
			H_c^R:&= \big\{ h\cap \scN_c^R: h\in H \big\}\backslash \{\emptyset\}, 
			\\
			\scL_c^R:& \scN_c^R \to \{1, ..., d\}, \quad \scL_c^R = \scL|_{\scN_c^R}. 
		\end{align*} 
		
		The prune of the cut $c$ is the Lions forest $(\scN_c^P, \scE_c^P, h_0\cap \scN_c^P, H_c^p, \scL_c^P)$ where $\scN_c^P = \scN\backslash \scN_c^R$ and $\scE_c^P = \scE\backslash (\scE_c^R \cup c)$, $H_c^P:= \{ h\cap \scN_c^P: h\in H \}\backslash \emptyset$ and $\scL_c^P = \scL|_{\scN_c^R}$.  
	\end{definition}
	
	Cuts are defined as in the usual Butcher-Connes-Kreimer setting but where the hyperedges of the roots (respectively prunes) are constructed accordingly by restricting the original hyperedges to the nodes of the roots (respectively prunes). Further, the root and prune do not form a pair, but a coupled pair with the coupling determined by the original collection of hyperedges. 
	
	\begin{example}
		Consider the following admissible cut applied to the tree:
		$$
		
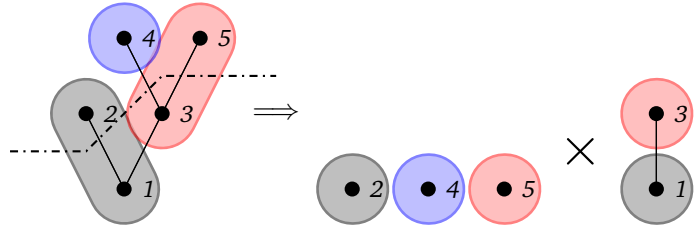
\begin{tikzpicture}
			\node[vertex, label=right:{\footnotesize 1}] at (0,0) {}; 
			\node[vertex, label=right:{\footnotesize 2}] at (-0.5,1) {}; 
			\node[vertex, label=right:{\footnotesize 3}] at (0.5,1) {}; 
			\node[vertex, label=right:{\footnotesize 4}] at (0,2) {}; 
			\node[vertex, label=right:{\footnotesize 5}] at (1,2) {}; 
			\draw[edge] (0,0) -- (-0.5, 1);
			\draw[edge] (0,0) -- (0.5, 1);
			\draw[edge] (0.5, 1) -- (0,2);
			\draw[edge] (0.5, 1) -- (1,2);
			\draw[thick,dash dot] (-1.5,0.5) -- (-0.5,0.5) -- (0.5,1.5) -- (2,1.5);
			\node at (2,1) {$\implies$}; 
			\node[vertex, label=right:{\footnotesize 2}] at (3,0) {}; 
			\node[vertex, label=right:{\footnotesize 4}] at (4,0) {}; 
			\node[vertex, label=right:{\footnotesize 5}] at (5,0) {}; 
			\node at (6,0.5) {{\huge$\times$}}; 
			\node[vertex, label=right:{\footnotesize 1}] at (7,0) {}; 
			\node[vertex, label=right:{\footnotesize 3}] at (7,1) {}; 
			\draw[edge] (7, 0) -- (7,1);
			\begin{pgfonlayer}{background}
				\draw[zhyedge1] (0,0) -- (-0.5,1);
				\draw[zhyedge2] (0,0) -- (-0.5,1);
				\draw[hyedge1, color=red] (0.5,1) -- (1,2);
				\draw[hyedge2] (0.5,1) -- (1,2);
				\draw[hyedge1, color=blue] (0,2) -- (0,2);
				\draw[hyedge2] (0,2) -- (0,2);
				\draw[zhyedge1] (3,0) -- (3,0);
				\draw[zhyedge2] (3,0) -- (3,0);
				\draw[hyedge1, color=blue] (4,0) -- (4,0);
				\draw[hyedge2] (4,0) -- (4,0);
				\draw[hyedge1, color=red] (5,0) -- (5,0);
				\draw[hyedge2] (5,0) -- (5,0);
				\draw[zhyedge1] (7,0) -- (7,0);
				\draw[zhyedge2] (7,0) -- (7,0);
				\draw[hyedge1, color=red] (7,1) -- (7,1);
				\draw[hyedge2] (7,1) -- (7,1);
			\end{pgfonlayer}
		\end{tikzpicture}
		$$
		The node $5$ and $3$ are both in a hyperedge that is coloured red. The colouring here represents that the hyperedges are coupled. 
	\end{example}
	
	A cut can be identified with a mapping that takes a Lions tree to a collection of coupled pairs of Lions forests. Critically, this coupling is necessarily a Lions coupling (see Definition \ref{definition:coupled_trees}):
	\begin{lemma}
		\label{lemma:Lions-coupling}
		Let $T = (\scN, \scE, h_0, H, \scL) \in \scT_d$ and let $c\in C(T)$. Then $T_c^R, T_c^P \in \scF_d$ and $H^T \in \lion\big( H_c^{P}, H_c^{R} \big)$. 
	\end{lemma}
	The proof of Lemma \ref{lemma:Lions-coupling} can be found in \cite{salkeld2022LionsTrees}. In the classical Butcher-Connes-Kreimer setting (see Equation \eqref{eq:francois:label:CK:co-product}), cuts of directed trees are intimately connected with the combinatorial properties of the Butcher-Connes-Kreimer coproduct. For Lions trees, the relationship is a little more complicated since the image of a coupled coproduct is the coupled tensor product from Definition \ref{definition:XTensorProduct2}. 
	\begin{definition}
		\label{definition:coproduct}
		Let 
		\begin{equation*}
			\Delta: \mbox{Span}\big( \scF_{0,d}[0] \big) \to \mbox{Span}\big( \scF_{0,d}[0] \tilde{\times} \scF_{0,d}[0] \big)
		\end{equation*}
		be the linear operator such that for $T=(\scN,\scE,h_0,H,\scL) \in \scT_{d}$, 
		\begin{equation}
			\label{eq:proposition:Coproduct2-Equivalence}
			\Delta\Big[ T \Big] = \rId \times^H T + T \times^H \rId + \sum_{c\in C(T)} T_c^P \times^{H^T} T_c^R. 
		\end{equation}
		For Lions forest $T=(\scN, \scE, h_0, H, \scL) \in \scF$ such that $|\fr[T]| = n > 1$, we have that the labelled directed forest $(\scN, \scE, \scL) \in \fF$ can be expressed as a product of labelled directed trees
		\begin{align*}
			(\scN, \scE, \scL) =& \tau_1 \odot ... \odot \tau_n
			= 
			\big( \scN^{\tau_1}, \scE^{\tau_1}, \scL^{\tau_1} \big) \odot ... \odot \big( \scN^{\tau_n}, \scE^{\tau_n}, \scL^{\tau_n} \big). 
		\end{align*}
		For each $i = 1, ..., |\fr[T]|$, we define the set
		\begin{align*}
			\fD_i(T) = \bigg\{& \Big( (\scN^{\tau_i}, \scE^{\tau_i}, \scL^{\tau_i}), (\emptyset, \emptyset, \emptyset) \Big), \quad
			\Big( (\emptyset, \emptyset, \emptyset), (\scN^{\tau_i}, \scE^{\tau_i}, \scL^{\tau_i}) \Big), 
			\\
			&\Big( \big((\scN^{\tau_i})_c^P, (\scE^{\tau_i})_c^P, (\scL^{\tau_i})_c^P \big), \big((\scN^{\tau_i})_c^R, (\scE^{\tau_i})_c^R, (\scL^{\tau_i})_c^R \big) \Big): c\in C(\tau_i) \bigg\} \subset \fF_0 \times \fF_0. 
		\end{align*}
		Next, we define the set
		\begin{align*}
			\fD(T) = \prod_{i=1}^n \fD_i(T) \subset \prod_{i=1}^n \fF_0 \times \fF_0
		\end{align*}
		and for each $\overline{\tau} = \big( (\tau_1^p, \tau_1^r), ..., (\tau_n^p, \tau_n^r)\big) \in \fD(T)$, we denote
		\begin{align*}
			&\scN_{\overline{\tau}}^p = \bigcup_{i=1}^n \scN^{\tau_i^p}, 
			\quad 
			\scN_{\overline{\tau}}^r = \bigcup_{i=1}^n \scN^{\tau_i^r},
			\qquad
			\scE_{\overline{\tau}}^p = \bigcup_{i=1}^n \scE^{\tau_i^p}, 
			\quad 
			\scE_{\overline{\tau}}^r = \bigcup_{i=1}^n \scE^{\tau_i^r},
			\\
			&H_{\overline{\tau}}^p := \big\{ h \cap \scN_{\overline{\tau}}^p : h\in H \big\} \backslash \{\emptyset\}, 
			\quad
			H_{\overline{\tau}}^r := \big\{ h \cap \scN_{\overline{\tau}}^r : h\in H \big\} \backslash \{\emptyset\}, 
			\\
			&\scL_{\overline{\tau}}^p = \scL\big|_{\scN_{\overline{\tau}}^p},
			\qquad
			\scL_{\overline{\tau}}^r = \scL\big|_{\scN_{\overline{\tau}}^r}. 
		\end{align*}
		
		Then we extend $\Delta$ to the set of Lions forests
		\begin{equation}
			\label{eq:proposition:Coproduct2-Equivalence2}
			\Delta\Big[T\Big] = \sum_{\overline{\tau} \in \fD(T)} \Big( \scN_{\overline{\tau}}^p, \scE_{\overline{\tau}}^p, h_0 \cap \scN_{\overline{\tau}}^p, H_{\overline{\tau}}^p, \scL_{\overline{\tau}}^p \Big) 
			\times^{H^T} 
			\Big( \scN_{\overline{\tau}}^r, \scE_{\overline{\tau}}^r, h_0 \cap \scN_{\overline{\tau}}^r, H_{\overline{\tau}}^r, \scL_{\overline{\tau}}^r \Big). 
		\end{equation}		
		We pair the operator $\Delta$ with the linear functional $\epsilon: \mbox{Span}\big( \scF_{0,d}[0] \big) \to \bR^e$ which satisfies $\epsilon(\rId) = 1$, and $\forall T\in \scF$, $\epsilon(T) = 0$. The operator $\epsilon$ is the counit of $\Delta$. 
	\end{definition}

	For efficient notation, we introduce the coproduct counting function $c:\scF_0 \times \scF_0 \times \scF_0 \to \bN_0$ defined for $T, \Upsilon, Y \in \scF_0$ by
	\begin{equation}
		\label{eq:Coproduct-counting}
		c \Big( T, \Upsilon, Y \Big) = \Big\langle \Delta[T], \Upsilon \times^{H^{T}} Y \Big\rangle. 
	\end{equation}
	This yields the expression 
	\begin{equation*}
		\Delta[ T ] = \sum_{\Upsilon, Y \in \scF_0} c \Big( T, \Upsilon, Y \Big) \Upsilon \times^{H^{T}} Y. 
	\end{equation*}
	
	\begin{example}
		\label{example:Coproduct}
		The elements of $\scF_0 \tilde{\times} \scF_0$ should be considered and treated as being very different from elements of the Cartesian product $\scF_0 \times \scF_0$. Thus, we visualise them differently to highlight this. 
		
		The coupled coproduct of the tree
		
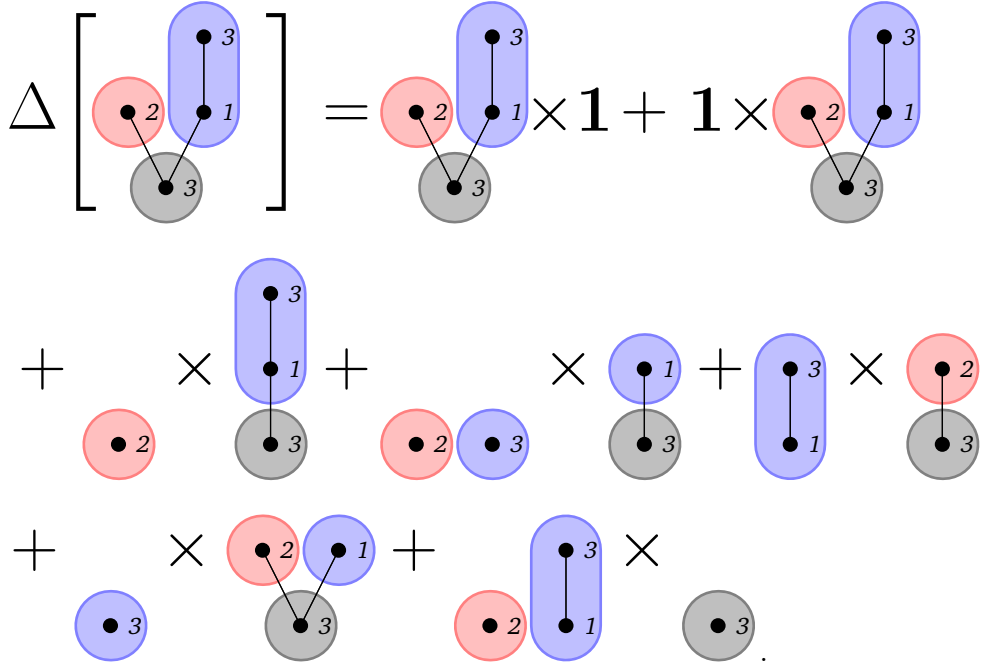
\begin{align*}
			\begin{tikzpicture}
				\node[vertex, label=right:{\footnotesize 3}] at (1,2) {}; 
				\node[vertex, label=right:{\footnotesize 2}] at (0,1) {}; 
				\node[vertex, label=right:{\footnotesize 3}] at (0.5,0) {}; 
				\node[vertex, label=right:{\footnotesize 1}] at (1,1) {}; 
				\node at (-1,1) {{\Huge$\Delta \Bigg[$}}; 
				\node at (2.5,1) {{\Huge$\Bigg]=$}}; 
				\draw[edge] (0.5,0) -- (0,1);
				\draw[edge] (0.5,0) -- (1,1);
				\draw[edge] (1,1) -- (1,2);
				\begin{pgfonlayer}{background}
					\draw[zhyedge1] (0.5,0) -- (0.5,0);
					\draw[zhyedge2] (0.5,0) -- (0.5,0);
					\draw[hyedge1, color=red] (0,1) -- (0,1);
					\draw[hyedge2] (0,1) -- (0,1);
					\draw[hyedge1, color=blue] (1,1) -- (1,2);
					\draw[hyedge2] (1,1) -- (1,2);
				\end{pgfonlayer}
			\end{tikzpicture}
			&
			\begin{tikzpicture}
				\node[vertex, label=right:{\footnotesize 3}] at (1,2) {}; 
				\node[vertex, label=right:{\footnotesize 2}] at (0,1) {}; 
				\node[vertex, label=right:{\footnotesize 3}] at (0.5,0) {}; 
				\node[vertex, label=right:{\footnotesize 1}] at (1,1) {}; 
				\node at (2,1) {{\Huge$\times \rId$}}; 
				\node at (3,1) {{\Huge$+$}}; 
				\draw[edge] (0.5,0) -- (0,1);
				\draw[edge] (0.5,0) -- (1,1);
				\draw[edge] (1,1) -- (1,2);
				\begin{pgfonlayer}{background}
					\draw[zhyedge1] (0.5,0) -- (0.5,0);
					\draw[zhyedge2] (0.5,0) -- (0.5,0);
					\draw[hyedge1, color=red] (0,1) -- (0,1);
					\draw[hyedge2] (0,1) -- (0,1);
					\draw[hyedge1, color=blue] (1,1) -- (1,2);
					\draw[hyedge2] (1,1) -- (1,2);
				\end{pgfonlayer}
			\end{tikzpicture}
			\begin{tikzpicture}
				\node[vertex, label=right:{\footnotesize 3}] at (1,2) {}; 
				\node[vertex, label=right:{\footnotesize 2}] at (0,1) {}; 
				\node[vertex, label=right:{\footnotesize 3}] at (0.5,0) {}; 
				\node[vertex, label=right:{\footnotesize 1}] at (1,1) {}; 
				\node at (-1,1) {{\Huge$\rId\times $}};  
				\draw[edge] (0.5,0) -- (0,1);
				\draw[edge] (0.5,0) -- (1,1);
				\draw[edge] (1,1) -- (1,2);
				\begin{pgfonlayer}{background}
					\draw[zhyedge1] (0.5,0) -- (0.5,0);
					\draw[zhyedge2] (0.5,0) -- (0.5,0);
					\draw[hyedge1, color=red] (0,1) -- (0,1);
					\draw[hyedge2] (0,1) -- (0,1);
					\draw[hyedge1, color=blue] (1,1) -- (1,2);
					\draw[hyedge2] (1,1) -- (1,2);
				\end{pgfonlayer}
			\end{tikzpicture}
			\\
			\begin{tikzpicture}
				\node[vertex, label=right:{\footnotesize 2}] at (0,0) {}; 
				\node at (1,1) {{\Huge$\times$}};
				\node at (3,1) {{\Huge$+$}};
				\node at (-1,1) {{\Huge$+$}};
				\node[vertex, label=right:{\footnotesize 3}] at (2,0) {}; 
				\node[vertex, label=right:{\footnotesize 1}] at (2,1) {}; 
				\node[vertex, label=right:{\footnotesize 3}] at (2,2) {}; 
				\draw[edge] (2,0) -- (2,1);
				\draw[edge] (2,1) -- (2,2);
				\begin{pgfonlayer}{background}
					\draw[hyedge1, color=red] (0,0) -- (0,0);
					\draw[hyedge2] (0,0) -- (0,0);
					\draw[zhyedge1] (2,0) -- (2,0);
					\draw[zhyedge2] (2,0) -- (2,0);
					\draw[hyedge1, color=blue] (2,1) -- (2,2);
					\draw[hyedge2] (2,1) -- (2,2);
				\end{pgfonlayer}
			\end{tikzpicture}
			&
			\begin{tikzpicture}
				\node[vertex, label=right:{\footnotesize 2}] at (0,0) {}; 
				\node[vertex, label=right:{\footnotesize 3}] at (1,0) {}; 
				\node at (2,1) {{\Huge$\times$}};
				\node at (4,1) {{\Huge$+$}};
				\node[vertex, label=right:{\footnotesize 3}] at (3,0) {}; 
				\node[vertex, label=right:{\footnotesize 1}] at (3,1) {}; 
				\draw[edge] (3,0) -- (3,1);
				\begin{pgfonlayer}{background}
					\draw[hyedge1, color=red] (0,0) -- (0,0);
					\draw[hyedge2] (0,0) -- (0,0);
					\draw[hyedge1, color=blue] (1,0) -- (1,0);
					\draw[hyedge2] (1,0) -- (1,0);
					\draw[zhyedge1] (3,0) -- (3,0);
					\draw[zhyedge2] (3,0) -- (3,0);
					\draw[hyedge1, color=blue] (3,1) -- (3,1);
					\draw[hyedge2] (3,1) -- (3,1);
				\end{pgfonlayer}
			\end{tikzpicture}
			\begin{tikzpicture}
				\node[vertex, label=right:{\footnotesize 1}] at (0,0) {}; 
				\node[vertex, label=right:{\footnotesize 3}] at (0,1) {}; 
				\node at (1,1) {{\Huge$\times$}};
				\node[vertex, label=right:{\footnotesize 3}] at (2,0) {}; 
				\node[vertex, label=right:{\footnotesize 2}] at (2,1) {}; 
				\draw[edge] (0,0) -- (0,1);
				\draw[edge] (2,0) -- (2,1);
				\begin{pgfonlayer}{background}
					\draw[hyedge1, color=blue] (0,0) -- (0,1);
					\draw[hyedge2] (0,0) -- (0,1);
					\draw[hyedge1, color=red] (2,1) -- (2,1);
					\draw[hyedge2] (2,1) -- (2,1);
					\draw[zhyedge1] (2,0) -- (2,0);
					\draw[zhyedge2] (2,0) -- (2,0);
				\end{pgfonlayer}
			\end{tikzpicture}
			\\
			\begin{tikzpicture}
				\node[vertex, label=right:{\footnotesize 3}] at (0,0) {}; 
				\node at (1,1) {{\Huge$\times$}};
				\node at (-1,1) {{\Huge$+$}};
				\node[vertex, label=right:{\footnotesize 3}] at (2.5,0) {}; 
				\node[vertex, label=right:{\footnotesize 2}] at (2,1) {}; 
				\node[vertex, label=right:{\footnotesize 1}] at (3,1) {}; 
				\draw[edge] (2.5,0) -- (2,1);
				\draw[edge] (2.5,0) -- (3,1);
				\begin{pgfonlayer}{background}
					\draw[zhyedge1] (2.5,0) -- (2.5,0);
					\draw[zhyedge2] (2.5,0) -- (2.5,0);
					\draw[hyedge1, color=blue] (0,0) -- (0,0);
					\draw[hyedge2] (0,0) -- (0,0);
					\draw[hyedge1, color=blue] (3,1) -- (3,1);
					\draw[hyedge2] (3,1) -- (3,1);
					\draw[hyedge1, color=red] (2,1) -- (2,1);
					\draw[hyedge2] (2,1) -- (2,1);
				\end{pgfonlayer}
			\end{tikzpicture}
			&
			\begin{tikzpicture}
				\node[vertex, label=right:{\footnotesize 2}] at (0,0) {}; 
				\node[vertex, label=right:{\footnotesize 1}] at (1,0) {}; 
				\node[vertex, label=right:{\footnotesize 3}] at (1,1) {}; 
				\node at (2,1) {{\Huge$\times$}};
				\node at (-1,1) {{\Huge$+$}};
				\node[vertex, label=right:{\footnotesize 3}] at (3,0) {}; 
				\draw[edge] (1,0) -- (1,1);
				\begin{pgfonlayer}{background}
					\draw[zhyedge1] (3,0) -- (3,0);
					\draw[zhyedge2] (3,0) -- (3,0);
					\draw[hyedge1, color=red] (0,0) -- (0,0);
					\draw[hyedge2] (0,0) -- (0,0);
					\draw[hyedge1, color=blue] (1,0) -- (1,1);
					\draw[hyedge2] (1,0) -- (1,1);
				\end{pgfonlayer}
			\end{tikzpicture}. 
		\end{align*}
		All the products are indeed elements of $\scF_0 \tilde{\times} \scF_0$, even though the product is denoted by $\times$. Implicitly, the hyperedges on either side of $\times$ are coupled if they have the same colour.
		
		What we see with this formula is that hyperedges on either side of the product can be coupled or disjoint and that these two scenarios are distinct. In particular, this means that
		\begin{equation*}
			\begin{tikzpicture}
				\node[vertex, label=right:{\footnotesize i}] at (0,0) {}; 
				\node at (1,0) {{\Huge$\times$}};
				\node[vertex, label=right:{\footnotesize j}] at (2,0) {}; 
				\node at (3,0) {{\Huge$\neq$}};
				\node[vertex, label=right:{\footnotesize i}] at (4,0) {};  
				\node at (5,0) {{\Huge$\times$}};
				\node[vertex, label=right:{\footnotesize j}] at (6,0) {}; 
				\begin{pgfonlayer}{background}
					\draw[hyedge1, color=red] (0,0) -- (0,0);
					\draw[hyedge2] (0,0) -- (0,0);
					\draw[hyedge1, color=blue] (2,0) -- (2,0);
					\draw[hyedge2] (2,0) -- (2,0);
					\draw[hyedge1, color=red] (4,0) -- (4,0);
					\draw[hyedge2] (4,0) -- (4,0);
					\draw[hyedge1, color=red] (6,0) -- (6,0);
					\draw[hyedge2] (6,0) -- (6,0);
				\end{pgfonlayer}
			\end{tikzpicture}. 
		\end{equation*}
		If this were a Cartesian product, these two elements would be equal since the trees on the left and right hand side of the product symbol are equal (colour only indicates correlation for the joint distribution of the two random variables we are producting). 
		
		The reader could object that if colour were a fixed feature of a hyperedge, then the standard Cartesian product would suffice for pairing two forests. However, this would be different from what we 
		do here: colours on the figure are just used to distinguish hyperedges, but they are ``free'' in the sense that the interpretation of all the terms above would be exactly the same if we replaced red by blue and conversely. Obviously, the latter is consistent with our choice not to label the hyperedges.
	\end{example}
	
	Now we use the natural extension to define a coproduct on our module indexed by Lions forests:
	\begin{definition}
		Let $(\Omega', \cF', \bP')$ be a probability space, let $\rp$ be an integrability functional and let $\hat{\rp}$ be a 2-integrability functional such that for all $T, \Upsilon, Y \in \scF_0$ such that $c(T, \Upsilon, Y)>0$, 
		\begin{equation}
			\label{eq:definition:coproduct2}
			\begin{aligned}
				\forall h \in H^Y 
				\quad & \quad
				\hat{p}[Y]_h \leq p[T]_{(\phi^{T, \Upsilon, Y})^{-1}[h]}, 
				\\
				\forall h \in H^\Upsilon: \varphi^{T, \Upsilon, Y}[h] = h 
				\quad & \quad
				\hat{p}[\Upsilon, T]_h \leq p[T]_{(\phi^{T, \Upsilon, Y})^{-1}[h]}. 
			\end{aligned}
		\end{equation}
		
		We define
		\begin{equation*}
			\Delta: \scH_{\rp}'(\Omega') \to \scH'(\Omega') \tilde{\otimes}_{\hat{\rp}} \scH'(\Omega')
		\end{equation*}
		for $X\in \scH_{\rp}'(\Omega')$ to be
		\begin{align*}
			\Delta\Big[ X \Big] =& \sum_{T\in \scF_0} \Big\langle X,  \Delta\big[T\big] \Big\rangle(\cdot_{\phi^{T, \Upsilon, Y}}) 
			= \sum_{Y\in \scF_0} \sum_{T, \Upsilon \in \scF_0} c\Big( T, \Upsilon, Y \Big) \Big\langle X, T\Big\rangle(\cdot_{\phi^{T, \Upsilon, Y}}). 
		\end{align*}
		
		Let $(\Omega, \cF, \bP)$ be a second probability space. Then we can extend $\Delta: \scH_{\rp}(\Omega, \Omega') \to \scH \tilde{\otimes}_{\hat{\rp}} \scH(\Omega, \Omega')$ canonically for $X\in \scH_{\rp}(\Omega, \Omega')$ to be
		\begin{align*}
			\Delta\Big[ X(\omega_0) \Big] =& \sum_{T\in \scF_0} \Big\langle X, \Delta\big[T\big] \Big\rangle(\omega_0, \cdot_{\phi^{T, \Upsilon, Y}}) 
			= 
			\sum_{Y\in \scF_0} \sum_{T, \Upsilon \in \scF_0} c\Big( T, \Upsilon, Y \Big) \cdot \Big\langle X, T \Big\rangle(\omega_0, \cdot_{\phi^{T, \Upsilon, Y}}) . 
		\end{align*}
		We pair this $\Delta$ with the linear functional $\epsilon: \scH_\rp(\Omega, \Omega') \to L^0(\Omega, \bP; \bR^e)$ which satisfies $\epsilon(\rId) = 1$, and $\forall T\in \scF$, $\epsilon(T) = 0$ (recall that $\scF$ does not contain the empty forest) so that $\epsilon$ is the counit of $\Delta$. 
	\end{definition}
	
	Let $T, \Upsilon, Y\in \scF$ and suppose that $c(T, \Upsilon, Y)>0$. Recalling Definition \ref{definition:psi-Mappings}, the collection of free variables is written as 
	\begin{align*}
		\omega_{\phi^{T, \Upsilon, Y}[H^T]} &= ( \underbrace{\omega_{\phi^{T, \Upsilon, Y}[h]}, ...}_{h\in H^T} ) = \big( \underbrace{\omega_{\phi^{T, \Upsilon, Y}[h]}, ...}_{h\in H^T: h \cap \scN^{Y} \neq \emptyset}, \underbrace{\omega_{\phi^{T, \Upsilon, Y}[h]}, ...}_{h\in H^T: h \cap \scN^{Y} = \emptyset} \big)
		\\
		&= \big( \underbrace{\omega_{h \cap \scN^{Y}}, ...}_{h\in H^T: h \cap \scN^{Y} \neq \emptyset}, \underbrace{\omega_{h}, ...}_{h\in H^T: h \cap \scN^{Y} = \emptyset} \big)
	\end{align*}
	so that the random variable 
	\begin{equation*}
		\Big\langle X, T \Big\rangle(\omega_0, \omega_{H^T}) \in L^0 \bigg( \Omega, \bP; L^{p[T]}\Big( (\Omega')^{\times|H^T|}, (\bP')^{\times|H^T|}; \lin\big( (\bR^d)^{\otimes |\scN^T|}, \bR^e \big) \Big) \bigg)
	\end{equation*}
	is embedded to the random variable
	\begin{align*}
		&\Big\langle X, T\Big\rangle(\omega_0, \omega_{\phi^{T, \Upsilon, Y}[H^T] }) 
		\\
		&\in L^0\Bigg( \Omega, \bP; L^{\hat{p}[Y]}\Bigg( (\Omega')^{\times|H^{Y}|}, (\bP')^{\times|H^{Y}|}; 
		\\
		&\qquad L^{\hat{p}[\Upsilon, G]}\bigg( (\Omega')^{\times (|H^T| -|H^{Y}| )}, (\bP')^{\times (|H^T| -|H^{Y}| )}; \lin\Big( (\bR^d)^{\otimes |\scN^{\Upsilon}|}, \lin\big( (\bR^d)^{\otimes |\scN^{Y}|}, \bR^e \big) \Big) \bigg) \Bigg) \Bigg). 
	\end{align*}

	\subsubsection{Iterative couplings and coassociativity}
	\label{subsubsec:Iterative-couplings}
	
	Following on from Definition \ref{definition:XTensorProduct2} we have the $\bR^e$-module
	\begin{align*}
		&\Big( \scH'(\Omega') \tilde{\otimes}_{\rp} \scH'(\Omega') \Big) \tilde{\otimes}_{\hat{\rp}} \scH'(\Omega') = \scH'(\Omega') \tilde{\otimes}_{\hat{\rp}} \Big( \scH'(\Omega') \tilde{\otimes}_{\rp} \scH'(\Omega') \Big)
		\\
		&= \bigoplus_{Y\in \scF_0}L^{p[Y]} \Bigg( (\Omega')^{\times |H^Y|}, (\bP')^{\times |H^Y|}; \bigoplus_{\substack{\Upsilon \in \scF_0 \\ G\in \lion(\Upsilon, Y)}} L^{p[\Upsilon, G]}\bigg( (\Omega')^{\times (|G| - |H^Y|)}, (\bP')^{\times (|G| - |H^Y|)}; 
		\\
		&\quad \bigoplus_{\substack{T \in \scF_0 \\ G' \in \lion(H^T, G)}} L^{p[T, G']} \Big( (\Omega')^{\times(|G'| - |G|)}, (\bP')^{\times(|G'| - |G|)}; 
		\\
		&\qquad \lin\big( (\bR^d)^{\otimes |\scN^T|}, \lin\big( (\bR^d)^{\otimes |\scN^{\Upsilon}|}, \lin\big( (\bR^d)^{\otimes |\scN^Y|}, \bR^e\big) \big) \big) \Big) \bigg) \Bigg)
	\end{align*}
	and the  $L^0(\Omega, \bP; \bR^e)$-modules
	\begin{align*}
		&\Big(\scH \tilde{\otimes}_{\rp} \scH \Big) \tilde{\otimes}_{\hat{\rp}} \scH (\Omega, \Omega') = \scH \tilde{\otimes}_{\hat{\rp}} \Big(\scH \tilde{\otimes}_{\rp} \scH \Big) (\Omega, \Omega')
		\\
		&=L^0\bigg( \Omega, \bP; \Big( \scH'(\Omega') \tilde{\otimes}_\rp \scH'(\Omega') \Big) \tilde{\otimes}_{\hat{\rp}} \scH'(\Omega') \bigg) = L^0\bigg( \Omega, \bP; \scH'(\Omega') \tilde{\otimes}_{\hat{\rp}} \Big( \scH'(\Omega') \tilde{\otimes}_\rp \scH'(\Omega') \Big) \bigg). 
	\end{align*}

	While this is discussed further in much greater detail in \cite{salkeld2022LionsTrees}, for the purposes of this paper we assert that the operators
	\begin{align*}
		\Delta \tilde{\otimes} I:& \scH'(\Omega') \tilde{\otimes}_\rp \scH'(\Omega') \to \Big( \scH'(\Omega) \tilde{\otimes}_{\rp} \scH'(\Omega') \Big) \tilde{\otimes}_{\hat{\rp}} \scH'(\Omega') 
		\\
		I \tilde{\otimes} \Delta:& \scH'(\Omega') \tilde{\otimes}_\rp \scH'(\Omega') \to \scH'(\Omega) \tilde{\otimes}_{\hat{\rp}} \Big( \scH'(\Omega') \tilde{\otimes}_{\rp} \scH'(\Omega') \Big)
	\end{align*}
	can be expressed as the linear extension of
	\begin{align*}
		\Delta \tilde{\otimes} I \Big[ T_1 \times^G T_2 \Big]
		=& 
		\sum_{\Upsilon, Y \in \scF_0} c\Big( T_1, \Upsilon, Y \Big) \cdot \big( \Upsilon_1, Y_1, T_2, G \big)
		\\
		I \tilde{\otimes} \Delta \Big[ T_1 \times^G T_2 \Big]
		=& 
		\sum_{\Upsilon, Y \in \scF_0} c\Big( T_2, \Upsilon, Y \Big) \cdot \big( T_1, \Upsilon, Y, G \big)
	\end{align*}
	where 
	\begin{equation*}
		\big( T_1, \Upsilon, Y, G \big), \big( \Upsilon, Y, T_2, G \big) \in \scF_0 \tilde{\times} \scF_0 \tilde{\times} \scF_0
	\end{equation*} 
	and with a similar extension for 
	\begin{align*}
		\Delta \tilde{\otimes} I:& \scH \tilde{\otimes}_\rp \scH(\Omega, \Omega') \to \Big( \scH \tilde{\otimes}_\rp \scH \Big) \tilde{\otimes}_{\hat{\rp}} \scH(\Omega, \Omega'), 
		\\
		I \tilde{\otimes} \Delta:& \scH \tilde{\otimes}_\rp \scH(\Omega, \Omega') \to \scH \tilde{\otimes}_{\hat{\rp}} \Big( \scH \tilde{\otimes}_\rp \scH \Big) (\Omega, \Omega'). 
	\end{align*} 

	\begin{proposition}
		\label{proposition:coassociativity}
		Let $(\Omega, \cF, \bP)$ and $(\Omega', \cF', \bP')$ be probability spaces. Let $\rp$ be an integrability functional.  
		Let $\scH_{\rp}'(\Omega')$ and $\scH_{\rp}(\Omega, \Omega')$ be the $\bR^e$-module and $L^0(\Omega, \bP; \bR^e)$-module respectively defined in Definition \ref{definition:AlgebraRVs}. 
		
		Then $\big( \scH_{\rp}'(\Omega'), \Delta, \epsilon \big)$ and $\big( \scH_{\rp}(\Omega, \Omega'), \Delta, \epsilon \big)$ are coassociative coupled coalgebras over the rings $(\bR^e, +, \centerdot)$ and $\big( L^0(\Omega, \bP; \bR^e), +, \centerdot \big)$ respectively. 
		
		In particular, $\Delta$ satisfies the commutative identity
		\begin{equation*}
			I \tilde{\otimes} \Delta \circ \Delta = \Delta \tilde{\otimes} I \circ \Delta
			\quad \mbox{and} \quad 
			\centerdot \circ \epsilon \tilde{\otimes} I \circ \Delta = \centerdot \circ I \tilde{\otimes} \epsilon \circ \Delta = I
		\end{equation*}
	\end{proposition}

	The proof of Proposition \ref{proposition:coassociativity} is found in \cite{salkeld2022LionsTrees}. 
	
	\subsection{Graded coupled bialgebras}
	\label{subse:2.5}
	
	Recalling again the classical Butcher-Connes-Kreimer setting, we have that $(\cH_d, \odot, \rId)$ and an associative algebra and $(\cH_d, \triangle, \epsilon)$ a co-associative coalgebra. Further, it is also well known that 
	\begin{equation*}
		\label{eq:Bialgebra-identity}
		\begin{split}
			\triangle \circ \odot =& \Big( \odot \otimes \odot \Big) \circ \Big(\mbox{Twist} \Big) \circ \Big( \triangle \otimes \triangle \Big)
			\\
			\epsilon \otimes \epsilon =& \epsilon \circ \odot \quad \rId \otimes \rId = \triangle \circ \rId 
			\quad \mbox{and}\quad 
			\epsilon \circ \rId = I. 
		\end{split}
	\end{equation*}
	where $\mbox{Twist}:(\cH_d)^{\otimes 4} \to (\cH_d)^{\otimes 4}$ is the linear map that satisfies that $\mbox{Twist}\big[ \tau_1 \otimes \tau_2 \otimes \tau_3 \otimes \tau_4 \big] = \tau_1 \otimes \tau_3 \otimes \tau_2 \otimes \tau_4$. 
	
	Thus $(\cH_d, \odot, \rId, \triangle, \epsilon)$ is a bialgebra over the ring $(\bR^e, +, \centerdot)$. Next, $\cH_d$ has a natural graded structure determined by $|\scN^\tau|$, 
	\begin{align*}
		&\cH_d = \bigoplus_{n \in \bN_0} \cH_{d, n}, \qquad \cH_{d, n} = \bigoplus_{\substack{\tau \in \fF_{0, d} \\ |\scN^{\tau}| = n}} \bR^e 
		\\
		&\cH_{d, m} \odot \cH_{d, n} \subseteq \cH_{d, m+n}, 
		\quad
		\triangle\big[ \cH_{d, n}\big] \subseteq \bigoplus_{p+q=n} \cH_{d, p} \otimes \cH_{d, q}. 
	\end{align*}
	
	We say that $\cH_{d}$ is connected because the subspace associated to the unit of the monoid $(\bN_0, +)$ is isomorphic to the ring of the bialgebra, that is $\cH_{d, 0} = \bR^e$. 

	Finally, for $0<\alpha<1$ and $\gamma>\alpha$, we denote $\fF_0^{\gamma, \alpha}:=\{ \tau \in \fF_0: \alpha\cdot |\scN^\tau|\leq \gamma\}$. Then
	\begin{equation*}
		\cH_d^{\gamma, \alpha} = \bigoplus_{\tau \in \fF_{0, d}}^{\gamma, \alpha} \bR^e = \bigoplus_{\tau \in \fF_{0, d}^{\gamma, \alpha}} \bR^e
	\end{equation*}
	is a quotient algebra of $\cH_d$ and a counital subcoalgebra of $(\cH_d, \triangle, \epsilon)$, making it a Hopf algebra too. The parameter $\gamma$ represents the truncation of Butcher-Connes-Kreimer expansion while $\alpha$ represents a choice of regularity. The redundancy of including $\alpha$ in this expression will become clear later.  
	
	\subsubsection{Coupled bialgebras}
	
	The purpose of this Subsection is to establish that the product and coupled coproduct operations interact with one another with the goal of introducing a coupled bialgebra. In order to do that, we need to introduce our own notion of $\mbox{Twist}$ operation that additionally keeps track of the couplings. Inspired by Equation \eqref{eq:Bialgebra-identity}, we introduce a twist operation for coupled tensor products:
	\begin{definition}
		Let $(\Omega', \cF', \bP')$ be a probability space and let $\rp$ be a 2-integrability functional. We denote the $\bR^e$ module
		\begin{align*}
			&\Big( \scH(\Omega') \otimes \scH(\Omega')\Big) \tilde{\otimes}_{(\rp, \rp)} \Big( \scH'(\Omega') \otimes \scH'(\Omega') \Big)
			\\
			&=\bigoplus_{(Y, \overline{Y}) \in \scF_0 \times \scF_0} L^{(\hat{p}[Y], \hat{p}[\overline{Y}])}\Bigg( (\Omega')^{\times|H^{Y}|} \times (\Omega')^{\times|H^{\hat{Y}}|}, (\bP')^{\times|H^{Y}|} \times (\bP')^{\times|H^{\hat{Y}}|}; 
			\\
			&\quad \bigoplus_{\substack{(\Upsilon, \overline{\Upsilon}) \in \scF_0 \times \scF_0 \\ G \in \lion(\Upsilon, Y) \\ \overline{G} \in \lion(\overline{\Upsilon}, \overline{Y}) }} L^{(\hat{p}[\Upsilon, G], \hat{p}[\overline{\Upsilon}, \overline{G}])}\bigg( (\Omega')^{\times (|G| - |H^{Y}|)} \times (\Omega')^{\times (|\overline{G}| - |H^{\overline{Y}}|)}, (\bP')^{\times (|G| - |H^{Y}|)} \times (\bP')^{\times (|\overline{G}| - |H^{\overline{Y}}|)}; 
			\\
			&\quad \lin\Big( (\bR^d)^{\otimes (|\scN^{\Upsilon}| + |\scN^{\overline{\Upsilon}}|)}, \lin\big( (\bR^d)^{\otimes (|\scN^{Y}| + |\scN^{\overline{Y}}|)}, \bR^e \big) \Big) \bigg) \Bigg)
		\end{align*}
		
		We define $\overline{\mbox{Twist}}$ to be the operator
		\begin{align*}
			\overline{\mbox{Twist}}:& \Big( \scH'(\Omega') \tilde{\otimes}_{\rp} \scH'(\Omega') \Big) \otimes \Big( \scH'(\Omega') \tilde{\otimes}_{\rp} \scH'(\Omega') \Big) 
			\\
			&\to \Big( \scH'(\Omega') \otimes \scH'(\Omega')\Big) \tilde{\otimes}_{(\rp, \rp)} \Big( \scH'(\Omega') \otimes \scH'(\Omega') \Big)
		\end{align*}
		that satisfies
		\begin{align*}
			&\overline{\mbox{Twist}}\Bigg[ \bigg( \sum_{\substack{\Upsilon \times^{G} Y \\ \in \scF_0 \tilde{\times} \scF_0 }} \Big\langle X, \Upsilon \times^{G} Y \Big\rangle(\omega_{\phi^{G, \Upsilon, Y}[G]} ) \bigg)
			\otimes \bigg( \sum_{\substack{\hat{\Upsilon} \times^{\hat{G}} \hat{Y} \\ \in \scF_0 \tilde{\times} \scF_0 }} \Big\langle \hat{X}, \hat{\Upsilon} \times^{\hat{G}} \hat{Y} \Big\rangle(\omega_{\phi^{\hat{G}, \hat{\Upsilon}, \hat{Y}}[\hat{G}]} ) \bigg) \Bigg]
			\\
			&= \sum_{\substack{(Y, \hat{Y}) \\ \in \scF_0 \times \scF_0 }} \sum_{\substack{(\Upsilon, \hat{\Upsilon}) \\ \in \scF_0 \times \scF \\ G \in \lion(\Upsilon, Y) \\ \hat{G} \in \lion(\hat{\Upsilon}, \hat{Y}) }} \mbox{Twist}\bigg[ \Big\langle X, \Upsilon\times^G Y \Big\rangle(\omega_{\phi^{G, \Upsilon, Y}[G]} ) \otimes \Big\langle \hat{X}, \hat{\Upsilon} \times^{\hat{G}} \hat{Y} \Big\rangle(\omega_{\phi^{\hat{G}, R, S}[\hat{G}]} ) \bigg]. 
		\end{align*}
		Finally, we can appropriately extend $\overline{\mbox{Twist}}$ to the $L^0(\Omega, \bP; \bR^e)$-modules
		\begin{align*}
			\overline{\mbox{Twist}}:& \Big( \scH\tilde{\otimes}_{\rp} \scH(\Omega, \Omega') \Big) \otimes \Big( \scH \tilde{\otimes}_{\rp} \scH(\Omega, \Omega') \Big) 
			\\
			&\to \Big( \scH(\Omega, \Omega') \otimes \scH(\Omega, \Omega') \Big) \tilde{\otimes}_{(\rp, \rp)} \Big( \scH(\Omega, \Omega') \otimes \scH(\Omega, \Omega') \Big). 
		\end{align*}
	\end{definition}
	Thus $\overline{\mbox{Twist}}$ has a similar nature to that of the $\mbox{Twist}$ operation, although we can see that it turns a pair of couples into a coupling between pairs. 

	\begin{proposition}
		\label{proposition:H-coupledBialgebra}
		Let $(\Omega, \cF, \bP)$ and $(\Omega', \cF', \bP')$ be probability spaces. Let $\rp$ be an integrability functional that satisfies Equation \eqref{eq:definition:AlgebraRVs-int+} and let $\hat{\rp}$ be a 2-integrability functional that satisfies Equation \eqref{eq:definition:coproduct2}. Let  $\scH_\rp'(\Omega')$ and $\scH_\rp(\Omega, \Omega')$ be the $\bR^e$-module and $L^0(\Omega, \bP; \bR^e)$-module respectively defined in Definition \ref{definition:AlgebraRVs}. 
		
		Then 
		\begin{align*}
			\big( \scH_\rp'(\Omega'), \circledast, \rId, \Delta, \epsilon \big)& \quad \mbox{is coupled bialgebras over the ring} \quad (\bR^e, +, \centerdot) \quad \mbox{and}
			\\
			\big( \scH_\rp(\Omega, \Omega'), \circledast, \rId, \Delta, \epsilon \big)& \quad \mbox{is coupled bialgebras over the ring} \quad \big( L^0(\Omega, \bP; \bR^e), +, \centerdot \big)
		\end{align*}
		In particular, $\big( \scH_\rp'(\Omega'), \circledast, \rId \big)$ and $\big( \scH_\rp(\Omega, \Omega'), \circledast, \rId \big)$ are associative algebras and $\big( \scH_\rp'(\Omega'), \Delta, \epsilon \big)$ and $\big( \scH_\rp(\Omega, \Omega'), \Delta, \epsilon \big)$ are co-associative coupled coalgebras and we have that
		\begin{equation}
			\label{eq:proposition:H-coupledBialgebra}
			\begin{aligned}
				\Delta \circ \circledast =& \Big( \circledast \tilde{\otimes} \circledast \Big) \circ \overline{\mbox{Twist}} \circ \Big( \Delta \otimes \Delta\Big)
				\\
				\epsilon \circ \circledast =& \centerdot \circ \epsilon \otimes \epsilon 
				\quad 
				\Delta \circ \rId \circ \centerdot = \rId \otimes \rId 
				\quad
				\epsilon \circ \rId = I
			\end{aligned}
		\end{equation}
	\end{proposition} 
	
	\begin{proposition}
		\label{proposition:Coproduct-Equivalence}
		The coproduct $\Delta: \scH_\rp(\Omega, \Omega') \to \scH \tilde{\otimes}_{\hat{\rp}} \scH(\Omega, \Omega')$ is the unique linear operator that satisfies the identities that
		\begin{equation*}
			\Delta[\rId] = \rId \times^{\emptyset} \rId, 
		\end{equation*}
		and $\forall T, T' \in \scF_0$, 
		\begin{equation}
			\label{eq:definition:coproduct}
			\begin{split}
				\Delta\Big[ \lfloor T\rfloor \Big] =& \lfloor T\rfloor \times^H \rId + \Big( I \times^H \lfloor \cdot \rfloor\Big)\circ \Delta \big[ T \big], 
				\qquad
				\Delta\Big[ \cE[T] \Big] = (\cE \tilde{\times} \cE) \Big[ \Delta [T] \Big], 
				\\
				\Delta\Big[ T \circledast T' \Big] =& \Big(\circledast \tilde{\otimes} \circledast \Big) \circ \overline{\mbox{Twist}} \circ \Big( \Delta \otimes \Delta\Big) \Big[ T \otimes T' \Big]. 
			\end{split}
		\end{equation}
	\end{proposition}
	
	The proof of Proposition \ref{proposition:H-coupledBialgebra} and Proposition \ref{proposition:Coproduct-Equivalence} along with many more details concerning this Section are found in \cite{salkeld2022LionsTrees}.
	
	\subsubsection{Gradings on Lions forests}
	
	Now we want to reproduce the graded structure of the Butcher-Connes-Kreimer setting. Suppose that $(\scJ, +)$ is a monoid. A coupled bialgebra $(\scH, \circledast, \rId, \Delta, \epsilon)$ over a ring $\cR$ is said to be $\scJ$-graded if it can be represented as
	\begin{equation*}
		\scH = \bigoplus_{j\in \scJ} \scH_{(j)}
	\end{equation*}
	and satisfies that
	\begin{equation}
		\label{eq:Grading-Definition}
		\scH_{(j_1)} \circledast \scH_{(j_2)} \subseteq \scH_{(j_1 + j_2)} 
		\quad \mbox{and} \quad
		\Delta\Big[ \scH_{(j)} \Big] \subseteq \bigoplus_{j_2 \in \scJ} \Big( \bigoplus_{\substack{j_1 \in \scJ \\ j_1 + j_2 = j}} \scH_{(j_1)}\Big) \tilde{\otimes} \scH_{(j_2)}. 
	\end{equation}
	Further, a graded bialgebra is said to be connected if for all units $0 \in \scJ$, we have $\scH_{(0)}= \cR$.

	The reader familiar with gradings on bialgebras will notice the asymmetry in Equation \eqref{eq:Grading-Definition}. This accounts for the asymmetry in the definition of the coupled tensor product (see \cite{salkeld2022LionsTrees}). 
	
	\begin{example}
		\label{example:regularitydiff}
		Let $N \in \bN$ and for each $i \in \{1, ..., N\}$ let $W^{i, N}:[0,1] \to \bR^d$ be a collection of independent Brownian motions on a probability space $(\Omega, \cF, \bP)$. It is widely known that for each $i\in \{1, ..., N\}$ and for any choice of $\alpha<\tfrac{1}{2}$, 
		\begin{equation}
			\label{eq:example:regularitydiff0}
			\sup_{s, t \in [0,1]} \frac{ \big| W_{s, t}^{i, N} \big|}{|t-s|^{\alpha}}< \infty \quad \bP-\mbox{almost surely. }
		\end{equation}
		Therefore, for any fixed choice of $N\in \bN$ we also have that
		\begin{equation}
			\label{eq:example:regularitydiff1}
			\sup_{s, t \in [0,1]} \frac{ \tfrac{1}{N} \sum_{j=1}^N \big| W_{s, t}^{j, N} \big|}{|t-s|^{\alpha}}< \infty \quad \bP-\mbox{almost surely. }
		\end{equation}
		However,
		\begin{equation}
			\label{eq:example:regularitydiff2}
			\sup_{s, t \in [0,1]} \frac{ \bE\Big[ \big| W_{s, t}^{i, N} \big| \Big] }{|t-s|^{1/2}}< \infty. 
		\end{equation}
		The contrast in regularity between Equation \eqref{eq:example:regularitydiff0} and Equation \eqref{eq:example:regularitydiff2} illustrates that is quite natural to suppose that the regularity of a driving signal may change depending on whether we are evaluating it on the tagged probability space or a free probability space. 
		
		Further, the fact that there is no change in the regularity between Equation \eqref{eq:example:regularitydiff0} and Equation \eqref{eq:example:regularitydiff1} shows us that we don't expect there to be any smoothing when we are working with systems of interacting equations as opposed to their mean-field limits. 
	\end{example}
	
	\begin{definition}
		\label{lemma:grading}
		We define $\scG:\scF_0 \to \bN_0^2$ such that for $T=(\scN^T, \scE^T, h_0^T, H^T, \scL^T)$, 
		\begin{equation}
			\label{eq:lemma:grading}
			\scG[T]:= \Big( |h_0^T|, |\scN^T\backslash h_0^T| \Big). 
		\end{equation}
	
		Let $(\Omega', \cF', \bP')$ be a probability space and let $\rp$ be an integrability functional that satisfies Equation \eqref{eq:definition:AlgebraRVs-int+}. For $(k, n) \in \bN_0^{\times 2}$, we define the $\bR^e$-module
		\begin{align*}
			\scH'_{\rp, (k, n)}\big( \Omega'\big) 
			:=& 
			\bigoplus_{\substack{T\in \scF_0 \\ \scG[T] = (k,n)}}
			L^{p[T]}\Big( (\Omega')^{\times |H^T|}, (\bP')^{\times |H^T|}; \lin \big( (\bR^d)^{\otimes |\scN^T|}, \bR^e\big) \Big). 
		\end{align*}
	
		Let $(\Omega, \cF, \bP)$ be a second probability space. For $(k, n) \in \bN_0^{\times 2}$, we define the $L^0(\Omega, \bP; \bR^e)$ module
		\begin{align*}
			\scH_{\rp, (k, n)}\big( \Omega, \Omega'\big) 
			:=& 
			L^0\bigg( \Omega, \bP; \bigoplus_{\substack{T\in \scF_0 \\ \scG[T] = (k,n)}}
			L^{p[T]}\Big( (\Omega')^{\times |H^T|}, (\bP')^{\times |H^T|}; \lin \big( (\bR^d)^{\otimes |\scN^T|}, \bR^e\big) \Big) \bigg).
		\end{align*}
	\end{definition}

	Thus for any $(k_1, n_1), (k_2, n_2) \in \bN_0^{\times 2}$ and $\hat{\rp}$ a 2-integrability functional that satisfies Equation \eqref{eq:definition:coproduct2}, we have that 
	\begin{align*}
		&\scH'_{\rp, (k_1, n_1)}(\Omega') \tilde{\otimes}_{\hat{\rp}} \scH'_{\rp, (k_2, n_2)}(\Omega')
		= \bigoplus_{\substack{Y \in \scF_0 \\ \scG[Y] = (k_2, n_2)}} L^{p[Y]}\Bigg( (\Omega')^{\times |H^Y|}, (\bP')^{\times |H^Y|}; 
		\\
		&\quad \bigoplus_{\substack{\Upsilon \in \scF_0 \\ \scG[\Upsilon] = (k_1, n_1) \\ G\in \lion(\Upsilon, Y)}}  L^{p[\Upsilon, G]} \bigg( (\Omega')^{\times (|G| - |H^Y|)}, (\bP')^{\times (|G| - |H^Y|)}; \lin\Big( (\bR^d)^{\otimes |\scN^{\Upsilon}|}, \lin\big( (\bR^d)^{\otimes |\scN^Y|}, \bR^e \big) \Big) \bigg) \Bigg). 
		\\
		&\scH_{\rp,(k_1, n_1)} \tilde{\otimes}_{\hat{\rp}} \scH_{\rp,(k_2, n_2)}(\Omega, \Omega') 
		= L^0\Bigg( \Omega, \bP; \bigoplus_{\substack{Y \in \scF_0 \\ \scG[Y] = (k_2, n_2)}} L^{p[Y]}\Bigg( (\Omega')^{\times |H^Y|}, (\bP')^{\times |H^Y|}; 
		\\
		&\quad \bigoplus_{\substack{\Upsilon \in \scF_0 \\ \scG[\Upsilon] = (k_1, n_1) \\ G\in \lion(\Upsilon, Y)}}  L^{p[\Upsilon, G]} \bigg( (\Omega')^{\times (|G| - |H^Y|)}, (\bP')^{\times (|G| - |H^Y|)}; \lin\Big( (\bR^d)^{\otimes |\scN^{\Upsilon}|}, \lin\big( (\bR^d)^{\otimes |\scN^Y|}, \bR^e \big) \Big) \bigg) \Bigg) \Bigg).
	\end{align*}
	
	\begin{proposition}
		\label{proposition:grading}
		Let $(\Omega, \cF, \bP)$ and $(\Omega', \cF', \bP')$ be probability spaces. Let $\rp$ be an integrability functional that satisfies Equation \eqref{eq:definition:AlgebraRVs-int+} and let $\hat{p}$ be a 2-integrability functional that satisfies Equation \eqref{eq:definition:coproduct2}. Let $\scH_\rp'(\Omega')$ and $\scH_\rp(\Omega, \Omega')$ be the $\bR^e$-module and $L^0(\Omega, \bP; \bR^e)$-module respectively defined in Definition \ref{definition:AlgebraRVs}. 
		
		Then the function $\scG:\scF_0 \to \bN_0^{\times 2}$ as defined in Equation \eqref{eq:lemma:grading} describes a $\bN_0^{\times 2}$-grading on the coupled bialgebras $\big( \scH_\rp'(\Omega'), \circledast, \rId, \Delta, \epsilon \big)$ and $\big( \scH_\rp(\Omega, \Omega'), \circledast, \rId, \Delta, \epsilon \big)$ and both coupled bialgebras are connected. 
		
		That is, 
		\begin{equation}
			\label{eq:lemma:grading_2}
			\begin{split}
				&\scH'_{\rp, (0,0)}(\Omega') = \bR^e,
				\\ 
				&\scH'_{\rp, (k_1, n_1)}(\Omega') \circledast \scH'_{\rp, (k_2, n_2)}(\Omega') \subseteq  \scH'_{\rp, (k_1+k_2, n_1+n_2)}(\Omega'),  
				\\
				&\Delta\Big[\scH'_{\rp, (k, n)}(\Omega')\Big] \subseteq  \bigoplus_{\substack{k', n' \in \bN_0}} \Big(  \scH'_{\rp, (k-k', n-n')}(\Omega') \Big) \tilde{\otimes}_{\hat{\rp}} \scH'_{\rp, (k', n')}(\Omega'), 
			\end{split}
		\end{equation}
		and
		\begin{equation}
			\label{eq:lemma:grading_3}
			\begin{split}
				&\scH_{\rp, (0,0)}(\Omega, \Omega') = L^0\big( \Omega, \bP; \bR^e \big), 
				\\
				&\scH_{\rp, (k_1, n_1)}(\Omega, \Omega') \circledast \scH_{\rp, (k_2, n_2)}(\Omega, \Omega') \subseteq  \scH_{\rp, (k_1+k_2, n_1+n_2)}(\Omega, \Omega'),  
				\\
				&\Delta\Big[\scH_{\rp, (k, n)}(\Omega, \Omega')\Big] \subseteq  \bigoplus_{\substack{k', n' \in \bN_0}} \Big(  \scH_{\rp, (k - k', n - n')}(\Omega, \Omega') \Big) \tilde{\otimes}_{\hat{\rp}} \scH_{\rp, (k', n')}(\Omega, \Omega'). 
			\end{split}
		\end{equation}
	\end{proposition}	
	The proof of Proposition \ref{proposition:grading} can be found in \cite{salkeld2022LionsTrees}. Although somewhat degenerate, we will often write 
	\begin{equation*}
		\scG_{\alpha, \beta}[T] = \alpha \cdot |h_0^T| + \beta \cdot |\scN\backslash h_0^T|
	\end{equation*}
	for $\alpha,\beta \in \bR^+$ and henceforth use the notation that for a set $\scA \subseteq \scF_0$, we will henceforward denote 
	\begin{equation}
		\label{eq:LionsTree-Set}
		\scA^{\gamma, \alpha, \beta} = \Big\{ T \in \scA: \scG_{\alpha, \beta}[T] \leq \gamma \Big\} 
		\quad \mbox{and} \quad
		\scA^{\gamma-, \alpha, \beta} = \Big\{ T \in \scA: \scG_{\alpha, \beta}[T] < \gamma \Big\}, 
	\end{equation}
	for some $\gamma\geq \alpha \wedge \beta$. We will frequently be studying summations that run over the elements of the sets \eqref{eq:LionsTree-Set} so to lighten the notation somewhat we will often denote
	\begin{equation*}
		\sum_{T\in\scA}^{\gamma, \alpha, \beta}
		:= 
		\sum_{\substack{T\in \scA\\ \scG_{\alpha, \beta}[T]\leq \gamma}}
		\quad\mbox{and}\quad 
		\sum_{T\in\scA}^{\gamma-, \alpha, \beta}
		:= 
		\sum_{\substack{T\in \scA\\ \scG_{\alpha, \beta}[T] < \gamma}}, 
	 \end{equation*}
 	
 	\begin{remark}
 		The reader may already have noticed that we could simplify sets such as Equation \eqref{eq:LionsTree-Set} by removing the value $\gamma$ and re-expressing these sets as
 		\begin{equation*}
	 		\scA^{\alpha/\gamma, \beta/\gamma} = \Big\{ T\in \scA: \scG_{\tfrac{\alpha}{\gamma}, \tfrac{\beta}{\gamma}}[T] \leq 1 \Big\}. 
 		\end{equation*}
 		Although this is more streamlined, we found that the choices of the values of $\alpha$ and $\beta$ became very unintuitive in later portions of this project and so the current notation remains. For this work, $\alpha$ represents the $\bP$-almost sure H\"older regularity of the probabilistic rough path while $\beta$ represents the mean-square regularity of the probabilistic rough path (Example \ref{example:regularitydiff} shows that these may be different)
 	\end{remark}
	
	\subsubsection{Ideals and Lions forests}
	
	The purpose of the algebraic structures that we have introduced so far is to better understand an abstract Lions-Taylor expansion. In practice, we never want to work with the infinite expansions so that we need to verify that we can truncate the expansions considered thus far and obtain something that retains the algebra/coupled coalgebra structure. It is important to expose the full details of these infinite expansions as it allows us to verify that the algebraic and combinatorial structures are valid but in practice we want to perform a truncation. 
	
	However, as explained in \cite{salkeld2022Lions}, our truncation needs to take into account the number of spacial as well as the number of Lions derivatives. We have already seen that our coupled bialgebra has a grading on the monoid of pairs of integers so that our perspective should be to truncate with respect to this grading. For notational purposes, we use $\alpha$ and $\beta$ to represent the ``regularity'' of terms on the tagged probability space and terms on detagged probability spaces respectively. The parameter $\gamma$ represents the choice of truncation point. Our next task is to verify that truncating these expansions can be performed meaningfully:
	\begin{definition}
		Let $(\Omega, \cF, \bP)$ and $(\Omega', \cF', \bP')$ be probability spaces and let $\rp$ be an integrability functional that satisfies Equation \eqref{eq:definition:AlgebraRVs-int+}. For $\alpha, \beta>0$ and $\gamma\geq \alpha \wedge \beta$, we define
		\begin{align*}
			\scI'(\Omega')_\rp^{\gamma, \alpha, \beta} = \bigoplus_{\substack{(i,j)\in \bN_0^2\\ \alpha i+ \beta j>\gamma}}^\infty \scH'_{\rp, (i, j)}(\Omega')
			\qquad& \mbox{and} \qquad
			\scI(\Omega, \Omega')_\rp^{\gamma, \alpha, \beta} = \bigoplus_{\substack{(i,j)\in \bN_0^2\\ \alpha i+ \beta j>\gamma}}^\infty \scH_{\rp, (i, j)}(\Omega, \Omega'). 
		\end{align*}
		
		We denote 
		\begin{equation*}
			\scH_\rp'(\Omega')^{\gamma, \alpha, \beta} = \scH_\rp'(\Omega') / \scI'(\Omega')_\rp^{\gamma, \alpha, \beta}
			\quad \mbox{and}\quad
			\scH_\rp^{\gamma, \alpha, \beta}(\Omega, \Omega') = \scH_\rp(\Omega, \Omega') / \scI'(\Omega')_\rp^{\gamma, \alpha, \beta}
		\end{equation*}
		as the quotient algebras (over the ring $\bR^e$ and $L^0(\Omega, \bP; \bR^e)$ respectively). 
	\end{definition}
	
	\begin{proposition}
		\label{proposition:Finite-grading}
		Let $\alpha, \beta>0$ and let $\gamma\geq \alpha\wedge \beta$. Let $(\Omega, \cF, \bP)$ and $(\Omega', \cF, \bP')$ be probability spaces. Let $\rp$ be an integrability functional that satisfies Equation \eqref{eq:definition:AlgebraRVs-int+} and let $\hat{\rp}$ be a 2-integrability functional that satisfies Equation \eqref{eq:definition:coproduct2}. Then 
		\begin{align*}
			\big( \scH_\rp'(\Omega')^{\gamma, \alpha, \beta}, \Delta, \epsilon \big)
			\quad& \mbox{is a co-associative sub-coupled coalgebras of} \quad 
			\big( \scH_\rp'(\Omega'), \Delta, \epsilon \big)
			\quad \mbox{and}
			\\
			\big( \scH_\rp^{\gamma, \alpha, \beta}(\Omega, \Omega'), \Delta, \epsilon \big)
			\quad& \mbox{is a co-associative sub-coupled coalgebras of} \quad
			\big( \scH_\rp(\Omega, \Omega'), \Delta, \epsilon \big).
		\end{align*} 
	\end{proposition}
	The proof of Proposition \ref{proposition:Finite-grading} can be found in \cite{salkeld2022LionsTrees}.
	
	\subsection{Group structures of a coupled bialgebra}
	\label{subse:2.6}
	
	As explained in \cite{connes1999hopf}, the Butcher-Connes-Kreimer Hopf algebra is a module that has characters taking values in the Butcher group, an infinite dimensional Lie group first identified in \cite{Hairer1974Butcher} in the context of approximating the solutions of non-linear differential equations. Formally, the \emph{Butcher group of characters} is the subset of
	\begin{equation*}
		\cG\big( \cH_d^{\gamma, \alpha}, \bR^e\big) \subseteq \bigoplus_{\tau \in \fF_0}^{\gamma, \alpha} (\bR^d)^{\otimes |\scN^T|}
	\end{equation*}
	that satisfy the character identity that $f\in G\big( \cH_d^{\gamma, \alpha}, \bR^e \big)$ if and only if
	\begin{equation*}
		\Big\langle f, \tau_1 \odot \tau_2 \Big\rangle = \Big\langle f, \tau_1 \Big\rangle \otimes \Big\langle f, \tau_2 \Big\rangle. 
	\end{equation*}

	More generally, let us recall that given a bialgebra $\big(\cH, \odot, \rId, \triangle, \epsilon \big)$ over a ring $(\cR, +, \centerdot)$ (such as but not specifically the Butcher-Connes-Kreimer Hopf algebra, for a reference see \cite{cartier2021hopf}) and a commutative algebra $(\cA, m_{\cA})$ over the ring $(\cR, +, \centerdot)$, we define the set of \emph{characters} to be
	\begin{equation*}
		\cG(\cH, \cA) = \Big\{ f\in \lin(\cH, \cA): \quad f\circ \odot = m_{\cA} \circ f \otimes f \Big\}. 
	\end{equation*}
	For two mappings $f,g\in \lin(\cH, \cA)$, we define the convolution product $\ast: \lin(\cH, \cA) \times \lin(\cH, \cA) \to \lin(\cH, \cA)$ by
	\begin{equation*}
		f \ast g = m_{\cA} \circ f \otimes g \circ \triangle
	\end{equation*}
	and note that it is well understood that $\big( \cG(\cH, \cA), \ast \big)$ is a monoid. 
	
	From classical theory of bialgebras we know that a bialgebra is graded if and only if it is a Hopf algebra, that is there exists an anti-automorphism $\cS$ such that
	\begin{equation}
		\label{eq:HopfAlgebraAntipode}
		\begin{aligned}
			&\odot \circ (\cS \otimes I) \circ \triangle = \odot \circ (I \otimes \cS) \circ \triangle = \rId \epsilon 
			\\
			&\mbox{and} \quad \forall n\in \bN_0 \quad \cS\big[ \cH_{d, n} \big] \subseteq \cH_{d, n}. 
		\end{aligned}
	\end{equation}
	When the bialgebra has an antipode, then an inverse operation exists and $\big( \cG(\cH, \cA), \ast \big)$ is a group since
	\begin{equation*}
		\Big( f \ast (f\circ \cS) \Big) [x] = f \Big[ \odot \circ (I \otimes \cS) \circ \triangle \big[ x \big] \Big] 
		=  
		f\big[ \rId \big] \epsilon[x] 
		= 
		\Big( (f\circ \cS) \ast f \Big) [x] 
	\end{equation*}
	Our goal is to extend all these concepts to the coupled coproduct. 
	
	\subsubsection{Dual spaces for coupled bialgebras}
	
	Let $(\Omega', \cF', \bP')$ be a probability space, let $\rp$ be an integrability functional and let $\rQ = (q[T])_{T\in \scF}$ be an integrability functional such that for all $T\in \scF$ and
	\begin{equation*}
		\forall h \in H^T \quad \frac{1}{p[T]_h} + \frac{1}{q[T]_h} = 1. 
	\end{equation*}
	
	Consider the tensor series modules 
	\begin{equation}
		\label{eq:Dual-LionTrees}
		\big(\scH_\rp'(\Omega') \big)^{\dagger} =
		\prod_{T\in \scF_0} L^{q[T]}\Big((\Omega')^{\times |H^T|}, (\bP')^{ \times |H^T|};  (\bR^d)^{\otimes |\scN^T|} \Big). 
	\end{equation}
	There is a natural choice of bilinear mapping
	\begin{equation}
		\label{eq:nonInt_BilinearForm2}
		\begin{split}
			\Big\langle \cdot, \cdot \Big\rangle:& \scH_\rp'(\Omega') \times \big( \scH_\rp'(\Omega') \big)^\dagger \to \bR^e
		\end{split}
	\end{equation}
	defined for $X\in \scH_\rp'(\Omega')$ and $W\in \big( \scH_\rp'(\Omega')\big)^{\dagger}$ by
	\begin{align*}
		\Big\langle X, W \Big\rangle = \sum_{T\in \scF_0} (\bE')^{\times |H^T|} \bigg[ \Big\langle X, T\Big\rangle(\omega'_{H^T}) \cdot \Big\langle W, T \Big\rangle(\omega'_{H^T}) \bigg]
	\end{align*}
	so that the module \eqref{eq:Dual-LionTrees} can be thought of as dual spaces and topologised accordingly using the projective topology induced by the bilinear forms. 
	
	Let $(\Omega, \cF, \bP)$ be a separate probability space and consider the tensor series modules
	\begin{equation}
		\label{eq:Dual-LionTrees-}
		\big(\scH_\rp(\Omega, \Omega') \big)^{\dagger} =
		L^0\bigg( \Omega, \bP; \prod_{T\in \scF_0} L^{q[T]}\Big((\Omega')^{\times |H^T|}, (\bP')^{ \times |H^T|};  (\bR^d)^{\otimes |\scN^T|} \Big) \bigg). 
	\end{equation}
	The bilinear mappings defined in \eqref{eq:nonInt_BilinearForm2} can be extended $\bP$-almost surely to give alternative bilinear mappings
	\begin{equation}
		\label{eq:nonInt_BilinearForm3}
		\begin{split}
			\Big\langle \cdot, \cdot \Big\rangle:& \scH_\rp(\Omega, \Omega') \times \big( \scH_\rp(\Omega, \Omega') \big)^\dagger \to L^0(\Omega, \bP; \bR^e)
		\end{split}
	\end{equation}
	defined for $X\in \scH_\rp'(\Omega')$ and $W\in \big( \scH_\rp'(\Omega')\big)^{\dagger}$ by
	\begin{align*}
		\Big\langle X, W \Big\rangle(\omega_0) = \sum_{T\in \scF_0} (\bE')^{\times |H^T|} \bigg[ \Big\langle X, T\Big\rangle(\omega_0, \omega'_{H^T}) \cdot \Big\langle W, T \Big\rangle(\omega_0, \omega'_{H^T}) \bigg]
	\end{align*}

	Finally, we topologise \eqref{eq:Dual-LionTrees-} according to the projective topology induced by the bilinear forms \eqref{eq:nonInt_BilinearForm3}. This detail is so important because it ensures that we want a sense of convergence in mean on the detagged probability spaces $(\Omega', \cF', \bP)$ but on the tagged probability space $(\Omega, \cF, \bP)$ we want a sense of $\bP$-almost sure convergence. 
	
	Let $\gamma, \alpha, \beta>0$ and suppose that $\gamma> \alpha\vee \beta$. Then in the same fashion
	\begin{equation}
		\label{eq:Dual-LionTrees+}
		\Big(\scH_\rp'(\Omega')^{\gamma, \alpha, \beta} \Big)^{\dagger} =
		\bigoplus_{T\in \scF_0}^{\gamma, \alpha, \beta} L^{q[T]}\Big((\Omega')^{\times |H^T|}, (\bP')^{ \times |H^T|};  (\bR^d)^{\otimes |\scN^T|} \Big). 
	\end{equation}
	Following on from the bilinear forms described in Equation \eqref{eq:nonInt_BilinearForm2}, we have that
	\begin{equation}
		\label{eq:nonInt_BilinearForm}
		\begin{split}
			\Big\langle \cdot, \cdot \Big\rangle:& \scH_\rp'(\Omega')^{\gamma, \alpha, \beta} \times \Big( \scH_\rp'(\Omega')^{\gamma, \alpha, \beta} \Big)^\dagger \to \bR^e
		\end{split}
	\end{equation}
	By extending to include the tagged probability space, finally we have
	\begin{equation}
		\label{eq:Dual-LionTrees+-}
		\Big(\scH_\rp^{\gamma, \alpha, \beta}(\Omega, \Omega') \Big)^{\dagger} = L^0\bigg( \Omega, \bP; 
		\bigoplus_{T\in \scF_0}^{\gamma, \alpha, \beta} L^{q[T]}\Big((\Omega')^{\times |H^T|}, (\bP')^{ \times |H^T|};  (\bR^d)^{\otimes |\scN^T|} \Big) \bigg)
	\end{equation}
	which (in the same fashion as Equation \eqref{eq:nonInt_BilinearForm3}) have bilinear form
	\begin{equation}
		\label{eq:nonInt_BilinearForm4}
		\begin{split}
			\Big\langle \cdot, \cdot \Big\rangle:& \scH_\rp^{\gamma, \alpha, \beta}(\Omega, \Omega') \times \Big( \scH_\rp^{\gamma, \alpha, \beta}(\Omega, \Omega') \Big)^\dagger \to L^0(\Omega, \bP; \bR^e)
		\end{split}
	\end{equation}

	\begin{remark}
		We have not taken any care in denoting the many different bilinear forms defined in Equation \eqref{eq:nonInt_BilinearForm2}, \eqref{eq:nonInt_BilinearForm3}, \eqref{eq:nonInt_BilinearForm} and \eqref{eq:nonInt_BilinearForm4} and will use these interchangeably where the context is clear.  
	\end{remark}
	
	\subsubsection{Characters of coupled bialgebras}
	
	As a first step, recalling Equations \eqref{eq:Dual-LionTrees} and \eqref{eq:Dual-LionTrees+} and following on from the bilinear forms Equation \eqref{eq:nonInt_BilinearForm2} and \eqref{eq:nonInt_BilinearForm} we note that
	\begin{align*}
		\lin\big( \scH_\rp'(\Omega'), \bR^e \big) \equiv \big( \scH_\rp'(\Omega') \big)^{\dagger}
		&\quad \mbox{and} \quad 
		\lin\Big( \scH_\rp'(\Omega')^{\gamma, \alpha, \beta}, \bR^e \Big) \equiv \Big( \scH_\rp'(\Omega')^{\gamma, \alpha, \beta} \Big)^{\dagger}. 
	\end{align*}
	Let $(\Omega, \cF, \bP)$ be a second probability space and let $U$ be a module over a ring $(\cR, +, \centerdot)$. Then we can think of $L^0(\Omega, \bP; U)$ as a module over the ring $\big( L^0(\Omega, \bP; \cR), +, \centerdot \big)$ with the ring product being defined $\bP$-almost surely. In contrast to the more commonly studied case where $L^0(\Omega, \bP; U)$ is a module over the ring $(\cR, +, \centerdot)$, we get that
	\begin{equation*}
		\lin\Big( L^0\big( \Omega, \bP; U \big), L^0\big( \Omega, \bP; \cR \big) \Big) \equiv L^0\Big(\Omega, \bP; \lin\big( U, \cR \big) \Big) \quad \mbox{over the ring} \quad L^0(\Omega, \bP; \cR). 
	\end{equation*}

	Recalling Equations \eqref{eq:Dual-LionTrees-} and \eqref{eq:Dual-LionTrees+-} and following on from a bilinear forms from Equations \eqref{eq:nonInt_BilinearForm3} and \eqref{eq:nonInt_BilinearForm4}, we get that
	\begin{equation}
		\label{eq:LinMap-onto-measurable}
		\left.
		\begin{aligned}
			\lin\Big( \scH_\rp(\Omega, \Omega'), L^0(\Omega, \bP; \bR^e) \Big) \equiv& \big( \scH_\rp(\Omega, \Omega') \big)^{\dagger} \quad \mbox{and} \qquad \qquad 
			\\
			\lin\Big( \scH_\rp^{\gamma, \alpha, \beta}(\Omega, \Omega'), L^0(\Omega, \bP; \bR^e) \Big) \equiv& \Big( \scH_\rp^{\gamma, \alpha, \beta}(\Omega, \Omega') \Big)^{\dagger}. 
		\end{aligned}
		\right\}
	\end{equation}
	
	\begin{definition}
		\label{definition:set-characters}
		Let $(\Omega', \cF', \bP')$ be a probability space and let $\rp$ be an integrability functional that satisfies Equation \eqref{eq:definition:AlgebraRVs-int+}. We define 
		\begin{align*}
			&\cG\big( \scH_\rp'(\Omega'), \bR^e \big) = \bigg\{ f \in \big(\scH_\rp'(\Omega') \big)^{\dagger} : \forall T_1, T_2 \in \scF, 
			\\
			&\quad \Big\langle f, T_1 \circledast T_2 \Big\rangle(\omega'_{H^{T_1}}, \omega'_{H^{T_2}})  = \Big\langle f, T_1 \Big\rangle(\omega'_{H^{T_1}}) \otimes \Big\langle f, T_2 \Big\rangle(\omega'_{H^{T_2}}) \quad (\bP')^{\times |H^{T_1 \circledast T_2}|}\mbox{-almost surely} \bigg\}. 
		\end{align*}
		Let $\alpha, \beta>0$ and let $\gamma>\alpha\wedge\beta$. We define
		\begin{align*}
			&\cG\big( \scH_\rp'(\Omega')^{\gamma, \alpha, \beta}, \bR^e \big) = \bigg\{ f \in \big(\scH_\rp'(\Omega')^{\gamma, \alpha, \beta} \big)^{\dagger} : \forall T_1, T_2 \in \scF \quad \mbox{s.t} \quad T_1 \circledast T_2 \in \scF^{\gamma, \alpha, \beta}, 
			\\
			&\quad \Big\langle f, T_1 \circledast T_2 \Big\rangle(\omega'_{H^{T_1}}, \omega'_{H^{T_2}})  = \Big\langle f, T_1 \Big\rangle(\omega'_{H^{T_1}}) \otimes \Big\langle f, T_2 \Big\rangle(\omega'_{H^{T_2}}) \quad (\bP')^{\times |H^{T_1 \circledast T_2}|}\mbox{-almost surely} \bigg\}. 
		\end{align*}
		Let $(\Omega, \cF, \bP)$ be a second probability space. Following on from Equation \eqref{eq:LinMap-onto-measurable}, we also define
		\begin{align*}
			\cG\Big( \scH_\rp(\Omega, \Omega'), &L^0(\Omega, \bP; \bR^e) \Big) = \bigg\{ f \in \big(\scH_\rp(\Omega, \Omega') \big)^{\dagger} : \forall T_1, T_2 \in \scF_0
			\\
			&\Big\langle f, T_1 \circledast T_2 \Big\rangle(\omega_0, \omega'_{H^{T_1}}, \omega'_{H^{T_2}}) = \Big\langle f, T_1 \Big\rangle(\omega_0, \omega'_{H^{T_1}}) \otimes \Big\langle f, T_2 \Big\rangle(\omega_0, \omega'_{H^{T_2}}) \\
			&\quad \bP \times (\bP')^{\times |H^{T_1 \circledast T_2}|}\mbox{-almost surely} \bigg\} \quad \mbox{and}
			\\
			\cG\Big( \scH_\rp^{\gamma, \alpha, \beta}(\Omega, \Omega'), &L^0(\Omega, \bP; \bR^e) \Big) = \bigg\{ f \in \big(\scH_\rp^{\gamma, \alpha, \beta}(\Omega, \Omega') \big)^{\dagger} : \forall T_1, T_2 \in \scF_0 \quad \mbox{s.t.} \quad T_1 \circledast T_2 \in \scF^{\gamma, \alpha, \beta}
			\\
			&\Big\langle f, T_1 \circledast T_2 \Big\rangle(\omega_0, \omega'_{H^{T_1}}, \omega'_{H^{T_2}}) = \Big\langle f, T_1 \Big\rangle(\omega_0, \omega'_{H^{T_1}}) \otimes \Big\langle f, T_2 \Big\rangle(\omega_0, \omega'_{H^{T_2}}) 
			\\
			&\quad \bP \times (\bP')^{\times |H^{T_1 \circledast T_2}|}\mbox{-almost surely} \bigg\}. 
		\end{align*}
	\end{definition}

	\begin{remark}
		\label{remark:Integrability-1}
		From the expansions in Equation \eqref{eq:Dual-LionTrees-}, we would expect that every character satisfies that for any two Lions trees $\Upsilon$ and $Y$, 
		\begin{equation*}
			\Big\langle f, \Upsilon \circledast Y\Big\rangle(\omega_0, \omega'_{H^\Upsilon}, \omega'_{H^Y}) = \Big\langle f, \Upsilon \Big\rangle(\omega_0, \omega'_{H^{\Upsilon}}) \otimes \Big\langle f, Y \Big\rangle(\omega_0, \omega'_{H^Y})
		\end{equation*}
		so that when we integrate we get
		\begin{equation*}
			(\bE')^{H^{Y\circledast \Upsilon}}\bigg[ \Big| \big\langle f, \Upsilon \circledast Y \big\rangle(\omega_0, \omega'_{H^\Upsilon}, \omega'_{H^Y}) \Big| \bigg]
			=
			(\bE)^{H^\Upsilon}\bigg[ \Big| \big\langle f, \Upsilon \big\rangle(\omega_0, \omega'_{H^\Upsilon}) \Big| \bigg]
			\otimes 
			(\bE)^{H^Y}\bigg[ \Big| \big\langle f, Y \big\rangle(\omega_0, \omega'_{H^Y}) \Big| \bigg].			
		\end{equation*}
		As such, even though we are considering the product of random variables the probability spaces over which we integrate are all orthogonal so that there are no concerns about the integrability of the product of random variables. 
	\end{remark}

	In Definition \ref{definition:AlgebraRVs}, we introduced integrability functionals as a way of capturing the different properties of each probability space in our module indexed by Lions forests. Finding the integrability functional that captures the integrability properties of the associated dual space is a relatively simple task, but we need to take care that the convolution product on the dual space (the dual of the coproduct) also remains meaningful:
	\begin{definition}
		\label{definition:integrability-functional}
		Let $\alpha,\beta>0$, let $\gamma > \alpha \wedge \beta$ and let $\scF^{\gamma, \alpha, \beta} = \big\{ T\in \scF: \scG_{\alpha, \beta}[T]\leq \gamma \big\}$. Let 
		\begin{equation*}
			\rp = \big( p[T] \big)_{T\in \scF^{\gamma, \alpha, \beta}} \in \bigsqcup_{T\in \scF}^{\gamma, \alpha, \beta} (1, \infty)^{\times |H^T|}
			\quad \mbox{and} \quad
			\rQ = \big( q[T] \big)_{T\in \scF^{\gamma, \alpha, \beta}} \in \bigsqcup_{T\in \scF}^{\gamma, \alpha, \beta} (1, \infty)^{\times |H^T|}
		\end{equation*}
		so that $\rp$ and $\rQ$ are integrability functionals and suppose that for all $T\in \scF^{\gamma, \alpha, \beta}$ and
		\begin{equation}
			\label{eqdefinition:integrability-functional.3}
			\forall h \in H^T,  \quad \frac{1}{p[T]_h} + \frac{1}{q[T]_h} = 1. 
		\end{equation}
		Then we have that $\rp$ satisfies Equation \eqref{eq:definition:AlgebraRVs-int+} if and only if $\rQ$ satisfies Equation \eqref{eq:definition:AlgebraRVs-int+}. 
		
		Next, suppose that the $\rQ$ satisfies that for all $T, \Upsilon, Y \in \scF_0$ such that $c(T, \Upsilon, Y)>0$, 
		\begin{equation}
			\label{eq:definition:integrability-functional.2}
			\left.
			\begin{aligned}
				&\forall h\in H^T \quad \mbox{s.t.} \quad h\cap \scN^\Upsilon = \emptyset &\implies& \quad \frac{1}{q[T]_h} \geq \frac{1}{q[Y]_{\phi^{T, \Upsilon, Y}[h]}}, 
				\\
				&\forall h\in H^T \quad \mbox{s.t.} \quad h\cap \scN^\Upsilon \neq \emptyset, h\cap \scN^Y \neq \emptyset  &\implies& \quad \frac{1}{q[T]_h} \geq \frac{1}{q[\Upsilon]_{h\cap \scN^{\Upsilon}}} + \frac{1}{q[Y]_{\phi^{T, \Upsilon, Y}[h]}}, 
				\\
				&\forall h\in H^T \quad \mbox{s.t.} \quad h\cap \scN^\Upsilon \neq \emptyset, h\cap \scN^Y = \emptyset  &\implies& \quad \frac{1}{q[T]_h} \geq \frac{1}{q[\Upsilon]_{h\cap \scN^{\Upsilon}}}. 
			\end{aligned}
			\right\}
		\end{equation}
	
		We say that $(\rp, \rQ)$ is called a \emph{dual integrability functional} if $\rp$ and $\rQ$ are integrability functionals that satisfy Equation \eqref{eq:definition:AlgebraRVs-int+}, Equation \eqref{eqdefinition:integrability-functional.3} and Equation \eqref{eq:definition:integrability-functional.2}. 
	\end{definition}
	
	\begin{example}
		\label{example:integrability-functional}
		Recall from \cite{2019arXiv180205882.2B} that the authors work with the module
		\begin{align*}
			L^0\bigg( \Omega, \bP;& \quad \bR^e \oplus L^q\Big( \Omega, \bP; \bR^d\Big) 
			\\
			&\oplus (\bR^d)^{\otimes 2} \oplus L^{q/2}\Big( \Omega, \bP; (\bR^d)^{\otimes 2} \Big) \oplus L^{q/2}\Big( \Omega, \bP; (\bR^d)^{\otimes 2} \Big) \oplus L^{q/2}\Big( \Omega^{\times 2}, \bP^{\times 2}; (\bR^d)^{\otimes 2} \Big) \bigg). 
		\end{align*}
		As stated following their definition of the $\omega$-controlled path (\cite{2019arXiv180205882.2B}*{page 12}) the choice of $q\geq 8$ is somewhat arbitrary but necessary. The purpose of the integrability functional introduced in Definition \ref{definition:integrability-functional} is to address this need for higher order expansions. 
		
		For instance, given $\alpha, \beta>0$ and $\gamma>\alpha\wedge \beta$ we could simply choose $n\in \bN$ such that 
		\begin{equation*}
			n > 2 \cdot \sup_{T\in \scF^{\gamma, \alpha, \beta}} \big| \scN^T \big|
		\end{equation*}
		and define $\rQ = \big( (q[T]_h)_{h\in H^T} \big)_{T\in \scF}$ by
		\begin{equation}
			\label{eq:example:integrability-functional}
			q[T] = \Big( \tfrac{n}{|h|} \Big)_{h\in H^T}. 
		\end{equation}
		We can verify that this integrability functional satisfies Equation \eqref{eq:definition:integrability-functional.2} and Equation \eqref{eq:definition:AlgebraRVs-int+}. By defining the integrability function $\rp$ according to Equation \eqref{eqdefinition:integrability-functional.3}, we obtain that $(\rp, \rQ)$ is a dual integrability functional. 
	\end{example}
	
	The bilinear forms from Equations \eqref{eq:nonInt_BilinearForm2}, \eqref{eq:nonInt_BilinearForm3}, \eqref{eq:nonInt_BilinearForm} and \eqref{eq:nonInt_BilinearForm4} provide a framework by which the coproduct $\Delta$ (which is defined on all of these modules) can induce a binary operation onto the dual spaces:
	\begin{definition}
		Let $(\Omega', \cF', \bP')$ be a probability space and let $(\rp, \rQ)$ be a dual integrability functional. Let $f, g\in \lin\big( \scH_\rp'(\Omega'), \bR^e \big)$. Then we define the convolution product $f\ast g$ according to the identity
		\begin{equation}
			\label{eq:definition:DualModule_ConvolProd}
			f \ast g = \centerdot \circ f \tilde{\otimes} g \circ \Delta. 
		\end{equation}
		That is, 
		\begin{align*}
			&f\ast g \in \bigoplus_{T \in \scF_0} L^{q[T]}\Big( (\Omega')^{\times |H^T|}, (\bP')^{\times |H^T|}; (\bR^d)^{\otimes |\scN^T|} \Big)
			\\
			\Big\langle &f \ast g, T\Big\rangle(\omega'_{H^T}) = \sum_{\Upsilon, Y \in \scF_0} c\Big( T, \Upsilon, Y \Big) \cdot \Big\langle f, \Upsilon \Big\rangle \otimes^{H^T} \Big\langle g, Y \Big\rangle (\omega'_{H^T}) \quad (\bP')^{\times |H^T|}\mbox{-almost surely}. 
		\end{align*}
		Further, let $\alpha, \beta>0$, let $\gamma>\alpha\wedge\beta$ and let $f, g\in \lin\big( \scH_{\rp}'(\Omega')^{\gamma, \alpha, \beta}, \bR^e \big)$. Then the convolution product $f\ast g$ satisfies Equation \eqref{eq:definition:DualModule_ConvolProd} and
		\begin{equation*}
			f\ast g \in \bigoplus_{T \in \scF_0}^{\gamma, \alpha, \beta} L^{q[T]}\Big( (\Omega')^{\times |H^T|}, (\bP')^{\times |H^T|}; (\bR^d)^{\otimes |\scN^T|} \Big). 
		\end{equation*}
	
		Let $(\Omega, \cF, \bP)$ be another probability space and now let $f, g \in \lin\big( \scH_{\rp}(\Omega, \Omega'), L^0(\Omega, \bP; \bR^e) \big)$. Then we extend the convolution product identity \eqref{eq:definition:DualModule_ConvolProd} to be
		\begin{align*}
			&f\ast g \in L^0\bigg( \Omega, \bP; \bigoplus_{T \in \scF_0} L^{q[T]}\Big( (\Omega')^{\times |H^T|}, (\bP')^{\times |H^T|}; (\bR^d)^{\otimes |\scN^T|} \Big) \bigg)
			\\
			&\Big\langle f \ast g, T\Big\rangle(\omega_0, \omega'_{H^T}) 
			\\
			&\qquad = \sum_{\Upsilon, Y \in \scF_0} c\Big( T, \Upsilon, Y \Big) \cdot \Big\langle f, \Upsilon \Big\rangle \otimes^{H^T} \Big\langle g, Y \Big\rangle (\omega_0, \omega'_{H^T}) \quad \bP \times (\bP')^{\times |H^T|}\mbox{-almost surely.} 
		\end{align*}
		Finally, for $f, g \in \lin\big( \scH_\rp^{\gamma, \alpha, \beta}(\Omega, \Omega'), L^0(\Omega, \bP; \bR^e) \big)$, we have the convolution product $f \ast g$ satisfies Equation \eqref{eq:definition:DualModule_ConvolProd} and
		\begin{equation*}
			f\ast g \in L^0\bigg( \Omega, \bP; \bigoplus_{T \in \scF_0}^{\gamma, \alpha, \beta} L^{q[T]}\Big( (\Omega')^{\times |H^T|}, (\bP')^{\times |H^T|}; (\bR^d)^{\otimes |\scN^T|} \Big) \bigg). 
		\end{equation*}
	\end{definition}
	
	In fact, with this work we only concern ourselves with the convolution product of characters:
	\begin{proposition}
		\label{proposition:Characters1!}
		Let $(\Omega, \cF, \bP)$ and $(\Omega', \cF', \bP')$ be probability spaces. Let $\alpha, \beta>0$ and let $\gamma>\alpha\wedge\beta$. Let $(\rp, \rQ)$ be a dual integrability functional. 

		Then
		\begin{equation}
			\label{eq:proposition:Characters!1}
			f, g\in \cG\big( \scH_\rp^{\gamma, \alpha, \beta}(\Omega, \Omega'), L^0(\Omega, \bP; \bR^e) \big) 
			\quad \implies \quad
			f \ast g \in \cG\big( \scH_{\rp}^{\gamma, \alpha, \beta}(\Omega, \Omega'), L^0(\Omega, \bP; \bR^e) \big)
		\end{equation}
		and
		\begin{equation}
			\label{eq:proposition:Characters!3}
			f\in \cG\big( \scH_\rp^{\gamma, \alpha, \beta}(\Omega, \Omega'), L^0(\Omega, \bP; \bR^e) \big) 
			\quad \implies \quad
			f \ast \epsilon = \epsilon \ast f = f. 
		\end{equation}
	\end{proposition}
	
	\begin{proof}
		Let us start by proving Equation \eqref{eq:proposition:Characters!1}: Let $f, g \in \cG\big( \scH_\rp^{\gamma, \alpha, \beta}(\Omega, \Omega'), L^0(\Omega, \bP; \bR^e) \big)$. 
		
		Using Proposition \ref{proposition:H-coupledBialgebra}, we have for any $T_1, T_2\in \scF_0$ such that $T_1 \circledast T_2 \in \scF_0^{\gamma, \alpha, \beta}$ that
		\begin{align*}
			\Big\langle f \ast g, T_1 \circledast T_2& \Big\rangle(\omega_0, \omega'_{H^{T_1 \circledast T_2}}) 
			= \Big\langle f \tilde{\otimes} g, \Delta\big[ T_1 \circledast T_2 \big] \Big\rangle(\omega_0, \omega'_{H^{T_1 \circledast T_2}}) 
			\\
			=& \Big\langle f \tilde{\otimes} g, \big( \circledast \tilde{\otimes} \circledast \big) \circ \overline{\mbox{Twist}} \circ \Delta \otimes \Delta \big[ T_1 \otimes T_2 \big] \Big\rangle(\omega_0, \omega'_{H^{T_1 \circledast T_2}}) 
			\\
			=& \Big\langle f \tilde{\otimes} g, \Delta\big[ T_1 \big] \Big\rangle(\omega_0, \omega'_{H^{T_1}}) \otimes \Big\langle f \tilde{\otimes} g, \Delta\big[ T_2 \big] \Big\rangle(\omega_0, \omega'_{H^{T_2}})
			\\
			=& \Big\langle f \ast g, T_1 \Big\rangle(\omega_0, \omega'_{H^{T_1}}) \otimes \Big\langle f \ast g, T_2 \Big\rangle(\omega_0, \omega'_{H^{T_2}}) \quad \bP \times (\bP')^{\times |H^{T_1 \circledast T_2}|}\mbox{-almost surely}
		\end{align*}
		which implies that $f \ast g \in \cG\big( \scH_\rp^{\gamma, \alpha, \beta}(\Omega, \Omega'), L^0(\Omega, \bP; \bR^e) \big)$. 
		
		Next, for any $T\in \scF_0^{\gamma, \alpha, \beta}$, we have that $\bP \times (\bP')^{\times |H^T|}$-almost surely
		\begin{align*}
			\Big\langle f\ast g, T \Big\rangle(\omega_0, \omega_{H^T}') = \sum_{\Upsilon, Y \in \scF_0} c\Big( T, \Upsilon, Y \Big) \cdot \Big\langle f, \Upsilon \Big\rangle \otimes^{H^T} \Big\langle g, Y \Big\rangle(\omega_0, \omega_{H^T}')
		\end{align*}
		and for each $\Upsilon, Y \in \scF^{\gamma, \alpha, \beta}$ we have that
		\begin{align*}
			\Big\langle f, \Upsilon\Big\rangle \in& L^0\bigg( \Omega, \bP; L^{q[\Upsilon]}\Big( (\Omega')^{\times |H^{\Upsilon}|}, (\bP')^{\times |H^{\Upsilon}|}; \lin\big( (\bR^d)^{\otimes |\scN^{\Upsilon}|}, \bR^e \big) \Big) \bigg) \quad \mbox{and}
			\\
			\Big\langle g, Y \Big\rangle \in& L^0\bigg( \Omega, \bP; L^{q[Y]}\Big( (\Omega')^{\times |H^{Y}|}, (\bP')^{\times |H^{Y}|}; \lin\big( (\bR^d)^{\otimes |\scN^{Y}|}, \bR^e \big) \Big) \bigg). 
		\end{align*}
		Applying the H\"older inequality with Equation \eqref{eq:definition:integrability-functional.2}, we conclude that
		\begin{align*}
			&\Big\langle f, \Upsilon \Big\rangle \otimes^{H^T} \Big\langle g, Y \Big\rangle \in L^0\bigg( \Omega, \bP; L^{q[T]}\Big( (\Omega')^{\times |H^{T}|}, (\bP')^{\times |H^{Y}|}; \lin\big( (\bR^d)^{\otimes |\scN^{T}|}, \bR^e \big) \Big) \bigg)
			\\
			&\mbox{so that}
			\\
			&\Big\langle f \ast g, T \Big\rangle \in L^0\bigg( \Omega, \bP; L^{q[T]}\Big( (\Omega')^{\times |H^{T}|}, (\bP')^{\times |H^{Y}|}; \lin\big( (\bR^d)^{\otimes |\scN^{T}|}, \bR^e \big) \Big) \bigg). 
		\end{align*}
		Thus we conclude that 
		\begin{equation*}
			f \ast g \in \cG\big( \scH_\rp^{\gamma, \alpha, \beta}(\Omega, \Omega'), L^0(\Omega, \bP; \bR^e) \big). 
		\end{equation*}
	\end{proof}
	
	\subsubsection{The set of characters is a group}
	
	Detailed results relating to the existence and derivation of the antipode for coupled Hopf algebras are not included in this work, for more information please see \cite{salkeld2022LionsTrees}. 
	
	\begin{definition}
		\label{definition:Antipode-McKean}
		Let $(\Omega', \cF', \bP')$ be a probability space and let $(\rp, \rQ)$ be a dual integrability functional. Let $\scS: \scH_\rp'(\Omega')^{\dagger} \to \scH_\rp'(\Omega')^{\dagger}$ be defined for $f\in \big( \scH_\rp'(\Omega') \big)^{\dagger}$ and $T \in \scT$ by
		\begin{align*}
			&\Big\langle \scS[f], \rId\Big\rangle = \Big\langle f, \rId \Big\rangle
			\\
			&\Big\langle \scS[f], T \Big\rangle(\omega'_{H^T}) 
			= - \Big\langle f, T\Big\rangle(\omega'_{H^T}) 
			- \sum_{c\in C(T)} \Big\langle \scS[f], T_c^P \Big\rangle \otimes^{H^T} \Big\langle f, T_c^R \Big\rangle(\omega'_{H^T}) 
		\end{align*}
		and extended to all Lions forests using the identity that $(\bP')^{\times |H^{T_1 \circledast T_2}|}$-almost surely, 
		\begin{equation*}
			\Big\langle \scS[f], T_1 \circledast T_2\Big\rangle(\omega'_{H^{T_1\circledast T_2}}) = \Big\langle \scS[f], T_1 \Big\rangle(\omega'_{H^{T_1}}) \otimes \Big\langle \scS[f], T_2\Big\rangle(\omega'_{H^{T_2}}). 
		\end{equation*}
		
		Let $(\Omega, \cF, \bP)$ be another probability space. We extend $\scS: \scH_\rp(\Omega, \Omega')^{\dagger} \to \big( \scH_\rp(\Omega, \Omega') \big)^{\dagger}$ to be defined inductively for $f\in \big( \scH_\rp(\Omega, \Omega') \big)^{\dagger}$ and $T\in \scT$ by
		\begin{equation}
			\label{eq:definition:Antipode-McKean-Bogoliubov-L}
			\left.
			\begin{aligned}
				&\Big\langle \scS[f], \rId \Big\rangle(\omega_0) = \Big\langle f, \rId \Big\rangle(\omega_0), 
				\\
				&\Big\langle \scS[f], T\Big\rangle(\omega_0, \omega'_{H^T}) 
				\\
				&\quad = 
				- \Big\langle f, T\Big\rangle(\omega_0, \omega_{H^T}) 
				- \sum_{c\in C(T)} \Big\langle \scS[f], T_c^P \Big\rangle \otimes^{H^T} \Big\langle f, T_c^R \Big\rangle(\omega_0, \omega'_{H^T}),  
				\\
				&\Big\langle \scS[f], T_1 \circledast T_2 \Big\rangle(\omega_0, \omega'_{H^{T_1 \circledast T_2}}) 
				= 
				\Big\langle \scS[f], T_1 \Big\rangle(\omega_0, \omega'_{H^{T_1}}) \otimes \Big\langle \scS[f], T_2 \Big\rangle(\omega_0, \omega'_{H^{T_2}}). 
			\end{aligned}
			\quad \right\}
		\end{equation}
	\end{definition}
	
	\begin{proposition}
		\label{proposition:Characters!}
		Let $(\Omega, \cF, \bP)$ and $(\Omega', \cF', \bP')$ be probability spaces. Let $\alpha, \beta>0$ and let $\gamma>\alpha\wedge\beta$. Let $(\rp, \rQ)$ be a dual integrability functional. Then
		\begin{align}
			\nonumber
			f \in \cG\big( \scH_\rp^{\gamma, \alpha, \beta}(\Omega, \Omega'), L^0(\Omega, \bP; \bR^e) \big) 
			\quad \implies \quad&
			\scS[f] \in \cG\big( \scH_\rp^{\gamma, \alpha, \beta}(\Omega, \Omega'), L^0(\Omega, \bP; \bR^e) \big)
			\\
			\label{eq:proposition:Characters!2}
			&\mbox{and} \quad
			f \ast \scS[f]  = \epsilon. 
		\end{align}
		That is, 
		\begin{equation*}
			\Big( \cG\big( \scH_\rp^{\gamma, \alpha, \beta}(\Omega, \Omega'), L^0(\Omega, \bP; \bR^e) \big), \ast \Big) \quad \mbox{ is a group with unit $\epsilon$. }
		\end{equation*}
	\end{proposition}

	\begin{proof}
		Let us start by supposing that $f\in \cG\big( \scH_\rp^{\gamma, \alpha, \beta}(\Omega, \Omega'), L^0(\Omega, \bP; \bR^e) \big)$. Firstly, 
		\begin{equation*}
			\Big\langle \scS[f], \rId \Big\rangle = \Big\langle f, \rId \Big\rangle \in L^0\Big( \Omega, \bP; \bR^d \Big). 
		\end{equation*}
		Now fix $n\in \bN$ and suppose that we have that for all $T\in \scF$ such that $|scN^T| \leq n$ we have that
		\begin{equation*}
			\Big\langle \scS[f], T\Big\rangle \in L^0\bigg( \Omega, \bP; L^{q[T]}\Big( (\Omega')^{\times |H^T|}, (\bP')^{\times |H^T|}; (\bR^d)^{\otimes |\scN^T|} \Big) \bigg). 
		\end{equation*}
		Consider $T\in \scF$ such that $\scN^T = n+1$. If $T$ is a Lions tree then thanks to Equation \eqref{eq:definition:Antipode-McKean-Bogoliubov-L} we have for any $c\in C(T)$ that
		\begin{align*}
			&\Big\langle f, T \Big\rangle \in L^0\bigg( \Omega, \bP; L^{q[T]}\Big( (\Omega')^{\times |H^T|}, (\bP')^{\times |H^T|}; (\bR^d)^{\otimes |\scN^T|} \Big) \bigg)
			\\
			&\Big\langle \scS[f], T_c^P \Big\rangle \otimes^{H^T} \Big\langle f, T_c^R \Big\rangle(\omega_0) \in L^0\bigg( \Omega, \bP; L^{q[T]}\Big( (\Omega')^{\times |H^T|}, (\bP')^{\times |H^T|}; (\bR^d)^{\otimes |\scN^T|} \Big) \bigg)
		\end{align*}
		so that 
		\begin{equation*}
			\Big\langle f, T \Big\rangle \in L^0\bigg( \Omega, \bP; L^{q[T]}\Big( (\Omega')^{\times |H^T|}, (\bP')^{\times |H^T|}; (\bR^d)^{\otimes |\scN^T|} \Big) \bigg)
		\end{equation*}
		On the other hand, if $T$ is a forest and not a tree then thanks to Equation \eqref{eq:definition:Antipode-McKean-Bogoliubov-L} and Equation \eqref{eq:definition:AlgebraRVs-int+} we also get that 
		\begin{align*}
			\Big\langle \scS[f], T \Big\rangle \in L^0\bigg( \Omega, \bP; L^{q[T]}\Big( (\Omega')^{\times |H^T|}, (\bP')^{\times |H^T|}; (\bR^d)^{\otimes |\scN^T|} \Big) \bigg)
		\end{align*}
		We conclude via induction that $\scS[f]$ satisfies the appropriate integrability. Further, we have that 
		\begin{equation*}
			\Big\langle \scS[f], T_1 \circledast T_2 \Big\rangle(\omega_0, \omega_{H^{T_1}}, \omega_{H^{T_2}} ) = \Big\langle \scS[f], T_1 \Big\rangle(\omega_0, \omega_{H^{T_1}} ) \otimes \Big\langle \scS[f], T_2 \Big\rangle(\omega_0, \omega_{H^{T_2}} )
		\end{equation*}
		$\bP \times (\bP')^{\times |H^{T_1 \circledast T_2}|}$-almost surely and for $T\in \scT$. As such, we have that
		\begin{equation*}
			\scS[f] \in \cG\big( \scH_\rp^{\gamma, \alpha, \beta}(\Omega, \Omega'), L^0(\Omega, \bP; \bR^e) \big). 
		\end{equation*}
		
		Next, recalling Equation \eqref{eq:definition:DualModule_ConvolProd} we get that for any $T\in \scT$ that
		\begin{align*}
			\Big\langle \scS[f] \ast f, T \Big\rangle(\omega_0, \omega'_{H^T}) =& \Big\langle f, T \Big\rangle(\omega_0, \omega'_{H^T}) + \Big\langle \scS[f], T \Big\rangle(\omega_0, \omega'_{H^T})
			\\
			&+ \sum_{c\in C(T)} \Big\langle \scS[f], T_c^P \Big\rangle \otimes^{H^T} \Big\langle f, T_c^R \Big\rangle(\omega_0, \omega'_{H^T})
			= 0
		\end{align*}
		$\bP \times (\bP')^{\times |H^T|}$-almost surely thanks to Equation \eqref{eq:definition:Antipode-McKean-Bogoliubov-L}. Further, 
		\begin{align*}
			\Big\langle \scS[f] \ast f, \rId \Big\rangle(\omega_0) = \Big\langle \scS[f], \rId\Big\rangle(\omega_0) \otimes \Big\langle f, \rId \Big\rangle(\omega_0) = 1 
		\end{align*}
		and applying the character property implies that for any $T\in \scF$ we get that
		\begin{align*}
			\Big\langle \scS[f] \ast f, T \Big\rangle(\omega_0, \omega'_{H^T}) = 0. 
		\end{align*}
		Hence $\scS[f] \ast f = \epsilon$ and we conclude with Proposition \ref{proposition:Characters1!} that 
		\begin{equation*}
			\Big( \cG\big( \scH_\rp^{\gamma, \alpha, \beta}(\Omega, \Omega'), L^0(\Omega, \bP; \bR^e) \big), \ast \Big) \quad \mbox{is a group. }
		\end{equation*}
	\end{proof}
	
	\subsubsection{McKean-Vlasov characters}

	In the context of McKean-Vlasov equations such as Equation \eqref{eq:meanfield:equation}, we turn our attention to an alternative subgroup of characters. Note that from this point on we use the convention that $(\Omega', \cF', \bP') = (\Omega, \cF, \bP)$ since our focus is now on McKean-Vlasov equations (where the probability space of the solution and the probability space of the lift of the solution law are the same) rather than other distribution dependent dynamics:
	\begin{definition}
		\label{definition:McKean-Vlasov-char}
		Let $(\Omega, \cF, \bP)$ be a probability space. Let $\alpha, \beta>0$, let $\gamma>\alpha\wedge \beta$ and let $(\rp, \rQ)$ be a dual integrability functional. We define
		\begin{equation*}
			G\Big( \scH_\rp^{\gamma, \alpha, \beta}(\Omega, \Omega), L^0(\Omega, \bP; \bR^e) \Big) \subseteq \cG\Big( \scH_\rp^{\gamma, \alpha, \beta}(\Omega, \Omega), L^0(\Omega, \bP; \bR^e) \Big)
		\end{equation*}
		to be the collection of characters that additionally satisfy that $\forall T \in \scF_0^{\gamma, \alpha, \beta}$ such that $h_0^T \neq \emptyset$, 
		\begin{equation}
			\label{eq:definition:DualModule2}
			\Big\langle f, T \Big\rangle(\omega_{h_0^T}, \omega_{H^T}) = \Big\langle f, \cE[T] \Big\rangle (\omega_0, \omega_{h_0^T}, \omega_{H^T}) \quad \bP \times \bP^{\times |H^{\cE[T]}|}\mbox{-almost surely. }
		\end{equation}
		When there is no ambiguity over the choice of probability space, we will denote 
		\begin{equation*}
			G_{\rQ}^{\gamma, \alpha, \beta}\big( L^0(\Omega, \bP; \bR^e)\big):= G\Big( \scH_\rp^{\gamma, \alpha, \beta}(\Omega, \Omega), L^0(\Omega, \bP; \bR^e)\Big).
		\end{equation*}
		We refer to $G_{\rQ}^{\gamma, \alpha, \beta}\big( L^0(\Omega, \bP; \bR^e)\big)$ as the set of \emph{McKean-Vlasov characters}. 
	\end{definition}
	To clarify a minor point here, Definition \ref{definition:McKean-Vlasov-char} only considers a single probability space. This is because we can identify the tagged probability space and the lifted probability space since they are the same. 
	
	The McKean-Vlasov group of characters describes the algebraic identities that the iterated integrals that occurred in Equation \eqref{eq:heuristic-Approx} need to satify. These are often referred to as  ``Chen's relation'' in the rough path literature.  
	
	\begin{remark}
		\label{remark:McKean-Vlasov-Characters}
		The intuition behind Equation \eqref{eq:definition:DualModule2} is actually quite natural. To justify this heuristically, consider the coefficient of a McKean-Vlasov equation (see Equation \eqref{eq:meanfield:equation}) and its canonical lift to a function on random variables (in the sense of Lions)
		\begin{equation*}
			f\Big( X_s(\omega), \cL_s^X \Big) = F\Big( X_s(\omega), \hat{X}_s(\omega, \hat{\omega}) \Big). 
		\end{equation*}
		The solution of the McKean-Vlasov equation $X_s$ is a random variable on the probability space $(\Omega, \cF, \bP)$ while the lift $\hat{X}_s$ is defined on the product probability space $(\Omega \times \hat{\Omega}, \cF \otimes \hat{\cF}, \bP \times \hat{\bP} )$ where $(\hat{\Omega}, \hat{\cF}, \hat{\bP})$ is the lifted probability space (which we remark is identical to the probability space $(\Omega, \cF, \bP)$ but we distinguish here for clarity) and $\bP$-almost surely the distribution 
		\begin{equation}
			\label{eq:remark:McKean-Vlasov-Characters}
			\hat{\bP} \circ \big( \hat{X}_s(\omega, \cdot) \big)^{-1}= \cL_s^X. 
		\end{equation}
		For any choice of $\omega$, the measure $\cL_s^X$ is the same so the random variable $\hat{X}$ with distribution $\cL_s^X$ is constant in $\omega$. Further, $X_s$ has distribution equal to $\cL_s^X$ so that we can choose $\hat{X}$ to satisfy that
		\begin{equation}
			\label{eq:remark:McKean-Vlasov-Characters-}
			\hat{X}_s(\omega, \hat{\omega}) = X_s(\hat{\omega}) \quad \bP \times \hat{\bP}\mbox{-almost surely. }
		\end{equation}
	
		as an application, consider some $a\in A[0]$ such that $a^{-1}[0] \neq \emptyset$. If we return to Equations \eqref{eq:heuristic-Approx} and \eqref{eq:heuristic-Approx-}, we can use Equation \eqref{eq:remark:McKean-Vlasov-Characters-} to rewrite the substitutions as 
		\begin{align*}
			\bigotimes_{i=1}^{|a|} \bigg( &X_{s, t}(\omega_0) \cdot \delta_{a_i = 0} + \hat{X}_{s, t}(\omega_0, \omega_{a_i}) \cdot \delta_{a_i > 0} \bigg) = \bigotimes_{i=1}^{|a|} X_{s, t}( \omega_{a_i} )
			\\
			&= \bigg( \bigotimes_{i=1}^{|a|} f\Big( X_s(\omega_{a_i}), \cL_s^X \Big) \bigg) \centerdot \bigg( \bigotimes_{i=1}^{|a|} W_{s, t}(\omega_{a_i}) \bigg)
		\end{align*}
		so that our probabilistic rough path satisfies that
		\begin{equation*}
			\Big\langle \rw_{s, t}, \cE^a\big[ \lfloor \rId \rfloor, ..., \lfloor \rId \rfloor \big] \Big\rangle(\omega_0, \omega_{m\{a\}} ) =  \bigotimes_{i=1}^{|a|} W_{s, t}(\omega_{a_i}). 
		\end{equation*}
		However, if we repeat this for some $\hat{a} \in A$ such that $\hat{a} = \llbracket a \rrbracket$ (see Definition \ref{def:a}) so that
		\begin{equation*}
			\cE\Big[ \cE^a\big[ \lfloor \rId \rfloor, ..., \lfloor \rId \rfloor \big] \Big] = \cE^{\hat{a}}\Big[ \lfloor \rId \rfloor, ..., \lfloor \rId \rfloor \Big], 
		\end{equation*}
		we get that the random variable
		\begin{equation*}
			\Big\langle \rw_{s, t}, \cE^{\hat{a}} \big[ \lfloor \rId \rfloor, ..., \lfloor \rId \rfloor \big] \Big\rangle(\omega_0, \omega_{m\{\hat{a}\}} ) =  \bigotimes_{i=1}^{|\hat{a}|} W_{s, t}(\omega_{\hat{a}_i}), 
		\end{equation*}
		which is constant in $\omega_0$. Further, by relabelling we get that
		\begin{align*}
			\Big\langle \rw_{s, t}&, \cE^a\Big[ \lfloor \rId \rfloor, ..., \lfloor \rId \rfloor \Big] \Big\rangle(\hat{\omega}_{0}, \hat{\omega}_{m\{a\}} ) 
			\\
			&= \Big\langle \rw_{s, t}, \cE^{\hat{a}}\Big[ \lfloor \rId \rfloor, ..., \lfloor \rId \rfloor \Big] \Big\rangle(\omega_0, \hat{\omega}_{\{0\}\cup m\{a\}} ) \quad \bP \times \bP \times \bP^{\times m[a]}\mbox{-almost surely. }
		\end{align*}
		This is just Equation \eqref{eq:definition:DualModule2}. 
		
		More generally, when we refer to a Lions forest with a non-empty $0$-hyperedge, we are saying that the appropriate object in question captures information about the tagged particle (which here refers to the solution $X_s(\omega)$). Equivalently, we are capturing pathwise information about the driving signal and the solution as opposed to information about their distribution. By contrast, a Lions forest with an empty $0$-hyperedge denotes that the appropriate object captures only information about the distribution (even if that distribution may be dependent on the tagged probability space). 
		
		For instance, a detagged hyperedge $h\in H$ identifies in the elementary differential that we derived in Equation \eqref{eq:heuristic-Approx-} a derivative in the measure variable so that the elementary differential is a function of the random variable $\hat{X}_s(\omega, \cdot)$. By contrast, the $0$-hyperedge identifies the derivatives in the elementary differential that correspond to spacial derivatives that carry $X_s(\omega)$ and thus indicate the presence of a dependency on the solution. 
		
		Hence, transforming a $0$-hyperedge into a detagged hyperedge equates with replacing spacial derivatives by Lions derivatives and replacing instances of $X_s(\omega_0)$ by $\hat{X}_s(\omega_0, \omega_h)$. When we translate this onto the characters (iterated integrals), what we see is this corresponds to replacing any instances of $\omega_0$ by $\omega_{h_0}$ which runs over a new probability space. Hence, the iterated integral corresponding to $\cE[T]$ is a random variable over a larger product probability space, but much like the lift $\hat{X}_s(\omega_0, \omega)$ it is constant in $\omega_0$. This inter-relationship is what Equation \eqref{eq:definition:DualModule2} describes. 
		
	\end{remark}

	\begin{proposition}
		\label{proposition:Characters!-MV}
		Let $(\Omega, \cF, \bP)$ be a probability space. Let $\alpha, \beta>0$, let $\gamma > \alpha \wedge \beta$ and let $(\rp, \rQ)$ be a dual integrability functional. Then the collection of McKean-Vlasov characters $G_{\rQ}^{\gamma, \alpha, \beta}\big( L^0(\Omega, \bP; \bR^e)\big)$ is a group with the convolution operation \eqref{eq:definition:DualModule_ConvolProd} and unit $\epsilon$. 
	\end{proposition}
	
	\begin{proof}
		Let $f, g \in G_{\rQ}^{\gamma, \alpha, \beta}\big( L^0(\Omega, \bP; \bR^e)\big)$. Courtesy of Proposition \ref{proposition:Characters!}, we have that
		\begin{equation*}
			f \ast g \in \cG^{\rQ}\big( \scH_\rp^{\gamma, \alpha,\beta} (\Omega, \Omega), L^0(\Omega, \bP; \bR^e) \big)
		\end{equation*}
		so that our focus is on verifying that $f\ast g$ satisfies Equation \eqref{eq:definition:DualModule2}. 
		
		For $T\in \scF^{\gamma, \alpha, \beta}$ such that $h_0^T \neq \emptyset$, 
		\begin{align*}
			\Big\langle f\ast g, T \Big\rangle(\omega_{h_0^T}, \omega_{H^T}) 
			&= \Big\langle f \tilde{\otimes} g, \Delta\big[ T \big] \Big\rangle(\omega_{h_0^T}, \omega_{H^T}) 
			\\
			&= \sum_{\Upsilon, Y \in \scF_0} c\Big( T, \Upsilon, Y\Big) \cdot \Big\langle f, \Upsilon \Big\rangle \otimes^{H^T} \Big\langle g, Y \Big\rangle (\omega_{h_0^T}, \omega_{H^T})
			\\
			&= \sum_{\Upsilon, Y \in \scF_0} c\Big( T, \Upsilon, Y\Big) \cdot \Big\langle f, \Upsilon \Big\rangle \otimes^{H^T\cup\{h_0^T\}} \Big\langle g, Y \Big\rangle (\omega_0, \omega_{h_0^T}, \omega_{H^T})
			\\
			&= \Big\langle f\ast g, \cE\big[ T \big] \Big\rangle(\omega_0, \omega_{h_0^T}, \omega_{H^T}) \quad \bP \times \bP^{\times |H^{\cE[T]}|}\mbox{-almost surely. }
		\end{align*}
		Hence $f\ast g\in G_{\rQ}^{\gamma, \alpha, \beta}\big( L^0(\Omega, \bP; \bR^e) \big)$. 
		
		Finally, for $f\in G_{\rQ}^{\gamma, \alpha, \beta}\big( L^0(\Omega, \bP; \bR^e) \big)$, we verify that the inverse $\scS[f] \in G_{\rQ}^{\gamma, \alpha, \beta}\big( L^0(\Omega, \bP; \bR^e) \big)$. Firstly, 
		\begin{equation*}
			\Big\langle \scS[f], \lfloor \rId \rfloor \Big\rangle(\omega_{h_0}) = \Big\langle \scS[f] , \cE\big[ \lfloor \rId \rfloor \big] \Big\rangle(\omega_0, \omega_{h_0}) \quad \bP \times \bP\mbox{-almost surely. } 
		\end{equation*}
		
		Next, let $n\in \bN$ and suppose that for all $T' \in \scF$ such that $|\scN^{T'}|<n$ and $h_0^T \neq \emptyset$, we have that 
		\begin{equation*}
			\Big\langle \scS[f], T'\Big\rangle(\omega_{h_0}, \omega_{H^T}) = \Big\langle \scS[f], \cE[T'] \Big\rangle(\omega_0, \omega_{h_0}, \omega_{H^T}) \quad \bP \times \bP^{\times |H^{\cE[T']}|}\mbox{-almost surely. }
		\end{equation*}
		
		Let $T\in \scT$ such that $h_0^T \neq \emptyset$ and $|\scN^T| = n$. Then
		\begin{align*}
			&\Big\langle \scS[f], T\Big\rangle(\omega_{h_0^T}, \omega_{H^T}) = - \Big\langle f, T\Big\rangle(\omega_{h_0^T}, \omega_{H^T}) - \sum_{c\in C(T)} \Big\langle \scS[f], T_c^P \Big\rangle \otimes^{H^T} \Big\langle f, T_c^R \Big\rangle(\omega_{h_0^T}, \omega_{H^T})
			\\
			&= - \Big\langle f, \cE[T] \Big\rangle(\omega_0, \omega_{h_0^T}, \omega_{H^T}) - \sum_{c\in C(T)} \Big\langle \scS[f], \cE[T_c^P] \Big\rangle \otimes^{H^T\cup \{h_0^T\}} \Big\langle f, \cE[T_c^R] \Big\rangle(\omega_0, \omega_{h_0^T}, \omega_{H^T}) 
			\\
			&= \Big\langle \scS[f], \cE[T] \Big\rangle(\omega_0, \omega_{h_0^T}, \omega_{H^T}) \quad \bP \times \bP^{\times |H^{\cE[T]}|}\mbox{-almost surely. }
		\end{align*}
		where the second line is thanks to the inductive hypothesis. 
	\end{proof}

	\begin{example}
		\label{example:TanakaTrick}
		Let $N \in \bN$ and let $(\Omega, \cF, \bP)$ be a probability space carrying $N$ independent copies of a driving signal $W$. 
		
		Consider the coefficient of the interacting equation (see Equation \eqref{eq:particle:system}) and its canonical lift to a function on random variables
		\begin{equation*}
			f\Big( X_s^{i, N}(\omega_0), \tfrac{1}{N} \sum_{j=1}^N \delta_{X_s^{j, N}(\omega_0)} \Big) = F \Big( X_s^{i, N}(\omega_0), \hat{X}_s(\omega_0, \omega') \Big)
		\end{equation*}
		where the random variable $\hat{X}_s(\omega_0, \omega') = X_s^{u(\omega'), N}(\omega_0)$ and $u$ is uniformly distributed on the set $\{1, ..., N\}$. Hence, $\bP$-almost surely, the random variable $\hat{X}_s(\omega_0, \cdot)$ is uniformly distributed on the set
		\begin{equation*}
			\Big\{ X_s^{1, N}(\omega_0), ..., X_s^{N, N}(\omega_0) \Big\}. 
		\end{equation*}
		Following the ideas of \cite{tanaka1984limit} (and the more recent work \cite{deuschel2017enhanced} which has a rough path perspective), we can consider the (random) probability space
		\begin{equation*}
			\Omega'_{\omega_0} = \Big\{ W^{1, N}(\omega_0), ..., W^{N, N}(\omega_0) \Big\} \subseteq C\big( [0,T]; \bR^d \big)
			\quad
			\cF' = 2^{\Omega'_{\omega_0}}
			\quad \mbox{and} \quad
			\bP'_{\omega_0} = \frac{1}{N} \sum_{j=1}^N \delta_{W^{j, N}(\omega_0)}. 
		\end{equation*}
		On the probability space $(\Omega'_{\omega_0}, \cF'_{\omega_0}, \bP'_{\omega_0})$, the dynamics of the McKean-Vlasov described in Equation \eqref{eq:meanfield:equation} are $\bP$-almost surely the same as the dynamics of \eqref{eq:particle:system}. This is sometimes referred to as the ``\emph{Tanaka trick}''. 
		
		Denoting by $u$ a uniformly distributed random variable on the set $\{1, ..., N\}$, let $X_s^{u(\omega'), N}(\omega_0)$ be the value of a uniformly chosen element of the interacting particle system as time $s$ and let $\hat{X}_s^{u(\omega'), N}(\omega_0, \hat{\omega}')$ be a random variable with distribution
		\begin{equation*}
			\bP_{\omega_0} \circ \Big( \hat{X}_s^{u(\omega'), N}(\omega_0, \cdot) \Big) = \frac{1}{N} \sum_{j=1}^N \delta_{X_s^{j, N}(\omega_0)}. 
		\end{equation*}
		Then we can choose the random variable $\hat{X}$ to satisfy
		\begin{equation*}
			\hat{X}_s^{u(\omega'), N}(\omega_0, \hat{\omega}') = X_s^{u(\hat{\omega}'), N}(\omega_0) \quad \bP'_{\omega_0} \times \hat{\bP}'_{\omega_0} \mbox{-almost surely. }
		\end{equation*}
		This matches Equation \eqref{eq:remark:McKean-Vlasov-Characters-} so that we are able to derive the same inter-relationship between the terms of the characters that we did in Remark \ref{remark:McKean-Vlasov-Characters} and we obtain Equation \eqref{eq:definition:DualModule2}. 
	\end{example}
	
	\newpage
	\section{Probabilistic and analytic structures}
	\label{section:Models}
	
	The theory of rough paths, first proposed in \cite{lyons1998differential}, is now a wide ranging, multi-disciplined field of research. Over the last twenty years, the field has developed and there are now many different approaches to defining what a rough path is with differing levels of abstraction. The concept of a branched rough paths was first introduced in \cite{gubinelli2010ramification}. However, here we use a definition closer to that of the recent work \cite{tapia2020geometry}. 
	
	Let $\alpha>0$ and let $\gamma = \alpha \lfloor \tfrac{1}{\alpha}\rfloor$. Let $\big( \cH, \odot, \rId, \triangle, \epsilon, \cS \big)$ be an $\bN$-graded Hopf algebra over a normed ring $\big( \cR, +, \centerdot, \| \cdot \| \big)$ with basis $\fF$. 
	
	We say that $\rw:[0,1] \to G(\cH, \cR)$ is a $(\cH, \alpha)$-rough path if it satisfies that $\forall s,t,u\in [0,1]$ such that $s\leq t\leq u$, 
	\begin{equation}
		\label{eq:BranchedRP1}
		\rw_{s, t} = (\rw_s)^{-1} \ast \rw_t, 
		\qquad
		\rw_{s,t} \ast \rw_{t,u} = \rw_{s,u}
	\end{equation}
	and $\forall \tau\in \fF$
	\begin{equation}
		\label{eq:BranchedRP2}
		\big\| \langle \rw_{s, t}, \tau \rangle \big\| \lesssim |t-s|^{\alpha \cdot |\tau |} . 
	\end{equation}
	The set of all $(\cH, \alpha)$-rough paths is denoted $\cC(\cH, \alpha)$
	
	In \cite{hairer2014theory}, this concept was generalised to solve singular stochastic partial differential equations. Regularity structures use an abstract Taylor expansion that best approximates the solution to determine the relevant necessary information about the driving noise to solve an equation. 
	
	In Subsection \ref{section:TaylorExpansions} we highlighted some of the new results (found in \cite{salkeld2022Lions}) relating to Lions calculus. In Section \ref{section:ProbabilisticRoughPaths}, we highlighted some of the properties of partition/Lions trees, identified  links to Lions calculus and explored some of the the arising algebraic properties relating to coupled Hopf algebras (with more details found in \cite{salkeld2022LionsTrees}). 
	
	In this section, we study the central application of these results to the understanding of probabilistic rough paths that drive McKean-Vlasov equations and similarly systems of interacting equations. We start with the definition of probabilistic rough paths in Subsection \ref{subsection:PRPs-defin}, drawing inspiration from \cite{gubinelli2010ramification} which describes a rough path as a path on the characters of a Hopf algebra (see also Subsection \ref{subse:2.6}) that satisfies a certain integrability condition . 
	
	Subsection \ref{subsection:Examples} focuses on important examples and Subsection \ref{subsection:StrongPRP} proves the existence of strong probabilistic rough paths, a critical concept first introduced in \cite{2019arXiv180205882.2B} that demonstrates the measurability of the lift from a driving signal to a probabilistic rough path. 
	
	\subsection{Probabilistic rough paths}
	\label{subsection:PRPs-defin}
	
	We start by introducing a general definition for a probabilistic rough path:
	\begin{definition}
		\label{definition:General-PRP}
		Let $(\Omega, \cF, \bP)$ and $(\Omega', \cF', \bP')$ be probability spaces, and let $\big( L^0(\Omega, \bP; \bR^e), +, \centerdot\big)$ be a ring of measurable functions. Let $(\rp, \rQ)$ be a dual integrability functional. Let $\alpha, \beta>0$ and let 
		\begin{equation}
			\label{definition:General-PRP-gamma}
			\gamma:= \inf\{ \scG_{\alpha, \beta}[T]: T \in \scF, \scG_{\alpha, \beta}[T]> 1-\alpha \}. 
		\end{equation}
		
		Let $\big( \scH_\rp^{\gamma, \alpha, \beta}(\Omega, \Omega'), \circledast, \rId, \Delta, \epsilon \big)$ be a coupled bialgebra over the ring $L^0(\Omega, \bP; \bR^e)$ with index set $\scF$ and grading $\scG_{\alpha, \beta}$.  
		
		Let $\big( \cG(\scH_\rp^{\gamma, \alpha, \beta}(\Omega, \Omega'), L^0(\Omega, \bP; \bR^e)), \ast \big)$ be the group of characters. We say that 
		\begin{equation*}
			\rw:[0,1] \to \cG\big( \scH_\rp^{\gamma, \alpha, \beta}(\Omega, \Omega'), L^0(\Omega, \bP; \bR^e) \big)
		\end{equation*}
		is a \emph{probabilistic rough path} if
		\begin{itemize}
			\item $\rw$ satisfies a coupled Chen's relationship: $\forall s, t, u\in[0,1]$, 
			\begin{align*}
				&\rw_{s, t}(\omega_0) = (\rw_s)^{-1}(\omega_0) \ast \rw_t(\omega_0)
				\quad\mbox{and}\quad
				\\
				&\rw_{s, t}(\omega_0) \ast \rw_{t, u}(\omega_0) = \rw_{s, u}(\omega_0) \quad \bP\mbox{-almost surely. } 
			\end{align*}
			\item For all $T\in \scF$, 
			\begin{equation*}
				\sup_{s, t\in[0,1]} \frac{(\bE')^{H^T}\bigg[ \Big| \big\langle \rw_{s, t}, T \big\rangle(\omega_0, \omega'_{H^T}) \Big| \bigg] }{\big| t-s \big|^{\scG_{\alpha, \beta}[T]}} < \infty
				\quad \mbox{$\bP$-almost surely. }
			\end{equation*}
		\end{itemize}
		
		We denote the set of probabilistic rough paths by
		\begin{equation*}
			\scC\Big( [0,1]; \cG\big( \scH_\rp^{\gamma, \alpha, \beta}(\Omega, \Omega'), L^0(\Omega, \bP; \bR^e) \big) \Big). 
		\end{equation*}
	\end{definition}
	In practice, one will want to specify a choice of coupled Hopf algebra tailored to the distributional dynamics under consideration and as we saw in Subsection \ref{subse:2.6} the group of characters is often much larger than the group on which the probabilistic rough path will run over. Therefore, we should only view Definition \ref{definition:General-PRP} as an umbrella under which relevant examples of probabilistic rough paths are included. In the rest of this Section, we will provide some more relevant examples of probabilistic rough paths. 
		
	\subsubsection*{Probabilistic rough paths for McKean-Vlasov equations}
	
	We consider a probabilistic rough path for McKean-Vlasov equations such as Equation \eqref{eq:meanfield:equation} (or alternatively using the Tanaka trick from Example \ref{example:TanakaTrick} to address Equation \eqref{eq:particle:system}):
	\begin{definition}
		\label{definition:ProbabilisticRoughPaths}
		Let $\alpha, \beta>0$ and let $\gamma$ satisfy Equation \eqref{definition:General-PRP-gamma}. Let $(\Omega, \cF, \bP)$ be a probability space. Let $(\rp, \rQ)$ be a dual integrability functional. 
		
		We say that $\rw:[0,1] \to G_{\rQ}^{\gamma, \alpha, \beta}\big( L^{0}(\Omega, \bP; \bR^e)\big)$ is a \emph{McKean-Vlasov probabilistic rough paths} if 
		\begin{enumerate}
			\item $\rw$ satisfies a coupled Chen's relationship: $\forall s, t, u \in [0,1]$, 
			\begin{equation}
				\label{eq:definition:ProbabilisticRoughPaths1}
				\begin{split}
					&\rw_{s, t}(\omega_0) = (\rw_s)^{-1}(\omega_0) \ast \rw_t(\omega_0)
					\quad\mbox{and}
					\\
					&\rw_{s,t}(\omega_0) \ast \rw_{t,u}(\omega_0) = \rw_{s,u}(\omega_0) \quad \mbox{$\bP$-almost surely.} 
				\end{split}
			\end{equation}
			\item
			For all $T \in \scF^{\gamma, \alpha, \beta}$, 
			\begin{align}
				\label{eq:definition:ProbabilisticRoughPaths2}
				\sup_{s, t\in[0,1]} \frac{\bE^{H^T}\bigg[ \Big| \big\langle \rw_{s, t} , T\big\rangle(\omega_{0}, \omega_{H^T}) \Big| \bigg]}{\big| t-s \big|^{\scG_{\alpha, \beta}[T]}}  < \infty \qquad \mbox{$\bP$-almost surely. }
			\end{align}
		\end{enumerate}
	\end{definition}

	\begin{remark}
		Thanks to Equation \eqref{eq:definition:DualModule2}, we have that for $T\in \scF$ and for any $s, t\in [0,1]$ that
		\begin{equation*}
			\frac{\bE^{\{h_0^T\} \cup H^{T}}\bigg[ \Big| \big\langle \rw_{s, t}, \cE[T] \big\rangle(\omega_0, \omega_{h_0^T}, \omega_{H^T}) \Big| \bigg] }{|t-s|^{\beta \cdot |\scN^T|}}
			=
			\bE^0\Bigg[ \frac{\bE^{H^{T}}\bigg[ \Big| \big\langle \rw_{s, t}, T \big\rangle( \omega_{0}, \omega_{H^T}) \Big| \bigg] }{|t-s|^{\beta \cdot |\scN^T|}} \Bigg]
		\end{equation*}
		Therefore, for any $T\in \scT_0$ we have that the random variable $\langle \rw_{s, t}, T\rangle$ satisfies both Equation \eqref{eq:definition:ProbabilisticRoughPaths2} and
		\begin{equation*}
			\sup_{s, t\in [0,1]} \bE^0\Bigg[ \frac{\bE^{H^T}\bigg[ \Big| \big\langle \rw_{s, t} , T\big\rangle(\omega_{0}, \omega_{H^T}) \Big| \bigg]}{\big| t-s \big|^{\beta\cdot|\scN^T|}} \Bigg]< \infty
		\end{equation*}
	\end{remark}

	\subsection{Examples of probabilistic rough paths}
	\label{subsection:Examples}
	
	Here, we introduce some examples of different types of probabilistic rough paths and give some context for why these are useful for solving  mean-field equations with applications. 
	
	\subsubsection{Probabilistic rough paths for systems of interacting Equations}
	
	Let us start by recalling Example \ref{example:regularitydiff}: as there was no change in regularity by comparing Equation \eqref{eq:example:regularitydiff0} and \eqref{eq:example:regularitydiff1}, for this section we will be working with the running assumption that $\alpha = \beta$. For $T\in \scF$ and recalling \eqref{eq:lemma:grading}, we could equivalently write $\scG_{\alpha, \alpha}[T] = \alpha |T| $ since 
	\begin{equation*}
		\scG_{\alpha, \alpha}[T] = \alpha |h_0^T| + \alpha |\scN^T \backslash h_0^T|. 
	\end{equation*}
	
	In this setting, because all of the random variables over the detagged probability spaces have finite support, we should not expect there to be any increase in regularity when we consider the mean-square H\"older continuity. 
	\begin{example}
		\label{example:smoothpaths-trees-2}
		Recall from Example \ref{example:smoothpaths-trees} that $C^1([0,1]; \bR^d)$ is the space of continuously differentiable paths and let $\cB$ be the associated Borel $\sigma$-field. Let $q\geq 2$ and let 
		\begin{equation*}
			\mu \in \cP\Big( C^1([0,1]; \bR^d) \Big) 
			\quad \mbox{such that} \quad 
			\int_{C^1([0,1]; \bR^d)} \| x \|_1^q d\mu(x) < \infty. 
		\end{equation*}
		Let $N\in \bN$, let $(\Omega, \cF, \bP)$ be a probability space and for $i=1, ..., N$ let $W^{i, N}:\Omega \to C^1([0,1]; \bR^d)$ be a sequence of independent, identically distributed random variables such that $\bP \circ (W^{i,N})^{-1} = \mu$. 
		
		Let $i\in \{1, ..., N\}$ and let $\nu$ be the uniform measure on $\big(\{1, .., N\}, 2^{\{1, .., N\}} \big)$. For each $T\in \scF$ and $s, t \in[0,1]$ such that $s\leq t$, we define the operators
		\begin{equation*}
			\cI_{s, t}^{i, T}: \Omega \times \{1, ..., N\}^{\times |H^T|} \to (\bR^d)^{\otimes |\scN^T|}
		\end{equation*}
		inductively using the relationships
		\begin{equation}
			\label{eq:example:smoothpaths-trees-2}
			\begin{split}
				&\cI_{s, t}^{i, \lfloor \rId \rfloor}(\omega) = W_{s, t}^{i, N}(\omega), 
				\\
				\forall T_1, T_2\in \scF, \quad 
				&\cI_{s, t}^{i, T_1 \circledast T_2}(\omega)(\iota_{H^{T_1}}, \iota_{H^{T_2}}) = \cI_{s, t}^{i, T_1}(\omega)(\iota_{H^{T_1}}) \otimes \cI_{s, t}^{i, T_2}(\omega)(\iota_{H^{T_2}}), 
				\\
				\forall T\in \scF \quad \mbox{s.t.} \quad h_0^T \neq \emptyset, \quad 
				&\cI_{s, t}^{i, \cE[T]}(\omega)(\iota_{h_0^T}, \iota_{H^{T}}) = \cI_{s, t}^{\iota_{h_0}, T}(\omega)(\iota_{H^{T}}),
				\\
				&\cI_{s, t}^{i, \lfloor T \rfloor}(\omega)(\iota_{H^{T}}) = \int_s^t \cI_{s, r}^{i, T}(\omega)(\iota_{H^{T}}) \otimes dW_r^{i, N}
			\end{split}
		\end{equation}
		with the usual convention that $\iota_{H} = \underbrace{(\iota_h, ... )}_{h \in H}$. 
		
		Then $\forall T \in \scF_0$ and $s, t\in [0,1]$ the mapping
		\begin{equation*}
			(i, \iota_{H^T} ) \mapsto \cI_{s, t}^{i, T}(\omega)(\iota_{H^T}) 
			\quad \mbox{is} \quad
			\Big( \big\{ 1, ..., N \big\}^{\times (1+|H^T|)}, 2^{\{1, ..., N\}^{\times (1+|H^T|)}} \Big)\mbox{-measurable}
			\quad  \bP\mbox{-almost surely}
		\end{equation*}
		so that for $\textbf{0} = (0)_{H^T}$, 
		\begin{equation*}
			\cI_{s, t}^{\cdot, T}\in L^0\bigg( \{1, ..., N\}, \nu; L^{\textbf{0}} \Big( \{1, ..., N\}^{\times |H^T|}, \nu^{\times |H^T|}; (\bR^d)^{\otimes |\scN^T|} \Big) \bigg). 
		\end{equation*}
		Let $\alpha>0$ and let $\gamma>\alpha$ and suppose that $\forall T \in \scF^{\gamma, \alpha, \alpha}$ that $\tfrac{q}{|\scN^T|}>1$. Then we define
		\begin{equation}
			\label{eq:trivialIntFunct}
			\rQ = (q[T])_{T\in \scF^{\gamma, \alpha, \beta}}, 
			\quad 
			q[T] = \big( q[T]_h \big)_{h\in H^T},
			\quad q[T]_h = \frac{q}{|h|}. 
		\end{equation}
		Thus, for each $T\in \scF^{\gamma, \alpha, \alpha}$ we have that $\bP$-almost surely
		\begin{equation*}
			\cI_{s, t}^{T}\in L^0\bigg( \{1, ..., N\}, \nu; L^{q[T]}\Big( \{1, .., N\}^{\times |H^T|}, \nu^{\times |H^T|}; (\bR^d)^{\otimes |\scN^T|} \Big) \bigg)
		\end{equation*}
		and for each $\omega \in \Omega$ we define
		\begin{equation*}
			\cI_{s, t}^{i, \gamma, \alpha}[\mu](\omega) = \sum_{T\in \scF_0}^{\gamma,\alpha, \alpha} \cI_{s, t}^{i, T}(\omega)(\cdot) \in \bigoplus_{T\in \scF_0}^{\gamma, \alpha, \alpha} L^{q[T]}\Big( \{1, ..., N\}^{\times |H^T|}, \nu^{\times |H^T|}; (\bR^d)^{\otimes |\scN^T|} \Big)
		\end{equation*}
		so that
		\begin{equation*}
			\cI_{s, t}^{\gamma, \alpha}[\mu](\omega) \in L^0\bigg( \{1, ..., N\}, \nu; \bigoplus_{T\in \scF_0}^{\gamma, \alpha, \alpha} L^{q[T]}\Big( \{1, ..., N\}^{\times |H^T|}, \nu^{\times |H^T|}; (\bR^d)^{\otimes |\scN^T|} \Big) \bigg). 
		\end{equation*}
	\end{example}
	
	\begin{proposition}
		\label{proposition:I_geoPRP2}
		Let $\alpha>0$, let $\gamma:= \sup\{ i\alpha: i\in \bN, i\alpha\leq 1 \}$ and let $q' > \big\lfloor \tfrac{\gamma}{\alpha}\big\rfloor \geq 1$. Let $(\rp, \rQ)$ be the dual integrability function that satisfies Equation \eqref{eq:trivialIntFunct}. 
		
		Let $s, t, u \in[0,1]$. Suppose that $\mu \in \cP_{q'}\big( C^1([0,1]; \bR^d)\big)$ and let $(\Omega, \cF, \bP)$ be a probability space carrying a collection of $N$ random variables $\big(W^{i, N} \big)_{i=1, ..., N}$ such that $\bP \circ (W^{i, N})^{-1}$. 
		
		Let $\nu$ be the uniform measure on $\{1, ..., N\}$ and define the (random) probability space
		\begin{equation*}
			\Omega'_{\omega_0} = \Big\{ W^{1, N}(\omega_0), ..., W^{N, N}(\omega_0) \Big\}, 
			\quad
			\cF_{\omega_0}' = 2^{\Omega'_{\omega_0}}
			\quad
			\bP'_{\omega_0} = \frac{1}{N} \sum_{j=1}^N \delta_{W^{j, N}(\omega_0)}. 
		\end{equation*}

		Then $\bP$-almost surely the path $\cI_{0, \cdot}^{i, \gamma, \alpha}[\mu](\omega)$ as defined in Equation \eqref{eq:example:smoothpaths-trees-2} are elements of the group 
		\begin{equation*}
			\overline{G}_{\rQ}^{\gamma, \alpha, \alpha} \Big( L^0\big(\{1, ..., N\}, \nu; \bR^e \big) \Big)
		\end{equation*}
		and satisfies the identity
		\begin{equation*}
			\cI_{s, t}^{i, \gamma, \alpha}[\mu](\omega) \ast \cI_{t, u}^{i, \gamma, \alpha}[\mu](\omega) = \cI_{s, u}^{i, \gamma, \alpha}[\mu](\omega)
			\quad
			\bP\mbox{-almost surely.}
		\end{equation*}
		As such, the path $\cI^{i, \gamma, \alpha}$ introduced in Example \ref{example:smoothpaths-trees-2} is a probabilistic rough path in the sense of Definition \ref{definition:ProbabilisticRoughPaths}. 
	\end{proposition}
	
	\begin{proof}
		Equation \eqref{eq:example:smoothpaths-trees-2} makes it a simple exercise to verify that for any $t\in [0,1]$, 
		\begin{equation*}
			\cI_{0, t}^{\gamma, \alpha}[\mu](\omega_0) \in \overline{G}_{\rQ}^{\gamma, \alpha, \alpha} \Big( L^0\big( \{1, ..., N\}, \nu; \bR^e \big) \Big)
		\end{equation*}
		and we prove via induction that $\forall T \in \scF_0^{\gamma, \alpha, \alpha}$ that $\bP$-almost surely
		\begin{equation}
			\label{eq:proposition:I_geoPRP-2}
			\Big\langle \cI_{s, t}^{i, \gamma, \alpha}[\mu](\omega) \ast \cI_{t, u}^{i, \gamma, \alpha}[\mu](\omega) , T \Big\rangle = \Big\langle \cI_{s, u}^{i, \gamma, \alpha}[\mu](\omega), T \Big\rangle. 
		\end{equation}
		Firstly, 
		\begin{align*}
			&\Big\langle \cI_{s, t}^{i, \gamma, \alpha}[\mu](\omega) \ast \cI_{t, u}^{i, \gamma, \alpha}[\mu](\omega), \lfloor \rId\rfloor \Big\rangle 
			\\
			&= 
			\Big\langle \cI_{s, t}^{i, \gamma, \alpha}[\mu](\omega), \lfloor \rId\rfloor \Big\rangle \otimes \Big\langle \cI_{t, u}^{i, \gamma, \alpha}[\mu](\omega), \rId \Big\rangle + \Big\langle \cI_{s, t}^{i, \gamma, \alpha}[\mu](\omega), \rId \Big\rangle \otimes \Big\langle \cI_{t, u}^{i, \gamma, \alpha}[\mu](\omega), \lfloor \rId\rfloor \Big\rangle
			\\
			&=\Big\langle \cI_{s, t}^{i, \gamma, \alpha}[\mu](\omega), \lfloor \rId\rfloor \Big\rangle + \Big\langle \cI_{t, u}^{i, \gamma, \alpha}[\mu](\omega), \lfloor \rId\rfloor \Big\rangle
			\\
			&=W_{s, t}^i(\omega) + W_{t, u}^i(\omega) 
			= W_{s, u}^i(\omega) = \Big\langle \cI_{s, u}^{i, \gamma, \alpha}[\mu](\omega), \lfloor \rId\rfloor \Big\rangle
		\end{align*}
		and
		\begin{align*}
			&\Big\langle \cI_{s, t}^{i, \gamma, \alpha}[\mu](\omega) \ast \cI_{t, u}^{i, \gamma, \alpha}[\mu](\omega), \cE\big[\lfloor \rId\rfloor\big] \Big\rangle (\iota)
			\\
			&= 
			\Big\langle \cI_{s, t}^{i, \gamma, \alpha}[\mu](\omega), \cE\big[\lfloor \rId\rfloor \big] \Big\rangle(\iota) \otimes \Big\langle \cI_{t, u}^{i, \gamma, \alpha}[\mu](\omega), \rId \Big\rangle + \Big\langle \cI_{s, t}^{i, \gamma, \alpha}[\mu](\omega), \rId \Big\rangle \otimes \Big\langle \cI_{t, u}^{i, \gamma, \alpha}[\mu](\omega), \cE\big[ \lfloor \rId\rfloor \big] \Big\rangle(\iota)
			\\
			&=\Big\langle \cI_{s, t}^{i, \gamma, \alpha}[\mu](\omega), \cE\big[ \lfloor \rId\rfloor \big] \Big\rangle(\iota) + \Big\langle \cI_{t, u}^{i, \gamma, \alpha}[\mu](\omega), \cE\big[ \lfloor \rId\rfloor \big] \Big\rangle(\iota)
			\\
			&=W_{s, t}^{\iota}(\omega) + W_{t, u}^{\iota}(\omega) 
			= W_{s, u}^{\iota}(\omega) = \Big\langle \cI_{s, u}^{i, \gamma, \alpha}[\mu](\omega), \cE\big[ \lfloor \rId\rfloor \big] \Big\rangle(\iota)
		\end{align*}
		so Equation \eqref{eq:proposition:I_geoPRP-2} holds for $\lfloor \rId \rfloor$ and $\cE\big[\lfloor \rId\rfloor\big]$. 
		
		Next, suppose that Equation \eqref{eq:proposition:I_geoPRP-2} holds for $T_1, T_2 \in \scF^{\gamma, \alpha, \alpha}$ and $T_1\circledast T_2 \in \scF^{\gamma, \alpha, \alpha}$. Using Sweedler notation, thanks to Equation \eqref{eq:proposition:H-coupledBialgebra} we get that 
		\begin{align*}
			&\Delta\big[ T_1 \big] = \sum_{j=0}^m T_1^{j, 1} \times^{H^{T_1}} T_1^{j, 2}, 
			\quad
			\Delta\big[ T_2 \big] = \sum_{k=0}^n T_2^{k, 1} \times^{H^{T_2}} T_k^{j, 2} \quad \mbox{and}
			\\
			&\Delta\Big[ T_1 \circledast T_2 \Big] 
			= 
			\sum_{j=0}^m \sum_{k=0}^n \big[ T_1^{j, 1} \circledast T_2^{k, 1} \big] \times^{H^{T_1} \cup H^{T_2}} \big[ T_1^{j, 2} \circledast T_2^{k, 2} \big]
		\end{align*}
		
		By denoting 
		\begin{align*}
			H_1^{j,k} = H^{T_1^{j, 1}} \cup H^{T_2^{k, 1}},
			\quad 
			H_2^{j, k} = H^{T_1^{j, 2}} \cup H^{T_2^{k, 2}}
			\quad \mbox{and}\quad
			H = H^{T_1}\cup H^{T_2}, 
		\end{align*}
		we get that
		\begin{align*}
			&\Big\langle \cI_{s, t}^{i, \gamma, \alpha}[\mu](\omega) \ast \cI_{t, u}^{i, \gamma, \alpha}[\mu](\omega), T_1 \circledast T_2 \Big\rangle( \iota_{H} )
			\\
			&=\sum_{i=0}^{m} \sum_{j=0}^{n} \Big\langle \cI_{s, t}^{i, \gamma, \alpha}[\mu](\omega), T_1^{j, 1} \circledast T_2^{k, 1} \Big\rangle(\iota_{\psi^{H_1^{j, k}, H}[H_1^{j, k}]}) 
			\\
			&\qquad \otimes \Big\langle \cI_{t, u}^{i, \gamma, \alpha}[\mu](\omega), T_1^{j,2} \circledast T_2^{k, 2} \Big\rangle(\iota_{\psi^{H_2^{j, k}, H}[H_2^{j, k}]})
			\\
			&=\sum_{j=0}^{m} \sum_{k=0}^{n} \bigg( \Big\langle \cI_{s, t}^{i, \gamma, \alpha}[\mu](\omega), T_1^{j, 1} \Big\rangle(\iota_{\psi^{T_1^{j, 1}, H}[H^{T_1^{j, 1}}]}) \otimes \Big\langle \cI_{t, u}^{i, \gamma, \alpha}[\mu](\omega), T_2^{k, 1} \Big\rangle(\iota_{\psi^{T_2^{k, 1}, H}[H^{T_2^{k, 1}}]}) \bigg)
			\\
			&\qquad \otimes
			\bigg( \Big\langle \cI_{s, t}^{i, \gamma, \alpha}[\mu](\omega), T_1^{j, 2} \Big\rangle(\iota_{\psi^{T_1^{j, 2}, H}[H^{T_1^{j,2}}]}) \otimes \Big\langle \cI_{t, u}^{i, \gamma, \alpha}[\mu](\omega), T_2^{k, 2} \Big\rangle(\iota_{\psi^{T_2^{k, 2}, H}[H^{T_2^{k,2}}]}) \bigg)
			\\
			&=\bigg( \sum_{j=0}^{m} \Big\langle \cI_{s, t}^{i, \gamma, \alpha}[\mu](\omega), T_1^{j, 1} \Big\rangle(\iota_{\psi^{T_1^{j, 1}, T_1}[H^{T_1^{j, 1}}]}) \otimes \Big\langle \cI_{t, u}^{i, \gamma, \alpha}[\mu](\omega), T_1^{j, 2} \Big\rangle(\iota_{\psi^{T_1^{j, 2}, T_1}[H^{T_1^{j, 2}}]}) \bigg)
			\\
			&\qquad \otimes
			\bigg( \sum_{j=0}^{n} \Big\langle \cI_{s, t}^{i, \gamma, \alpha}[\mu](\omega), T_2^{j, 1} \Big\rangle(\iota_{\psi^{T_2^{j, 1}, T_2}[H^{T_2^{j, 1}}]}) \otimes \Big\langle \cI_{t, u}^{i, \gamma, \alpha}[\mu](\omega), T_2^{j, 2} \Big\rangle(\iota_{\psi^{T_2^{j, 2}, T_2}[H^{T_2^{j, 2}}]}) \bigg)
			\\
			&= \Big\langle \cI_{s, u}^{i, \gamma, \alpha}[\mu](\omega), T_1 \Big\rangle(\iota_{H^{T_1}}) \otimes \Big\langle \cI_{s, u}^{i, \gamma, \alpha}[\mu](\omega), T_2 \Big\rangle(\iota_{H^{T_2}})
			= 
			\cI_{s, u}^{T_1 \circledast T_2}[\mu](\omega), 
		\end{align*}
		so Equation \eqref{eq:proposition:I_geoPRP-2} holds for $T_1\circledast T_2$. 
		
		Next, suppose that Equation \eqref{eq:proposition:I_geoPRP-2} holds for $T=(\scN, \scE, h_0, H) \in \scF^{\gamma, \alpha, \alpha}$ and suppose additionally that $\lfloor T\rfloor \in \scT^{\gamma, \alpha, \alpha}$. 
		
		Using Sweedler notation, courtesy of Definition \ref{definition:coproduct} we have that
		\begin{align*}
			\Delta\big[T\big] = \sum_{j=0}^n T^{j, 1} \times^{H^T} T^{j, 2} 
			\quad \mbox{and} \quad 
			\Delta\Big[ \lfloor T \rfloor \Big] = \lfloor T\rfloor \times^{H} \rId + \sum_{j=0}^n T^{j, 1} \times^{H} \lfloor T^{j, 2} \rfloor
		\end{align*}
		Then 
		\begin{align*}
			&\Big\langle \cI_{s, t}^{i, \gamma, \alpha}[\mu](\omega) \ast \cI_{t, u}^{i, \gamma, \alpha}[\mu](\omega), \lfloor T \rfloor \Big\rangle(\iota_{H}) 
			\\
			&= 
			\Big\langle \cI_{s, t}^{i, \gamma, \alpha}[\mu](\omega), \lfloor T \rfloor \Big\rangle(\iota_{H}) 
			\\
			&\qquad + \sum_{j=0}^{n} \Big\langle \cI_{s, t}^{i, \gamma, \alpha}[\mu](\omega), T^{j, 1} \Big\rangle(\iota_{\psi^{H^{j, 1}, H}[H^{j, 1}]}) \otimes \Big\langle \cI_{t, u}^{i, \gamma, \alpha}[\mu](\omega), \lfloor T^{j, 2}\rfloor \Big\rangle(\iota_{\psi^{H^{j, 2}, H}[H^{j, 2}]}) 
			\\
			&= 
			\int_s^t \Big\langle \cI_{s, r}^{i, \gamma, \alpha}[\mu](\omega), T \Big\rangle(\iota_{H}) \otimes dW_r^{i, N}
			\\
			&\qquad + 
			\int_{t}^u \sum_{j=0}^{n} \Big\langle \cI_{s, t}^{i, \gamma, \alpha}[\mu](\omega), T^{j, 1} \Big\rangle(\iota_{\psi^{H^{j,1}, H}[H^{j,1}]}) \otimes \Big\langle \cI_{t, r}^{i, \gamma, \alpha}[\mu](\omega), T^{j, 2} \Big\rangle(\iota_{\psi^{H^{j,2}, H}[H^{j,2}]}) \otimes dW_r^{i, N}
			\\
			&=
			\int_s^t \Big\langle \cI_{s, r}^{i, \gamma, \alpha}[\mu](\omega), T \Big\rangle(\iota_{H}) \otimes dW_r^{i, N}
			+ 
			\int_t^u \Big\langle \cI_{s, r}^{i, \gamma, \alpha}[\mu](\omega), T \Big\rangle(\iota_{H}) \otimes dW_r^{i, N} 
			\\
			&=
			\int_s^u \Big\langle \cI_{s, r}^{i, \gamma, \alpha}[\mu](\omega), T \Big\rangle(\iota_{H}) \otimes dW_r^{i, N}
			= 
			\cI_{s, u}^{\lfloor T\rfloor}[\mu](\omega)(\iota_{H}), 
		\end{align*}
		so Equation \eqref{eq:proposition:I_geoPRP-2} holds for $\lfloor T \rfloor$. 
		
		Finally, we recall Proposition \ref{proposition:CompletionOfTrees} to conclude via induction that Equation \eqref{eq:proposition:I_geoPRP-2} holds for all $T\in \scF^{\gamma, \alpha, \alpha}$. 
	\end{proof}
	
	A natural question to address at this point is how are probabilistic rough paths connected to branched rough paths as first described in \cite{gubinelli2010ramification}. Recall from Subsection \ref{subsection:2.1LionsTrees} that we use the notation $T\in \scF$ to denote a Lions forest and $\tau \in \fF_{N}$ to denote a labelled directed forest. 
	\begin{example}
		\label{example:Branched=Probab}
		Let $N\in \bN$ and let $i \in \{1, ..., N\}$. Let $\alpha>0$ and let $\gamma = \sup\{i \alpha: i \in \bN, i\alpha \leq 1\}$. Suppose that $N> \tfrac{\gamma}{\alpha}$. Let
		\begin{equation*}
			\cH_{d, N}^{\gamma, \alpha} = \bigoplus_{\tau \in \fF_{0, N }^{\gamma, \alpha}} \lin\big( (\bR^d)^{\otimes |\scN^T|}, \bR^e \big) 
		\end{equation*}
		and let $\cC\big( \cH_{d, N}^{\gamma, \alpha}, \alpha \big)$ be the set of all $(\cH_{d, N}^{\gamma, \alpha}, \alpha)$-(standard) rough paths, i.e., the set of paths $\rv:[0,1] \to G\big( \cH_{d, N}^{\gamma, \alpha} \big)$ that satisfy Equations \eqref{eq:BranchedRP1} and \eqref{eq:BranchedRP2}. 
		
		Let $(\Omega, \cF, \bP)$ be a probability space and let $i\in \{1, ..., N\}$. We define 
		\begin{equation*}
			\fL_i: \bigsqcup_{T\in \scF_0}^{\gamma, \alpha, \alpha} \Big( \Omega^{\times |H^T|} \times \scN^T \Big) \to \{1, ..., N\}
		\end{equation*}
		to be the mapping such that $\fL_i[T](\omega_{H^T}):\scN^T \to \{1, ..., N\}$ satisfies
		\begin{equation*}
			\fL_i\big[ T \big](\omega_{H^T})[x] = 
			\begin{cases}
				i \quad&\quad \mbox{if $x\in h_0^T$}
				\\
				u(\omega_h) \quad&\quad \mbox{if $x\in h$},
			\end{cases}
		\end{equation*}
		where $u:\Omega \to \{1, ..., N\}$ is uniformly distributed. 
		
		Next, we define
		\begin{equation*}
			\fR_i: \bigsqcup_{T\in \scF_0}^{\gamma, \alpha, \alpha} \Omega^{\times |H^T|} \to \fF_{0, N}
		\end{equation*}
		to be the mapping
		\begin{equation*}
			\fR_i\big[ T, \omega_{H^T} \big] = \fR_i\Big[ \big( \scN^T, \scE^T, h_0^T, H^T \big), \omega_{H^T} \Big] = \big( \scN^T, \scE^T, \fL_i[T](\omega_{H^T}) \big). 
		\end{equation*}
		Notice that the mapping $\fR_i[T, \cdot]$ is $\sigma(u)^{\otimes |H^T|}$-measurable. 
		
		For each $s, t\in[0,1]$ such that $s\leq t$ and $T\in \scF^{\gamma, \alpha, \alpha}$, we define
		\begin{equation}
			\label{eq:example:Branched=Probab}
			\Big\langle \rw_{s, t}, T\Big\rangle(i, \omega'_{H^T}) = \Big\langle \rv_{s, t}, \fR_i\big[ T, \omega'_{H^T} \big] \Big\rangle
		\end{equation}
		and for each $i \in \{1, ..., N\}$ we define
		\begin{equation}
			\label{eq:example:Branched=Probab-}
			\rw_{s, t}(i) = \sum_{T\in \scF_0}^{\gamma, \alpha, \alpha} \Big\langle \rw_{s, t}, T \Big\rangle(i, \cdot). 
		\end{equation}
		By the measurability of Equation \eqref{eq:example:Branched=Probab}, we conclude that for any integrability fucntional $\rQ$ and for each $s, t\in [0,1]$
		\begin{equation*}
			\rw_{s, t} \in L^0\bigg( \{1, ..., N\}, \bP\circ (u)^{-1}; \bigoplus_{T\in \scF_0}^{\gamma, \alpha, \alpha} L^{q[T]}\Big( \Omega^{\times |H^T|}, \bP^{\times |H^T|}; (\bR^d)^{\otimes |\scN^T|} \Big) \bigg). 
		\end{equation*}
	\end{example}
	
	Recalling Lemma \ref{lemma:Link-infin_Syst-Part} from earlier, we have the following statement, which asserts that the above construction provides a probabilistic rough path:
	\begin{theorem}
		\label{theorem:Existence-EmpProbRP}
		Let $\alpha>0$ and let $\gamma = \sup\{ i\cdot \alpha: i\in \bN, i\cdot \alpha \leq 1\}$. Let $N\in \bN$ and suppose that $N>\tfrac{\gamma}{\alpha}$. 
		
		Let $(\Omega, \cF, \bP)$ be a probability space carrying a random variable $u: \Omega \to \{1, ..., N\}$ such that $(\bP) \circ u^{-1}$ is uniformly distributed. Let $(\rp, \rQ)$ be a dual integrability functional.
		
		Then for each $s, t\in [0,1]$, the random variable 
		\begin{equation}
			\label{eq:theorem:Existence-EmpProbRP}
			\rw_{s, t}( u ) \in G_{\rQ}^{\gamma, \alpha, \alpha}\big( L^0(\Omega, \bP; \bR^e) \big). 
		\end{equation}
	
		Further, we have that for each $s, t\in [0,1]$, $\bP$-almost surely the collection of random variables $\big( \rw_{s, t}(u) \big)_{s, t \in [0,1]}$ defined in Equation \eqref{eq:example:Branched=Probab-} satisfies Equation \eqref{eq:definition:ProbabilisticRoughPaths1} and \eqref{eq:definition:ProbabilisticRoughPaths2} so that $\rw(u)$ is a probabilistic rough path in the sense of Definition \ref{definition:ProbabilisticRoughPaths}. 
	\end{theorem}
	
	\begin{proof}
		Firstly, for $s, t\in [0,1]$ and any two Lions forests $T_1, T_2 \in \scF_0$ such that $T_1\circledast T_2 \in \scF^{\gamma, \alpha, \alpha}$ we have that 
		\begin{align}
			\nonumber
			\Big\langle \rw_{s, t}(u), T_1 \circledast T_2 \Big\rangle(\omega_0, \omega_{H^{T_1}}, \omega_{H^{T_2}}) =& \Big\langle \rv_{s, t}, \fR_{u(\omega_0)}\big[ T_1 \circledast T_2, \omega_{H^{T_1\circledast T_2}} \big] \Big\rangle
			\\ 
			\nonumber
			=& \Big\langle \rv_{s, t}, \fR_{u(\omega_0)}\big[ T_1 , \omega_{H^{T_1}} \big] \Big\rangle \otimes \Big\langle \rv_{s, t}, \fR_{u(\omega_0)}\big[ T_2 , \omega_{H^{T_2}} \big] \Big\rangle 
			\\
			\label{eq:theorem:Existence-EmpProbRP1.1}
			=&\Big\langle \rw_{s, t}(u), T_1 \Big\rangle(\omega_0, \omega_{H^{T_1}}) \otimes \Big\langle \rw_{s, t}(u), T_2 \Big\rangle(\omega_0, \omega_{H^{T_2}})
		\end{align}
		so that Equation $\rw_{s, t}(u)$ is a character. Since the measures have finite support, for any choice of integrability functional $\rQ$, we have that
		\begin{equation*}
			\Big\langle \rw_{s, t}(u), T \Big\rangle \in L^0\bigg( \Omega, \bP; \bigoplus_{T\in \scF_0}^{\gamma, \alpha, \alpha} L^{\rQ[T]}\Big( \Omega^{\times |H^T|}, \bP^{\times |H^T|}; (\bR^d)^{\otimes |\scN^T|} \Big) \bigg). 
		\end{equation*}
		
		Secondly, thanks to the definition of $\fR_u$, for any $T\in \scF_0^{\gamma, \alpha, \alpha}$ and any $s, t\in[0,1]$ such that $h_0^T \neq \emptyset$, we have that
		\begin{align}
			\nonumber
			&\Big\langle \rw_{s, t}(u), \cE[T] \Big\rangle(\omega_0, \omega_{h_0^T}, \omega_{H^T}) 
			=
			\Big\langle \rv_{s, t}, \fR_{u(\omega_0)}\Big[ \cE[T], \omega_{h_0^T}, \omega_{H^T} \Big] \Big\rangle 
			\\
			\label{eq:theorem:Existence-EmpProbRP1.2}
			&=\Big\langle \rv_{s,t}, \fR_{u(\omega_{h_0^T})}\Big[ T, \omega_{H^T} \Big] \Big\rangle
			=
			\Big\langle \rw_{s,t}(u), T\Big\rangle(\omega_{h_0^T}, \omega_{H^T}) \quad \bP \times (\bP)^{\times |H^{\cE[T]}|}\mbox{-almost surely} 
		\end{align}
		so that Equation \eqref{eq:definition:DualModule2} is satisfied. Hence, for each $s, t\in[0,1]$ we have that Equation \eqref{eq:theorem:Existence-EmpProbRP} is satisfied. 
		
		Next, for any $T\in \scF_0^{\gamma, \alpha, \alpha}$ and $r, s, t\in [0,1]$ such that $r \leq s \leq t$ we have that 
		\begin{align*}
			\Big\langle \rw_{r, s}(u)& \ast \rw_{s, t}(u), T \Big\rangle(\omega_0, \omega'_{H^T}) 
			=
			\Big\langle \rw_{r, s}(u) \tilde{\otimes} \rw_{s, t}(u), \Delta\big[ \cE[T] \big] \Big\rangle(\omega_0, \omega'_{H^T})
			\\
			=&\sum_{\Upsilon, Y \in \scF_0} c\Big( T, \Upsilon, Y \Big) \cdot \Big\langle \rw_{r, s}(u), \Upsilon \Big\rangle \otimes^{H^T} \Big\langle \rw_{s, t}(u), Y \Big\rangle(\omega_0, \omega'_{H^T} ) 
			\\
			=&\sum_{\Upsilon, Y\in \scF_0} c\Big( T, \Upsilon, Y \Big) \cdot \Big\langle \rv_{r, s}, \fR_{u(\omega_0)}\big[ \Upsilon, \omega'_{\psi^{\Upsilon, T}[H^{\Upsilon}]} \big] \Big\rangle \otimes \Big\langle \rv_{s, t}, \fR_{u(\omega_0)}\big[ Y, \omega'_{\psi^{Y, T}[H^{Y}]} \big] \Big\rangle
			\\
			=& \Big\langle \rv_{r, s} \otimes \rv_{s, t}, \triangle\Big[ \fR_{u(\omega_0)}\big[ T, \omega'_{H^T} \big] \Big] \Big\rangle
			\\
			=& \Big\langle \rv_{r, s} \ast \rv_{s, t}, \fR_{u(\omega_0)}\big[ T, \omega_{H^T} \big] \Big\rangle
			=
			\Big\langle \rv_{r, t}, \fR_{u(\omega_0)}\big[ T, \omega_{H^T} \big] \Big\rangle
			=
			\Big\langle \rw_{r, t}(u), T\Big\rangle( \omega_0, \omega_{H^T})
		\end{align*}
		so that $\rw_{s, t}$ satisfies Equation \eqref{eq:definition:ProbabilisticRoughPaths1}. 
	\end{proof}
	
	\subsubsection{Probabilistic rough paths in continuum}
	\label{subsubsection:GeometricPRPs}
	
	In order to provide a trivial example of a probabilistic rough path (in the setting where the cloud of particles is a continuum) that satisfies Definition \ref{definition:ProbabilisticRoughPaths}, we start by developing the concept of \emph{geometric rough path} from classical rough path theory (see for example \cite{frizhairer2014}). 
	\begin{example}
		\label{example:smoothpaths-trees}
		Let $C^1([0,1]; \bR^d)$ be the space of paths $x:[0,1] \to \bR^d$ that are continuously differentiable with the norm
		\begin{equation*}
			\| x\|_1 = |x_0| + \sup_{s \in [0,1]} |x'_s|. 
		\end{equation*}
		Let $\cB$ be the Borel $\sigma$-field. Let $q\geq 2$ and let $\mu$ be a measure defined on $\Big( C^1([0,1]; \bR^d), \cB \Big)$ such that
		\begin{equation*}
			\int_{C^1([0,1]; \bR^d)} \| x\|_1^q d\mu(x)< \infty. 
		\end{equation*}
		Let $\Omega \subseteq C^1([0,1]; \bR^d)$ be the support of the measure $\mu$ (the smallest closed set with full measure), denoted $\Omega=\supp(\mu)$ and let $\cF$ be the restriction of $\cB$ to $\Omega$. Since $C^1([0,1]; \bR^d)$ is Polish, $\mu(\Omega)=1$. 
		
		Then for all $T \in \scF$ and $s, t\in [0,1]$, we define the operators
		\begin{equation*}
			\cJ_{s, t}^T: \Omega \times \Omega^{\times |H^T|} \to (\bR^d)^{\otimes |\scN^T|}
		\end{equation*}
		inductively using the relations
		\begin{equation}
			\label{eq:example:smoothpaths-trees}
			\begin{split}
				&\cJ_{s, t}^{\lfloor \rId \rfloor}(x) = x_{s, t}, 
				\quad
				\cJ_{s, t}^{T_1 \circledast T_2}(x)(y_{H^{T_1}}, y_{H^{T_2}}) = \cJ_{s, t}^{T_1}(x)(y_{H^{T_1}}) \otimes \cJ_{s, t}^{T_2}(x)(y_{H^{T_2}}), 
				\\
				&\cJ_{s,t}^{\cE[T]}(x)(y_{H^{\cE[T]}}) = \cJ_{s,t}^{T}(y_{h_0^T})(y_{H^{T}}), 
				\quad
				\cJ_{s, t}^{\lfloor T \rfloor}(x)(y_{H^T}) = \int_s^t \cJ_{s, r}^{T}(x)(y_{H^T}) \otimes dx_r
			\end{split}
		\end{equation}
		with the usual convention that $y_H = (\underbrace{y_h, ...}_{h\in H})$. Then $\forall T\in \scF_0$ and $\forall s, t\in[0,1]$ the mapping
		\begin{equation}
			\label{eq:example:smoothpaths-trees_map}
			\big( x, y_{H^T} \big) \mapsto \cJ_{s, t}^T(x)(y_{H^T})
		\end{equation}
		is measurable so that
		\begin{equation*}
			\cJ_{s, t}^T \in L^0 \bigg( \Omega, \mu; L^0\Big( \Omega^{\times|H^T|}, \mu^{\times |H^T|};  (\bR^d)^{\otimes |\scN^T|} \Big) \bigg). 
		\end{equation*}
		
		Let $\alpha, \beta>0$, $\gamma>\alpha\wedge \beta$ and suppose that, for all $T \in \scF^{\gamma, \alpha, \beta}$, $\tfrac{q}{|\scN^T|}>1$. Let $\rQ$ be the integrability functional that satisfies Equation \eqref{eq:trivialIntFunct}. Then for each $T\in \scF^{\gamma, \alpha, \beta}$ we have that 
		\begin{equation*}
			\cJ_{s, t}^T \in L^0 \bigg( \Omega, \mu; L^{q[T]}\Big( \Omega^{\times|H^T|}, \mu^{\times |H^T|};  (\bR^d)^{\otimes |\scN^T|} \Big) \bigg). 
		\end{equation*}
		and for each $x\in \Omega$ we define 
		\begin{align*}
			\cJ_{s, t}^{\gamma, \alpha, \beta}[\mu](x) 
			= 
			\sum_{T \in \scF_0}^{\gamma, \alpha, \beta} \cJ_{s, t}^T(x)(\cdot) 
			\in 
			\bigoplus_{T\in \scF_0}^{\gamma, \alpha, \beta} L^{q[T]}\Big( \Omega^{\times |H^T|}, \mu^{\times |H^T|}; (\bR^d)^{\otimes |\scN^T|} \Big) 
		\end{align*}
		so that
		\begin{align*}
			\cJ_{s, t}^{\gamma, \alpha, \beta}[\mu]
			\in 
			L^0\bigg( \Omega, \mu; \bigoplus_{T\in \scF_0}^{\gamma, \alpha, \beta} L^{q[T]}\Big( \Omega^{\times |H^T|}, \mu^{\times |H^T|}; (\bR^d)^{\otimes |\scN^T|} \Big) \bigg). 
		\end{align*}
	\end{example}

	\begin{proposition}
		\label{proposition:I_geoPRP}
		Let $\alpha, \beta>0$ and let $\gamma$ satisfy Equation \eqref{definition:General-PRP-gamma}. Let $q' \geq \big\lfloor \tfrac{\gamma}{\alpha \wedge \beta}\big\rfloor>1$. Let $(\rp, \rQ)$ be the dual integrability functional that satisfies Equation \eqref{eq:trivialIntFunct}. 
		
		Let $s, t, u \in[0,1]$. Suppose that $\mu \in \cP_{q'}\big( C^1([0,1]; \bR^d)\big)$ and let $x \in \supp(\mu)$. Then the increments $\cJ_{\cdot, \cdot}^{\gamma, \alpha, \beta}$ as defined in Equation \eqref{eq:example:smoothpaths-trees} are elements of the group $G_{\rQ}^{\gamma, \alpha, \beta}\big( L^0(\Omega, \mu; \bR^e) \big)$ and satisfies the identity
		\begin{equation*}
			\cJ_{s, t}^{\gamma, \alpha, \beta}[\mu](x) \ast \cJ_{t, u}^{\gamma, \alpha, \beta}[\mu](x) = \cJ_{s, u}^{\gamma, \alpha, \beta}[\mu](x). 
		\end{equation*}
		As such, the path $\cJ^{\gamma, \alpha, \beta}$ introduced in Example \ref{example:smoothpaths-trees} is a McKean-Vlasov probabilistic rough path in the sense of Definition \ref{definition:ProbabilisticRoughPaths}. 
	\end{proposition}
	
	\begin{proof}
		The proof of Proposition \ref{proposition:I_geoPRP} proceeds in a similar fashion to that of the proof of Proposition \ref{proposition:I_geoPRP2}: Firstly, courtesy of Equation \eqref{eq:example:smoothpaths-trees} it is simple enough to verify that for any $t\in [0,1]$, 
		\begin{equation*}
			\cJ_{0, t}^{\gamma, \alpha, \beta}[\mu] \in G_{\rQ}^{\gamma, \alpha, \beta}\big( L^0(\Omega, \bP; \bR^e) \big). 
		\end{equation*}
		Therefore, our main goal is to prove via induction that $\forall T \in \scF^{\gamma, \alpha, \beta}$ that
		\begin{equation}
			\label{eq:proposition:I_geoPRP}
			\Big\langle \cJ_{s, t}^{\gamma, \alpha, \beta}[\mu](x) \ast \cJ_{t, u}^{\gamma, \alpha, \beta}[\mu](x) , T \Big\rangle = \Big\langle \cJ_{s, u}^{\gamma, \alpha, \beta}[\mu](x), T \Big\rangle. 
		\end{equation}
		Firstly, we can have that
		\begin{align*}
			&\Big\langle \cJ_{s, t}^{\gamma, \alpha, \beta}[\mu](x) \ast \cJ_{t, u}^{\gamma, \alpha, \beta}[\mu](x), \lfloor \rId\rfloor \Big\rangle 
			\\
			&= 
			\Big\langle \cJ_{s, t}^{\gamma, \alpha, \beta}[\mu](x), \lfloor \rId\rfloor \Big\rangle \otimes \Big\langle \cJ_{t, u}^{\gamma, \alpha, \beta}[\mu](x), \rId \Big\rangle + \Big\langle \cJ_{s, t}^{\gamma, \alpha, \beta}[\mu](x), \rId\Big\rangle \otimes \Big\langle \cJ_{t, u}^{\gamma, \alpha, \beta}[\mu](x), \lfloor \rId\rfloor \Big\rangle
			\\
			&=\Big\langle \cJ_{s, t}^{\gamma, \alpha, \beta}[\mu](x), \lfloor \rId\rfloor \Big\rangle + \Big\langle \cJ_{t, u}^{\gamma, \alpha, \beta}[\mu](x), \lfloor \rId\rfloor \Big\rangle
			\\
			&=x_{s, t} + x_{t, u} 
			= x_{s, u} = \Big\langle \cJ_{s, u}^{\gamma, \alpha, \beta}[\mu](x), \lfloor \rId\rfloor \Big\rangle. 
		\end{align*}
		
		Next, suppose that Equation \eqref{eq:proposition:I_geoPRP} holds for $T\in \scF^{\gamma, \alpha, \beta}$ and that $\cE[T] \in \scF^{\gamma, \alpha, \beta}$ too. For $T=(\scN^T, \scE^T, h_0^T, H^T, \scL^T)$ we use the Sweedler notation
		\begin{equation}
			\label{eq:SweedlerNotation}
			\Delta\big[ T \big] = \sum_{i=0}^n T^{i, 1} \times^{H^T} T^{i, 2}
		\end{equation}
		so that by Equation \eqref{eq:definition:coproduct}, 
		\begin{equation*}
			\Delta\Big[ \cE\big[T\big] \Big] = \sum_{i=0}^n \cE[T^{i, 1}] \times^{(H^T)'} \cE[T^{i, 2}]. 
		\end{equation*}
		Then
		\begin{align*}
			\Big\langle \cJ_{s, t}^{\gamma, \alpha, \beta}&[\mu](x) \ast \cJ_{t, u}^{\gamma, \alpha, \beta}[\mu](x), \cE[T] \Big\rangle(y_{h_0}, y_{H^T})
			\\
			&=\sum_{i=0}^n \Big\langle \cJ_{s, t}^{\gamma, \alpha, \beta}[\mu](x), \cE\big[ T^{i,1} \big] \Big\rangle \otimes^{(H^T)'} \Big\langle \cJ_{t, u}^{\gamma, \alpha, \beta}[\mu](x) , \cE\big[ T^{i,2} \big] \Big\rangle(y_{(H^T)'})
			\\
			&=\sum_{i=0}^n \cJ_{s, t}^{T^{i,1}}(y_{h_0^{T}})(y_{\psi^{T^{i,1}, T}[H^{T^{i,1}}]}) \otimes \cJ_{t, u}^{T^{i,2}}(y_{h_0^{T}})(y_{\psi^{T^{i,2}, T}[H^{T^{i,2}}]}) 
			\\
			&=
			\cJ_{s, u}^T(y_{h_0^T})(y_{H^T}) = \cJ_{s, u}^{\cE[T]}(x)(y_{H^{\cE[T]}}), 
		\end{align*}
		so Equation \eqref{eq:proposition:I_geoPRP} holds for $\cE[T]$. 
		
		Next, suppose that Equation \eqref{eq:proposition:I_geoPRP} holds for $T_1, T_2 \in \scF^{\gamma, \alpha, \beta}$ and $T_1\circledast T_2\in \scF^{\gamma, \alpha, \beta}$. Using Sweedler notation again, thanks to Equation \eqref{eq:definition:coproduct} we get that
		\begin{align*}
			&\Delta\big[ T_1 \big] = \sum_{i=0}^m T_1^{i, 1} \times^{H^{T_1}} T_1^{i, 2}, 
			\quad
			\Delta\big[ T_2 \big] = \sum_{j=0}^n T_2^{j, 1} \times^{H^{T_2}} T_2^{j, 2} \quad \mbox{and}
			\\
			&\Delta\Big[ T_1 \circledast T_2 \Big] = \sum_{i=0}^m \sum_{j=0}^n \big( T_1^{i, 1} \circledast T_2^{j, 1} \big) \times^{H^{T_1} \cup H^{T_2}} \big( T_1^{i, 2} \circledast T_2^{j, 2} \big). 
		\end{align*}
		Then
		\begin{align*}
			&\Big\langle \cJ_{s, t}^{\gamma, \alpha, \beta}[\mu](x) \ast \cJ_{t, u}^{\gamma, \alpha, \beta}[\mu](x), T_1 \circledast T_2 \Big\rangle (y_{H^{T_1 \circledast T_2}})
			\\
			&=\sum_{i=0}^{m} \sum_{j=0}^{n} \Big\langle \cJ_{s, t}^{\gamma, \alpha, \beta}[\mu](x), T_1^{i, 1} \circledast T_2^{j, 1} \Big\rangle \otimes^{H^{T_1} \cup H^{T_2}} \Big\langle \cJ_{t, u}^{\gamma, \alpha, \beta}[\mu](x), T_1^{i,2} \circledast T^{j, 2} \Big\rangle(y_{H^{T_1}}, y_{H^{T_2}})
			\\
			&=\sum_{i=0}^{m} \sum_{j=0}^{n} \bigg( \Big\langle \cJ_{s, t}^{\gamma, \alpha, \beta}[\mu](x), T_1^{i, 1} \Big\rangle \otimes \Big\langle \cJ_{t, u}^{\gamma, \alpha, \beta}[\mu](x), T_2^{j, 1} \Big\rangle \bigg)
			\\
			&\qquad \otimes^{H^{T_1}\cup H^{T_2}}
			\bigg( \Big\langle \cJ_{s, t}^{\gamma, \alpha, \beta}[\mu](x), T_1^{i, 2} \Big\rangle \otimes \Big\langle \cJ_{t, u}^{\gamma, \alpha, \beta}[\mu](x), T_2^{j, 2} \Big\rangle \bigg)(y_{H^{T_1}}, y_{H^{T_2}})
			\\
			&=\bigg( \sum_{i=0}^{m} \Big\langle \cJ_{s, t}^{\gamma, \alpha, \beta}[\mu](x), T_1^{i, 1} \Big\rangle \otimes^{H^{T_1}} \Big\langle \cJ_{t, u}^{\gamma, \alpha, \beta}[\mu](x), T_1^{i, 2} \Big\rangle(y_{H^{T_1}}) \bigg)
			\\
			&\qquad \otimes
			\bigg( \sum_{j=0}^{n} \Big\langle \cJ_{s, t}^{\gamma, \alpha, \beta}[\mu](x), T_2^{j, 1} \Big\rangle \otimes^{H^{T_2}} \Big\langle \cJ_{t, u}^{\gamma, \alpha, \beta}[\mu](x), T_2^{j, 2} \Big\rangle(y_{H^{T_2}}) \bigg)
			\\
			&= \Big\langle \cJ_{s, u}^{\gamma, \alpha, \beta}[\mu](x), T_1 \Big\rangle(y_{H^{T_1}}) \otimes \Big\langle \cJ_{s, u}^{\gamma, \alpha, \beta}[\mu](x), T_2 \Big\rangle(y_{H^{T_2}}) 
			= 
			\cJ_{s, u}^{T_1 \circledast T_2}[\mu](x)(y_{H^{T_1 \circledast T_2}}), 
		\end{align*}
		so Equation \eqref{eq:proposition:I_geoPRP} holds for $T_1\circledast  T_2$. 
		
		Finally, suppose that Equation \eqref{eq:proposition:I_geoPRP} for $T \in \scF^{\gamma, \alpha, \beta}$. Using Sweedler notation, thanks to Equation \eqref{eq:definition:coproduct} we have that
		\begin{align*}
			\Delta\big[T\big] = \sum_{j=0}^n T^{j, 1} \times^{H^T} T^{j, 2} 
			\quad \mbox{and} \quad 
			\Delta\Big[ \lfloor T \rfloor \Big] = \lfloor T\rfloor \times^{H^T} \rId + \sum_{j=0}^n T^{j, 1} \times^{H^T} \lfloor T^{j, 2} \rfloor
		\end{align*}
		
		Then
		\begin{align*}
			&\Big\langle \cJ_{s, t}^{\gamma, \alpha, \beta}[\mu](x) \ast \cJ_{t, u}^{\gamma, \alpha, \beta}[\mu](x), \lfloor T \rfloor \Big\rangle(y_{H^T}) 
			\\
			&= 
			\Big\langle \cJ_{s, t}^{\gamma, \alpha, \beta}[\mu](x), \lfloor T \rfloor \Big\rangle(y_{H^T}) + \sum_{j=0}^{n} \Big\langle \cJ_{s, t}^{\gamma, \alpha, \beta}[\mu](x), T^{j, 1} \Big\rangle \otimes^{H^T} \Big\langle \cJ_{t, u}^{\gamma, \alpha, \beta}[\mu](x), \lfloor T^{j, 2}\rfloor \Big\rangle(y_{H^T}) 
			\\
			&= 
			\int_s^t \Big\langle \cJ_{s, r}^{\gamma, \alpha, \beta}[\mu](x), T \Big\rangle(y_{H^T}) \otimes dx_r
			\\
			&\qquad + 
			\int_{t}^u \sum_{j=0}^{n} \Big\langle \cJ_{s, t}^{\gamma, \alpha, \beta}[\mu](x), T^{j, 1} \Big\rangle \otimes^{H^T} \Big\langle \cJ_{t, r}^{\gamma, \alpha, \beta}[\mu](x), T^{j, 2} \Big\rangle(y_{H^T}) \otimes dx_r
			\\
			&=
			\int_s^t \Big\langle \cJ_{s, r}^{\gamma, \alpha, \beta}[\mu](x), T \Big\rangle(y_{H^T}) \otimes dx_r
			+ 
			\int_t^u \Big\langle \cJ_{s, r}^{\gamma, \alpha, \beta}[\mu](x), T \Big\rangle(y_{H^T}) \otimes dx_r 
			\\
			&=
			\int_s^u \Big\langle \cJ_{s, r}^{\gamma, \alpha, \beta}[\mu](x), T \Big\rangle(y_{H^T}) \otimes dx_r
			= 
			\cJ_{s, u}^{\lfloor T\rfloor}[\mu](x)(y_{H^T}), 
		\end{align*}
		so Equation \eqref{eq:proposition:I_geoPRP} holds for $\lfloor T \rfloor$. 
		
		Finally, we recall Proposition \ref{proposition:CompletionOfTrees} to conclude via induction that Equation \eqref{eq:proposition:I_geoPRP} holds for all $T\in \scF^{\gamma, \alpha, \beta}$.
	\end{proof}
	
	\subsection{Strong probabilistic rough paths}
	\label{subsection:StrongPRP}
	
	The next step in the analysis of probabilistic rough paths is to extend the notion of \emph{strong rough path} introduced in \cite{2019arXiv180205882.2B}. The idea is to require all the components of a probabilistic rough path (corresponding to iterated integrals) to be constructed as measurable functions of the underlying (first level) trajectory. This property is a natural pre-requisite for addressing the question of a mean-field limit of empirical probabilistic rough paths converging to McKean-Vlasov probabilistic rough paths. 
	
	In order to construct a strong probabilistic rough path of any kind, we need to show the existence of a measurable mapping of an $\alpha$-H\"older continuous path that satisfies the following:
	\begin{definition}
		\label{dfn:Enhancement-BRP}
		Let $0<\alpha<1$ and let $\gamma = \alpha \cdot \big\lfloor \tfrac{1}{\alpha} \big\rfloor$. Let $\Omega \subseteq C^{\alpha}([0,1]; \bR^d)$ be a complete and separable metric space with respect to the H\"older metric. 
		
		Let $T=(\scN, \scE, h_0, H)\in \scF^{\gamma, \alpha, \alpha}$ and recall that we denote 
		$H'= \big( H \cup \{h_0\} \big)\backslash\{\emptyset\}$. For each $T\in \scF^{\gamma, \alpha, \alpha}$ and $s, t \in [0,1]$,  we define the collection of operators 
		\begin{equation*}
			\cM_{s, t}^T: \Omega^{\times |H'|} \to (\bR^d)^{\otimes |\scN|}
		\end{equation*}
		to be Borel measurable operators such that:
		\begin{enumerate}
			\item For each $T=(\scN, \scE, h_0, H) \in \scF^{\gamma, \alpha, \alpha}$ and for $s, t\in[0,1]$, we have the inductive relationship
			\begin{equation}
				\label{eq:dfn:Enhancement-BRP1}
				\left.
				\begin{aligned}
					&\cM_{s, t}^{\lfloor \rId\rfloor }\big( w \big) = w_{s, t}
					\\
					&\cM_{s, t}^{\cE[T]} \big( w_{(H^{\cE[T]})'} \big) = \cM_{s, t}^{T} \big( w_{(H^T)'} \big)
					\\ 
					&\cM_{s, t}^{T_1 \circledast T_2}\big( w_{(H^{T_1\circledast T_2})'} \big) = \cM_{s, t}^{T_1} \otimes^{G'} \cM_{s, t}^{T_2} \big( w_{G'} \big)
				\end{aligned}
				\right\}
			\end{equation}
			where we used that
			\begin{equation}
				\label{eq:dfn:Enhancement-BRP3}
				G' = \big( H^{T_1} \cup H^{T_2} \cup \{h_0^{T_1} \cup h_0^{T_2} \} \big)\backslash \{\emptyset \} \in \lion\big( (H^{T_1} )', (H^{T_2})' \big)
			\end{equation} 
			\item For each $T=(\scN, \scE, h_0, H) \in \scF^{\gamma, \alpha, \alpha}$ and $\forall s, t, u \in [0,1]$
			\begin{equation}
				\label{eq:dfn:Enhancement-BRP1*}
				\cM_{s, u}^{T}\big( w_{H'} \big) = \sum_{\Upsilon, Y \in \scF_0}^{\gamma, \alpha, \alpha} c\Big( T, \Upsilon, Y \Big) \cdot \cM_{s, t}^{\Upsilon} \otimes^{(H)'} \cM_{t, u}^{Y}\big( w_{H'} \big)
			\end{equation}
			where again we used that
			\begin{equation}
				\label{eq:dfn:Enhancement-BRP4}
				H \in \lion\big( H^{\Upsilon} , H^{Y} \big)  \implies H' \in \lion\big( (H^{\Upsilon})' , (H^{Y})' \big).  
			\end{equation}
			\item For each $T = (\scN, \scE, h_0, H) \in \scF_{0}^{\gamma, \alpha, \alpha}$ and $\forall s, t\in[0,1]$, $\exists C_T:\Omega^{\times |H'|} \to \bR^+$, a measurable function of $x_{H'} \in \Omega^{\times|H'|}$ such that
			\begin{equation}
				\label{eq:dfn:Enhancement-BRP2}
				\sup_{s, t \in [0,1]} \frac{ \Big| \cM_{s, t}^T \big( w_{H'} \big) \Big|}{|t-s|^{\alpha |\scN|}} \leq C_{T} \big( w_{H'} \big) . 
			\end{equation}
		\end{enumerate}
	\end{definition}

	Equations \eqref{eq:dfn:Enhancement-BRP3} and \eqref{eq:dfn:Enhancement-BRP4} are easy to verify and capture the coupling (or lack thereof) between the tagged hyperedges $h_0^{T_1}$, $h_0^{T_2}$ that arises via the hyperedge $h_0^T \in (H^T)'$. 
	
	Notice that $\Omega$ in Definition \ref{dfn:Enhancement-BRP} cannot be chosen as the entire space $C^{\alpha}([0,1]; \bR^d)$ since the latter is not separable. A convenient choice is for instance the closure of smooth paths from $[0,1]$ to $\bR^d$ under the $\alpha$-H\"older metric, see for instance \cite{friz2010differential}*{Chapter 5}. 
	
	In particular, we highlight that in Example \ref{example:smoothpaths-trees} we required that the space $\Omega$ was a subset of continuously differentiable paths whereas in Definition \ref{dfn:Enhancement-BRP} the only assumption we make about the regularity is \eqref{eq:dfn:Enhancement-BRP2} and allows for paths that are $\alpha$-H\"older continuous paths. 
	
	\begin{remark}
		The reader may recall certain similarities between Equation \eqref{eq:example:smoothpaths-trees} and Equation \eqref{eq:dfn:Enhancement-BRP1}. While these operations are defined for different collections of inputs for the same tree, we want to emphasise that the additional identity
		\begin{equation*}
			\cJ_{s, t}^{\lfloor T\rfloor} \big( w, w_{H} \big) = \int_s^t \cJ_{s, r}^T \big( w, w_{H} \big) \otimes dw_r
		\end{equation*}
		implies that $\forall s, t, u\in [0,1]$, 
		\begin{equation*}
			\cJ_{s, u}^{T}\big( w_{h_0}, w_{H} \big) = \sum_{\Upsilon, Y \in \scF_0}^{\gamma, \alpha, \alpha} c\Big( T, \Upsilon, Y \Big) \cdot \cJ_{s, t}^{\Upsilon} \otimes^{H} \cJ_{t, u}^{Y}\big( w_{h_0}, w_{H} \big)
		\end{equation*}
		whereas Equation \eqref{eq:dfn:Enhancement-BRP1*} ensures this but does not require any relationship between $\cM_{s, t}^T$ and $\cM_{s, t}^{\lfloor T\rfloor}$. 
	\end{remark}
	
	The following result guarantees that Definition \ref{dfn:Enhancement-BRP} is not empty by proving the existence of a branched rough path with given marginals:
	\begin{theorem}
		For any fixed choice of $\Omega$, there exists at least one operator that satisfies Definition \ref{dfn:Enhancement-BRP}. 
	\end{theorem}
	
	In fact, there are many choices of such a lift from a path to a signature. Observe that when $\Omega \subseteq C^1([0,1]; \bR^d)$, we could choose $\cJ$ as introduced in Example \ref{example:smoothpaths-trees}. However, our focus is on when $\Omega \subseteq C^{\alpha}([0,1]; \bR^d)$ for $\alpha<1$:
	\begin{proof}
		Fix $T=(\scN, \scE, h_0, H) \in \scF^{1, \alpha, \alpha}$. Let $w_{H'} \in \Omega^{\times |H'|}$, so we can equivalently say
		\begin{equation*}
			w_{H'} = \big( w_h \big)_{h\in H'} \in C^{\alpha}\big( [0,1]; (\bR^d)^{\oplus |H'|} \big)
		\end{equation*}
		
		Following the ideas of \cite{tapia2020geometry}*{Theorem 3.4}, we can find an $\alpha$-H\"older continuous path taking values on the Lie group 
		\begin{equation*}
			G\big( \cH_{d\times |H'|}^{1, 1/2-}\big) \subseteq \bigoplus_{\tau \in \fF_{0, |H'|}}^{|\tau|\leq 2} (\bR^d)^{\otimes |\scN^{\tau}|}
		\end{equation*} 
		that agrees with $w_{H'}$ when restricted to the first level of the tensor expansion. Note this is a sum over all the directed rooted forests with labellings taking values on the set $\{1, ..., |H'|\}$. The path is constructed explicitly from $w_{H'}$ so that this $\alpha$-H\"older continuous path is a measurable function of $w_{H'}$. 
		
		By repeating this argument, we can find an $\alpha$-H\"older continuous path $\rw$ taking values on the Lie group $G\big( \cH_{d\times |H'|}^{1, \alpha}\big)$ that agrees with $w_{H'}$ when restricted to the first level of the tensor expansion. Then $\rw$ is a $\big(\cH_{d\times |H'|}^{1, \alpha}, \alpha \big)$-rough path. Further, $\rw$ is a measurable function of the original path $w_{H'}$. We define
		\begin{equation*}
			\cM_{s, t}^{T} \big( w_{H'} \big) = \Big\langle \rw, \tau' \Big\rangle
		\end{equation*}
		where $\tau'=(\scN, \scE, \scL) \in \fF_{0, |H'|}^{\gamma, \alpha}$ and $\scL:\scN\to H'$ such that $\{\scL^{-1}[h]: h\in H'\} = H'$. That is, 
		\begin{equation*}
			\scL[y] = h \quad \mbox{s.t.} \quad y \in h. 
		\end{equation*}
		The key point here is to check that the right-hand side in the above identity only depends on the label $\scL$ of $\tau'$ through the pre-images of $\scL$. This can be verified by following the inductive construction achieved in \cite{tapia2020geometry}. Then by construction, $\cM_{s, t}^{T}$ satisfies the properties \eqref{eq:dfn:Enhancement-BRP1} and $\cM_{s, t}^{T}$ is measurable. 
		
		Secondly, since $\rw$ is a $\big(\cH_{d\times |H'|}^{1, \alpha}, \alpha \big)$-rough path, the regularity condition Equation \eqref{eq:dfn:Enhancement-BRP2} is satisfied. 
	\end{proof}
	
	Finally, we emphasise that the lift constructed in \cite{tapia2020geometry}*{Theorem 3.4} is not unique and specific circumstances may be inappropriate given a choice of application for the associated rough differential equations. 
	
	\subsubsection{Strong probabilistic rough paths for interacting systems of equations}
	
	The proof of the existence of a strong probabilistic rough path is much more straightforward when the underlying measure $\bP$ is finite support:
	\begin{proposition}
		\label{proposition:Enhancement-BRP}
		Let $\alpha>0$, let $\gamma = \alpha \cdot \big\lfloor\tfrac{1}{\alpha}\big\rfloor$ and let $N\in \bN$ such that $N> \big\lfloor \tfrac{1}{\alpha} \big\rfloor$. 
		
		Let
		\begin{equation*}
			\big( \cM_{s, t}^T \big)_{s,t\in[0,1], T\in \scF_0^{\gamma, \alpha, \alpha}}
		\end{equation*}
		be a collection of operators that satisfy Definition \ref{dfn:Enhancement-BRP}. 
		
		Let $(\Omega, \cF, \bP)$ be a probability space equipped with the mapping $u:\Omega \to \{1, ..., N\}$ such that $\bP\circ (u)^{-1}$ is uniformly distributed. Let $(\rp, \rQ)$ be a dual pair of integrability functionals. For each $i\in \{1, ..., N\}$, let
		\begin{equation*}
			w^{i, N} \in C^{\alpha}\big([0,1]; \bR^d \big). 
		\end{equation*}
		
		Then there exists a probabilistic rough path $\rw:[0,1] \to G_{\rQ}^{\gamma, \alpha, \alpha}\big( L^0(\Omega, \bP; \bR^e ) \big)$ such that
		\begin{align*}
			&\Big\langle \rw_{s, t}, \lfloor \rId \rfloor \Big\rangle(\omega_0) = w_{s, t}^{u(\omega_0), N}
			\quad \mbox{and}\quad 
			\Big\langle \rw_{s, t}, \cE\big[ \lfloor \rId \rfloor \big] \Big\rangle (\omega_0, \omega) = w_{s, t}^{u(\omega), N}. 
		\end{align*}
	\end{proposition}
	
	\begin{proof}
		The first step is to apply \cite{tapia2020geometry}*{Theorem 3.4} to any path $\boldsymbol{w} = (w^{i, N})\in C^{\alpha}\big( [0,1]; \bR^e \big)^{\otimes |N|}$ and construct a $\big( \cH_{d\times N}^{\gamma, \alpha}, \alpha \big)$-rough path
		\begin{align*}
			t\mapsto \rv_t = \sum_{\tau \in \fF_{0, N}} \Big\langle \rv_t, \tau \Big\rangle \in \bigoplus_{\tau \in \fF_{0, N}} (\bR^d)^{\otimes |\scN^\tau|}. 
		\end{align*}
		Then we apply Theorem \ref{theorem:Existence-EmpProbRP} to conclude. 
	\end{proof}
	
	\subsubsection{Strong probabilistic rough paths for the continuum of equations}
	
	Now the existence of a collection of operators that satisfies Definition \ref{dfn:Enhancement-BRP} is established, we can use them to construct a strong probabilistic rough path:
	\begin{proposition}
		Let $0<\alpha<1$ and let $\gamma = \alpha \lfloor \tfrac{1}{\alpha}\rfloor$. Let $\Omega\subseteq C^{\alpha}([0,1]; \bR^d)$ be a complete and separable metric space with respect to the H\"older metric. Let $\cF$ be the Borel $\sigma$-algebra of $\Omega$. 
		
		Suppose that 
		\begin{equation*}
			\Big( \cM_{s, t}^T \Big)_{s, t\in [0,1], T\in \scF^{\gamma, \alpha, \alpha}}
		\end{equation*}
		to be a collection of operators that satisfy Definition \ref{dfn:Enhancement-BRP}. For each $T\in \scF^{\gamma, \alpha, \alpha}$ we define
		\begin{equation*}
			\Big\langle \rw_{t}, T \Big\rangle(w_0, w_{H^T}) = \left\{
			\begin{aligned}
				&\cM_{0, t}^T\big( w_0, w_{H^T} \big) 
				\qquad &\mbox{if $h_0^T \neq \emptyset$,}
				\\
				&\cM_{0, t}^T\big( w_{H^T} \big)
				\qquad &\mbox{if $h_0^T = \emptyset$. }
			\end{aligned}
			\right.
		\end{equation*}
		and
		\begin{equation*}
			\rw_t(w_0) = \sum_{T\in \scF_0}^{\gamma, \alpha, \alpha} \Big\langle \rw_{t}, T \Big\rangle(w_0, \cdot) . 
		\end{equation*}
		
		Finally, let $(\rp, \rQ)$ be a dual integrability functional and let $\bP$ be a $(\Omega, \cF)$-probability measure such that for every $T\in \scF^{\gamma, \alpha, \alpha}$ and $\forall t\in [0,1]$, 
		\begin{equation*}
			\big\langle \rw_t, T\big\rangle \in L^0\bigg( \Omega, \bP; L^{q[T]}\Big( \Omega^{\times |H^T|}, \bP^{\times |H^T|}; (\bR^d)^{\otimes |\scN^T|} \Big) \bigg)
		\end{equation*}
		Then $\rw:[0,1] \to G_{\rQ}^{\gamma, \alpha, \alpha}\big( L^0(\Omega, \bP; \bR^e)\big)$ is a strong probabilistic rough path in the sense of Definition \ref{definition:ProbabilisticRoughPaths} and $\forall s, t\in[0,1]$, 
		\begin{equation*}
			\Big\langle \rw_{s, t}, \lfloor \rId\rfloor \Big\rangle(w) = w_{s, t}. 
		\end{equation*}
	\end{proposition}

	\begin{proof}
		The key steps here are that for $\bP$-almost every $w_0 \in \Omega$, the mapping
		\begin{equation*}
			\Big\langle \rw_t, T\Big\rangle(w_0, \cdot) \in L^{q[T]}\Big( \Omega^{\times |H^T|}, \bP^{\times |H^T|}; (\bR^d)^{\otimes |\scN^T|} \Big)
		\end{equation*}
		so that
		\begin{equation*}
			\rw_t(w_0) = \sum_{T\in \scF_0}^{\gamma, \alpha, \alpha} \Big\langle \rw_{t}, T \Big\rangle(w_0, \cdot) \in \bigoplus_{T\in \scF}^{\gamma, \alpha, \alpha} L^{q[T]}\Big( \Omega^{\times |H^T|}, \bP^{\times |H^T|}; (\bR^d)^{\otimes |\scN^T|} \Big)
		\end{equation*}
		and further
		\begin{equation*}
			\rw_t \in L^0\bigg( \Omega, \bP; \bigoplus_{T\in \scF}^{\gamma, \alpha, \alpha} L^{q[T]}\Big( \Omega^{\times |H^T|}, \bP^{\times |H^T|}; (\bR^d)^{\otimes |\scN^T|} \Big) \bigg). 
		\end{equation*}
		Verifying that $\rw$ is a path on group of McKean-Vlasov characters that satisfies Equation \eqref{eq:definition:ProbabilisticRoughPaths1} and Equation \eqref{eq:definition:ProbabilisticRoughPaths2} is then just a task of verification via the properties that arise from $\cM$. 
	\end{proof}
	
	

	
	\appendix
	
	\section{Notation}
	\label{subsection:Notation}
	
	Let $\bN$ be the set of positive integers and 
	$\bN_{0}=\bN \cup \{0\}$. Let $\bR$ be the field of real numbers and for $d \in \bN$, let $\bR^d$ be the $d$-dimensional vector space over the field $\bR$. Let $\langle \cdot,\cdot \rangle$ be the Euclidean inner product over the vector space $\bR^d$. 
	
	For modules $U$ and $V$ over a ring $\cR$, we define $\lin(U, V)$ to be the collection of linear operators from $U$ to $V$ (which is a module over the ring $\cR$). Let $U \oplus V$ and $U \otimes V$ be the direct sum and tensor product of two modules. 

	For a topological module $U$, let $\cB(U)$ be the Borel $\sigma$-algebra. Let $(\Omega, \cF, \bP)$ be a probability space. For $p \in (1, \infty)$, let $L^p(\Omega, \cF, \bP; U)$ be the space of $p$-integrable random variable taking values in $U$. When the $\sigma$-algebra is not ambiguous, we will simply write $L^p(\Omega, \bP; U)$. Further, let $L^0(\Omega, \bP; U)$ be the space of measurable mappings $(\Omega, \cF) \mapsto (U, \cB(U))$ and $L^{\infty}(\Omega, \bP; U)$ the space of essentially bounded mappings $(\Omega, \bF) \mapsto (U, \cB(U))$. 
	
	For a set $\scN$, we call $2^\scN$ the collection of subsets of $\scN$ and $\scP(\scN)$ the set of all partitions of the set $\scN$. This means $\scP(\scN) \subseteq 2^{2^\scN}$. A partition $P \in \scP(\scN)$
	if and only if the following three properties are satisfied:
	\begin{equation*}
		\forall x \in \scN, \quad \exists p \in P: x \in p; 
		\qquad
		\forall p, q \in P, \quad p \cap q = \emptyset; 
		\qquad
		\emptyset \notin P.
	\end{equation*}
	The set of partitions $\scP(\scN)$ has a partial ordering $\subseteq$ where $P \subseteq Q$ if and only if
	\begin{equation*}
		\forall q\in Q, \quad \exists p \in P: q\subseteq p
	\end{equation*}
	or in words $Q$ is finer than $P$. 
\end{document}